\documentclass[11pt,reqno]{amsart}
\pdfoutput=1

\usepackage[utf8]{inputenc}
\usepackage[T1]{fontenc}
\IfFileExists{lmodern.sty}{\usepackage{lmodern}}{}
\usepackage[a4paper,margin=1.15in]{geometry}

\usepackage{amsmath,amssymb,amsfonts,amsthm}
\usepackage{mathtools}
\usepackage{booktabs}
\usepackage{array}
\usepackage{longtable}
\usepackage{enumitem}
\usepackage[final,expansion=false]{microtype}
\usepackage{xurl}
\usepackage{xcolor}
\usepackage[colorlinks=true,
linkcolor={blue!70!black},
citecolor={blue!70!black},
urlcolor={blue!80!black}]{hyperref}

\hypersetup{
	pdftitle={Superspecial curves of genus two over prime fields: Frobenius types and rational (2,2)-kernels},
	pdfauthor={Hung Dang}}

\newtheorem{theorem}{Theorem}[section]
\newtheorem{proposition}[theorem]{Proposition}
\newtheorem{lemma}[theorem]{Lemma}
\newtheorem{corollary}[theorem]{Corollary}
\newtheorem{conjecture}[theorem]{Conjecture}
\theoremstyle{definition}
\newtheorem{definition}[theorem]{Definition}
\theoremstyle{remark}
\newtheorem{remark}[theorem]{Remark}

\DeclareMathOperator{\Aut}{Aut}
\DeclareMathOperator{\End}{End}
\DeclareMathOperator{\Jac}{Jac}
\DeclareMathOperator{\Frob}{Frob}
\DeclareMathOperator{\ord}{ord}
\DeclareMathOperator{\Mass}{Mass}
\DeclareMathOperator{\Hom}{Hom}

\DeclareMathOperator{\Fib}{Fib}
\providecommand{\PGL}{\operatorname{PGL}}
\newcommand{\F}{\mathbb{F}}
\newcommand{\Fp}{\mathbb{F}_p}
\newcommand{\Z}{\mathbb{Z}}
\newcommand{\Q}{\mathbb{Q}}
\newcommand{\PP}{\mathbb{P}}
\newcommand{\kbar}{\bar k}
\newcommand{\eps}{\varepsilon}

\newcommand{\leg}[2]{\bigl(\tfrac{#1}{#2}\bigr)}

\begin{document}

\title[Superspecial curves of genus two over prime fields]{Superspecial curves of genus two over prime fields: Frobenius types and rational \((2,2)\)-kernels}

\author{Hung Dang}
\address{Department of Mathematics, University of Phuong Dong, Hanoi, Vietnam}
\email{hung.dt@phuongdong.edu.vn}
\thanks{ORCID 0009-0006-3272-0573.}


\begin{abstract}
	Superspecial curves of genus two over a prime field \(\Fp\), and the rational
	\((2,2)\)-isogenies between their Jacobians, are the vertices and edges of the
	graph on which genus-two isogeny-based cryptography is built. Let \(C\) be such
	a curve, with Jacobian \(J\) and Frobenius endomorphism \(\pi\). We show that
	\(\pi^{2}=p\,\alpha_*\) for a unique automorphism \(\alpha\) of \(C\) over
	\(\Fp\) whose action on the Weierstrass points is the square of that of
	Frobenius. For \(p\ge7\) this gives the list of Frobenius types of the
	Weierstrass points: eight types, each with its Weil polynomial, a congruence
	condition on \(p\), and its number \(15,3,3,7,1,1,3\) or \(1\) of
	\(\Fp\)-rational \((2,2)\)-kernels; when all Weierstrass points are rational,
	the \(2\)-Tate pairing on \(J[2]\) is determined by \(p\bmod8\). We then prove
	closed formulas, in terms of class numbers of imaginary quadratic orders and of
	\(\zeta_{\Q(\sqrt p)}(-1)\), for the mass of the curves of each type and for
	the mass of the rational kernels whose quotient is a product of elliptic
	curves. For \(\pi^{2}=\pm p\) the total masses are due to Ibukiyama and to Xue
	and Yu; our proofs use the \(2\)-adic Tate modules, the symmetry of the
	Richelot correspondence and Eichler's mass formula, and, for the other curves,
	modules over orders of a quartic CM field. The kernels are counted through
	gluings of supersingular elliptic curves. As a consequence, the proportion of
	rational kernels with decomposed quotient is of order
	\(p^{-1/2}(1+L(1,\chi))\), where \(\chi\) is the quadratic character of
	\(\Q(\sqrt{-p})\). All
	formulas agree with an exact enumeration for the \(37\) primes
	\(7\le p\le173\). Finally we describe the graph of \(\Fp\)-rational Richelot
	isogenies, which is far from connected in the range computed, and discuss what these
	results say about its use in isogeny-based cryptography.
\end{abstract}

\maketitle

\section{Introduction}\label{sec:intro}

A curve \(C\) of genus two over a field of characteristic \(p\) is
\emph{superspecial} if its Jacobian is isomorphic, over an algebraic closure and
as an unpolarised abelian variety, to a product of two supersingular elliptic
curves. Superspecial curves and the \((2,2)\)-isogenies between their Jacobians
are the vertices and edges of the graphs used in genus-two isogeny-based
cryptography \cite{CastryckDecruSmith2020,CostelloSmith2020,FloritSmith2022auto};
over an algebraically closed field their number and the number of their Richelot
isogenies with decomposed codomain are known
\cite{IbukiyamaKatsuraOort1986,KatsuraTakashima2020}.

This paper is about the curves that are defined over the prime field \(\Fp\) and
about the isogenies that are defined over \(\Fp\). This is the genus-two
counterpart of the set of supersingular elliptic curves over \(\Fp\) with its
rational \(2\)-isogenies \cite{DelfsGalbraith2016}. On it the Frobenius
endomorphism is not a scalar, and it acts on the \(2\)-torsion in a way that is
visible on the six Weierstrass points. The main results are the following.
\begin{enumerate}[label=\textup{(\roman*)},leftmargin=*]
	\item \(\pi^{2}=p\,\alpha_*\) for an automorphism \(\alpha\) of the curve, which
	gives eight Frobenius types of the Weierstrass points, each with its Weil
	polynomial and its number of rational \((2,2)\)-kernels
	(Theorems~\ref{thm:local} and~\ref{thm:class}).
	\item Closed formulas for the mass of the curves of each type
	(Theorem~\ref{thm:vertex}) and for the mass of the rational kernels whose
	quotient is a product of elliptic curves (Theorem~\ref{thm:kernel}), and the
	order of magnitude of the proportion of such kernels
	(Corollary~\ref{cor:prop}).
	\item The graph of \(\Fp\)-rational Richelot isogenies is far from connected in
	the range computed, with a conjectural count of its components
	(Section~\ref{sec:graph}).
\end{enumerate}

\subsection*{Notation}
Throughout, \(p\) is an odd prime, \(k=\Fp\), \(\kbar\) is an algebraic closure
and \(\Frob\colon x\mapsto x^{p}\). For a curve \(C/k\) of genus two with
hyperelliptic involution \(\iota\) we let
\(W=\{P_1,\dots,P_6\}\subset C(\kbar)\) be its Weierstrass points, \(\sigma\) the
permutation of \(W\) induced by \(\Frob\), and \(\lambda\) the cycle type of
\(\sigma\), the \emph{Frobenius type} of \(C\). We write \(J=\Jac(C)\), with
canonical principal polarisation \(\theta\), Rosati involution \(\dagger\), Weil
pairings \(e_n\) on \(J[n]\) for \(p\nmid n\), and Frobenius endomorphism \(\pi\).
For \(\alpha\in\Aut(C_{\kbar})\) we write \(\alpha_*\) for the induced
automorphism of \(J_{\kbar}\) and \(\bar\alpha\) for its image in the reduced
automorphism group \(\Aut(C_{\kbar})/\langle\iota\rangle\subset\PGL_2(\kbar)\). A
subgroup \(K\subset J[2]\) of order four is \emph{maximal isotropic} if \(e_2\) is
trivial on \(K\times K\), and \emph{rational} if \(\pi(K)=K\); the rational
maximal isotropic subgroups are the kernels of the \((2,2)\)-isogenies
\(J\to J/K\) of principally polarised surfaces defined over \(\Fp\). We put
\(\chi_p=(-1)^{(p-1)/2}\).

\subsection*{Frobenius types}
Our first result describes \(\pi^{2}\) on the curve.

\begin{theorem}\label{thm:local}
	Let \(p\) be an odd prime and \(C/\Fp\) a superspecial curve of genus two.
	\begin{enumerate}[label=\textup{(\roman*)},leftmargin=*]
		\item There is a unique \(u\in\End_{\Fp}(J)\) with \(\pi^{2}=pu\); it is an
		automorphism of \((J,\theta)\) of finite order and commutes with \(\pi\).
		There is a unique \(\alpha\in\Aut_{\Fp}(C)\) with \(u=\alpha_*\), and
		\(\alpha|_W=\sigma^{2}\). The following are equivalent: \(\bar\alpha=1\);
		\(u=\pm1\); \(\pi^{2}=\eps p\) with \(\eps\in\{\pm1\}\); \(\sigma^{2}=1\);
		\(\lambda\in\{1^{6},1^{4}2,1^{2}2^{2},2^{3}\}\). If they hold, the Weil
		polynomial of \(J\) is \((T^{2}-\eps p)^{2}\), \(\#C(\Fp)=p+1\),
		\(\#J(\Fp)=(p-\eps)^{2}\) and \(\pi^{\dagger}=\eps\pi\).
	\end{enumerate}
	Assume \(\pi^{2}=\eps p\) in \textup{(ii)--(v)}.
	\begin{enumerate}[label=\textup{(\roman*)},leftmargin=*,start=2]
		\item The maximal isotropic subgroups of \(J[2]\) correspond bijectively to
		the partitions of \(W\) into three pairs, rational subgroups corresponding
		to \(\sigma\)-stable partitions. Their number is \(15\), \(3\), \(3\),
		\(7\) for \(\lambda=1^{6},1^{4}2,1^{2}2^{2},2^{3}\). Frobenius fixes each of
		the three pairs for all of them when \(\lambda\in\{1^{6},1^{4}2\}\), for one
		of the three when \(\lambda=1^{2}2^{2}\), and for one of the seven when
		\(\lambda=2^{3}\); for the others it exchanges two pairs.
		\item If \(\lambda\in\{1^{6},1^{4}2\}\), then \(\eps p\equiv1\pmod 4\).
		\item If \(\eps=+1\), then \(e_2(x,\pi x)=1\) for all \(x\in J[2]\), and
		\(\lambda\in\{1^{6},2^{3}\}\).
		\item If \(\lambda=1^{6}\), then \(\eps=\chi_p\), \(N:=(\pi-1)/2\) lies in
		\(\End_{\Fp}(J)\),
		\[
		N^{2}+N+\frac{1-\eps p}{4}=0,\qquad
		N^{\dagger}=\begin{cases}N,&\eps=+1,\\-1-N,&\eps=-1,\end{cases}
		\]
		and the reduced Frey--R\"uck pairing
		\(\tau\colon J[2]\times J(\Fp)/2J(\Fp)\to\mu_2\) is
		\(\tau(x,y)=e_2(x,Ny)\).
	\end{enumerate}
\end{theorem}

Here \(\tau([D],[E])=f(E)^{(p-1)/2}\) for \(2D=\operatorname{div}f\) with
\(f\in\Fp(C)\) and \(E\) a rational divisor with support disjoint from that of
\(D\). Identifying \(\mu_2\) with \(\F_2\), we regard \(e_2\) and \(\tau\)
as \(\F_2\)-valued bilinear forms. In a basis of \(J[2]\) over \(\F_2\), let
\(B\) and \(\bar N\) be the
matrices of \(e_2\) and of \(N|_{J[2]}\), and \(T\) the matrix of the form
\((x,y)\mapsto\tau(x,y\bmod 2J(\Fp))\) on \(J[2]\times J[2]\), that is, of \(\tau\)
pulled back along \(J[2]\to J(\Fp)/2J(\Fp)\); then statement (v) reads
\(T=B\bar N\). The pairing \(\tau\) itself is perfect, and the map
\(J[2]\to J(\Fp)/2J(\Fp)\) need not be injective.

Without the assumption on \(\sigma^{2}\) we obtain the following list.

\begin{theorem}\label{thm:class}
	Let \(p\ge7\) and let \(C/\Fp\) be a superspecial curve of genus two. Then
	\(\lambda\), \(u\) and the Weil polynomial \(P\) of \(J\) are as in one of the
	rows of Table~\textup{\ref{tab:class}}, where \(\omega=\alpha_{0*}\) for the
	lift \(\alpha_0\) of order three of \(\bar\alpha\), and
	\(\omega^{2}+\omega+1=0\). In particular
	\(\lambda\notin\{1^{3}3,\,1\,2\,3,\,1\,5\}\), the group \(J[2]\) always has a
	rational maximal isotropic subgroup, and \(\#C(\Fp)=p+1\).
\end{theorem}

\begin{table}[ht]
	\centering
	\caption{Frobenius types of superspecial curves of genus two over \(\Fp\),
		\(p\ge7\) (Theorem~\ref{thm:class}). The last column is the number of
		rational maximal isotropic subgroups of \(J[2]\).}\label{tab:class}
	\begin{tabular}{@{}llllll@{}}
		\toprule
		\(\lambda\) & \(\ord\bar\alpha\) & \(u\) & \(P(T)\) & only if & kernels\\
		\midrule
		\(1^{6}\) & \(1\) & \(\chi_p\) & \((T^{2}-\chi_pp)^{2}\) & & \(15\)\\
		\(1^{4}2\) & \(1\) & \(-1\) & \((T^{2}+p)^{2}\) & \(p\equiv3\bmod4\) & \(3\)\\
		\(1^{2}2^{2}\) & \(1\) & \(-1\) & \((T^{2}+p)^{2}\) & & \(3\)\\
		\(2^{3}\) & \(1\) & \(\pm1\) & \((T^{2}\mp p)^{2}\) & & \(7\)\\
		\(1^{2}4\) & \(2\) & \(u^{2}=-1\) & \(T^{4}+p^{2}\) & \(p\equiv1\bmod4\) & \(1\)\\
		\(2\,4\) & \(2\) & \(u^{2}=-1\) & \(T^{4}+p^{2}\) & \(p\equiv3\bmod4\) & \(1\)\\
		\(3^{2}\) & \(3\) & \(\chi_p\,\omega\) & \(T^{4}+\chi_ppT^{2}+p^{2}\) & & \(3\)\\
		\(6\) & \(3\) & \(\pm\omega\) & \(T^{4}\pm pT^{2}+p^{2}\) & & \(1\)\\
		\bottomrule
	\end{tabular}
\end{table}

\begin{corollary}\label{cor:tate}
	In the situation of Theorem~\textup{\ref{thm:local}(v)} the matrix \(\bar N\)
	satisfies \(\bar N^{2}+\bar N+c=0\) with \(c=(1-\eps p)/4\bmod 2\), has rank
	\(4\) if \(p\equiv3,5\pmod 8\) and rank \(2\) if \(p\equiv1,7\pmod8\), and
	\(T+T^{t}=B\) if \(p\equiv3\pmod4\), while \(T\) is symmetric with zero
	diagonal if \(p\equiv1\pmod4\). The triple \((J[2],e_2,\tau)\) is determined up
	to isomorphism by \(p\bmod 8\).
\end{corollary}

The conditions in Table~\ref{tab:class} are necessary conditions. The enumeration
described below shows that each row occurs for some prime \(7\le p\le43\);
Table~\ref{tab:examples} in Section~\ref{sec:crypto} gives an explicit curve for
each row and each Weil polynomial. The
last two rows can be sharpened by a known result: the isogeny class with Weil
polynomial \(T^{4}-pT^{2}+p^{2}\) contains no principally polarised surface
when \(p\equiv1\pmod3\) \cite[Thm.~1.4]{HoweNartRitzenthaler2009}, a restatement
of the main theorem of \cite{HoweMaisnerNartRitzenthaler2008}. Hence
\(u=-\omega\) in the row \(\lambda=6\) occurs only if \(p\equiv2\pmod3\), and
the row \(\lambda=3^{2}\) is empty if \(p\equiv7\pmod{12}\).

\subsection*{What is known}
As far as we could determine from the sources below, what is new here is the
description of \(\pi^{2}\) on the curve and its dictionary with the \(2\)-adic
classes, the counts of rational kernels, the two families of mass formulas, and
the data on the graph. Parts of Theorem~\ref{thm:local} are known in the language
of descent data and of Hermitian lattices, and for these parts our contribution
is a direct proof on the curve. The relation \(\pi^{2}=pu\) with \(u\) an automorphism of finite order of
the polarised variety is implicit in the descent argument of Ibukiyama and
Katsura \cite[proof of Lemma~2.2 and \S3.1]{IbukiyamaKatsura1994}. For
\(\eps=-1\), assertion (iii) corresponds to the statement that the \(2\)-torsion
is rational exactly on the modules over the maximal order
\cite[Lemma~5.11 and the paragraph after its proof]{Robert2024module},
\cite[Cor.~3.2]{Yu2012}; for \(\eps=+1\),
assertion (iv) corresponds to the classification of the \(2\)-adic Tate modules
by Xue and Yu \cite[Prop.~3.2.5]{XueYu2020}. The identity
\(\tau(x,y)=e_2(x,(\pi-1)y')\) with \(2y'=y\) behind (v) is due to Bruin
\cite[Thm.~4.1]{Bruin2011}; see also \cite[Def.~3.21]{Robert2023tate}. The
Weil polynomials of Table~\ref{tab:class}, and with them \(\#C(\Fp)=p+1\), are
among those of the supersingular curves of genus two over \(\Fp\), which are
tabulated in \cite[Table~10.1]{HoweNartRitzenthaler2009}. Corte-Real Santos,
Lorenzon and Reijnders \cite[Thm.~3.26]{CorteRealSantosLorenzonReijnders2026}
characterise the superspecial Jacobians over \(\Fp\) by their Weil polynomial and
the reduced automorphism group of \(C_{\kbar}\): Weil polynomial
\((T^{2}\mp p)^{2}\); or \(T^{4}+p^{2}\) and reduced automorphism group
\(\Z/2\times\Z/2\), \(S_4\) or dihedral of order \(12\); or
\(T^{4}\pm pT^{2}+p^{2}\) and reduced automorphism group \(S_3\), \(S_4\) or
dihedral of order \(12\). Their proof twists a
minimal or maximal curve over \(\F_{p^{2}}\) by an automorphism, using a result
of Howe, Nart and Ritzenthaler
\cite[Lemma~3.24 and Prop.~3.25]{CorteRealSantosLorenzonReijnders2026}; over
\(\kbar\) this is the relation between \(\ord\bar\alpha\) and \(P\) in
Table~\ref{tab:class}. They also recall the group structures of \(J(\Fp)\)
\cite{Zhu2000} and the rank of \(J[2](\Fp)\) in terms of the factorisation of
\(f\) \cite{Cornelissen2001}
\cite[Thm.~3.5 and Lemma~3.7]{CorteRealSantosLorenzonReijnders2026}; together
these two facts give assertion (iii) and the equality \(\eps=\chi_p\) for
\(\lambda=1^{6}\). In these sources we did not find the description of \(u\) on
the curve (\(\alpha|_W=\sigma^{2}\)) as an automorphism defined over \(\Fp\),
the dictionary between the \(2\)-adic classes and the cycle types \(\lambda\),
the counts in (ii), assertion (iv), the matrix form of (v) on the whole cell
\(\lambda=1^{6}\) with Corollary~\ref{cor:tate}, or, in
Theorem~\ref{thm:class}, the complete correspondence between \(\lambda\), \(u\)
and \(P\), the congruence conditions and the excluded types. Robert names the graph of \(\Fp\)-rational isogenies in
dimension two and its levels \cite[Remark~5.8]{Robert2024module}. Our search covered the sources cited here and the recent work on supersingularity
testing \cite{CorteRealSantosLorenzonReijnders2026}, but was not exhaustive, and
the theorem numbers quoted from
\cite{XueYu2020,Robert2024module,Robert2023tate} are those of the preprint
versions.

\subsection*{Masses of rational kernels}
The \emph{mass} of a set of curves over \(\Fp\) is the sum of
\(1/\#\Aut_{\Fp}(C)\) over its \(\Fp\)-isomorphism classes, and similarly for
pairs \((C,K)\). A rational maximal isotropic subgroup \(K\) is
\emph{decomposed} if \(J/K\) is a product of two elliptic curves with the
product polarisation, \emph{split} if the two factors are defined over \(\Fp\)
and \emph{nonsplit} if they are conjugate over \(\F_{p^{2}}\). Let \(H_p\) be
the number of \(\Fp\)-isomorphism classes of supersingular elliptic curves over
\(\Fp\).

Our second result, Theorem~\ref{thm:kernel}, gives the mass of the pairs
\((C,K)\) with \(C\) superspecial and \(K\) rational and decomposed, for each
Weil polynomial, Frobenius type and type of \(K\) (Table~\ref{tab:kernel}). In
particular, on the curves with \(\pi^{2}=\eps p\) the nonsplit kernels have mass
\((p-1)/12-H_p/8\) for each sign, the split kernels occur only for
\(\eps=-1\), with mass \((6h^{2}+2h'^{2})/8-(h+h')/4\) for \(p\equiv3\pmod4\)
and \((h'^{2}-h')/4\) for \(p\equiv1\pmod4\), where \(h=h(-p)\) and
\(h'=h(-4p)\) are class numbers. The proof passes to the dual isogeny: a
decomposed kernel on \(J\) is the same as a gluing of two supersingular elliptic
curves along their \(2\)-torsion \cite{Kani1997,HoweLeprevostPoonen2000} whose
result is a Jacobian, and these are counted on the elliptic curves. The
nonsplit case is the construction of Scholten, used by Costello
\cite{CostelloKummer2018}; we did not find the masses in the literature. As a
by-product we obtain the masses of the superspecial principally polarised
surfaces over \(\Fp\) that are not Jacobians (Proposition~\ref{prop:nonjac}).

\subsection*{Masses of curves}
Our third result, Theorem~\ref{thm:vertex}, gives the mass of the superspecial
curves of genus two over \(\Fp\) in each row of Table~\ref{tab:class}, as an
explicit expression in \(p\), the class numbers \(h(-p)\), \(h(-4p)\),
\(h(-3p)\), \(h(-8p)\), \(h(-12p)\) and \(\zeta_{\Q(\sqrt p)}(-1)\).

For \(\pi^{2}=\pm p\) the total mass of the principally polarised superspecial
surfaces over \(\Fp\) is known: for \(\pi^{2}=-p\) it was computed by Ibukiyama
\cite[Thm.~2]{Ibukiyama1993}, whose proof has one term for each \(2\)-adic
genus, with local data stated without proof, and for \(\pi^{2}=p\) the class
numbers of the genera are due to Xue and Yu \cite[Thm.~1.1]{XueYu2020}. What
Theorem~\ref{thm:vertex} adds is the identification of the \(2\)-adic classes
with the Frobenius types, the masses of the Jacobians alone, and a proof of the
masses of the classes that does not use these results: we classify the
\(2\)-adic Tate modules with their Frobenius and Weil form, obtain the ratios of
the masses of the classes from the symmetry of the Richelot correspondence, and
compute one mass in each connected component from Eichler's mass formula
(Appendix~\ref{sec:vertexpm}). The table of these ratios
(Proposition~\ref{prop:neighbours}) is obtained by a finite exact enumeration
on explicit normal forms of the Tate modules, on which parts (a) and (b) of
Theorem~\ref{thm:vertex} therefore depend. The sum of our masses agrees with
Ibukiyama's theorem.

For \(\pi^{2}\ne\pm p\) the curves have Weil polynomial \(T^{4}+p^{2}\) or
\(T^{4}\pm pT^{2}+p^{2}\), and the surfaces of the isogeny class correspond to
modules over an order of the quartic field \(\Q(\pi)\) by the theorem of
Centeleghe and Stix \cite{CentelegheStix2015}; the superspecial ones are the
modules over \(\Z[\pi,\pi^{2}/p]\), as shown by Xue, Yang and Yu
\cite{XueYangYu2016}. We count the principally polarised ones
(Appendix~\ref{sec:vertexcm}); the class numbers come from Kuroda's class number
formula for biquadratic fields, and the vanishing of one of the masses for
\(p\equiv1\pmod3\) is the theorem of Howe, Maisner, Nart and Ritzenthaler
\cite{HoweMaisnerNartRitzenthaler2008}.

\subsection*{Enumeration}
The formulas of Theorems~\ref{thm:kernel} and~\ref{thm:vertex} were found from
an enumeration. For the \(37\) primes
\(7\le p\le173\) we enumerated all superspecial curves of genus two over \(\Fp\)
up to \(\Fp\)-isomorphism, with their automorphism groups over \(\Fp\), Weil
polynomials and Frobenius types, and all their rational maximal isotropic
subgroups \(K\), with the type of the quotient \(J/K\)
(Section~\ref{sec:comp}). The enumeration works on the curves alone and agrees
with both theorems in every cell.

In this range the proportion of rational kernels with
decomposed quotient lies between \(0.80\,p^{-1/2}\) and \(1.63\,p^{-1/2}\)
(Table~\ref{tab:prop}). For comparison, over \(\bar\F_p\) the proportion of
decomposed Richelot isogenies out of superspecial Jacobians, counted up to
isomorphism, is \(4(p+3)/((p+2)(p+7))\) up to bounded terms in numerator and
denominator \cite[Thm.~6.2]{KatsuraTakashima2020}, and among all pairs of a
genus-two curve over \(\Fp\) and a rational kernel, weighted by automorphisms,
it is \(1/p\) \cite{DangNguyen2026}.

The \(\Fp\)-rational superspecial locus is the object set of the rank-two
module action of \cite{Robert2024module}, and the statements above describe its
\(2\)-torsion and its rational \((2,2)\)-kernels. Section~\ref{sec:crypto}
collects what they imply for the graph of \(\Fp\)-rational Richelot isogenies as
a setting for isogeny-based cryptography: the degrees of its vertices, how to
find all rational kernels, how often products of elliptic curves occur, and the
fact that the graph is far from connected and mixes slowly in the range we
computed. We make no claim about the security or the efficiency of any
cryptographic scheme.

\subsection*{Organisation}
Section~\ref{sec:proofs} contains the proofs of Theorems~\ref{thm:local}
and~\ref{thm:class} and of Corollary~\ref{cor:tate}. Section~\ref{sec:mass}
states the theorem on the masses of curves and compares it with the published
mass formulas. Section~\ref{sec:kernels} states the theorem on the rational kernels and the
data on their proportion, and Section~\ref{sec:kernelproof} proves it.
Appendices~\ref{sec:vertexpm} and~\ref{sec:vertexcm} prove the theorem on the
masses of curves, for \(\pi^{2}=\pm p\) and for \(\pi^{2}\ne\pm p\); they can be
read independently of Sections~\ref{sec:kernels} to~\ref{sec:comp}, and
Theorem~\ref{thm:vertex} is stated in Section~\ref{sec:mass}.
Section~\ref{sec:graph} describes the graph of the \(\Fp\)-rational Richelot
isogenies, Section~\ref{sec:crypto} discusses its relevance for isogeny-based
cryptography, and Section~\ref{sec:comp} describes the computations.

\section{Proofs of the structure theorems}\label{sec:proofs}

\begin{lemma}\label{lem:frob2}
	Let \(A/\Fp\) be an abelian variety such that \(A_{\kbar}\) is isomorphic to a
	product of supersingular elliptic curves. Then \(\pi_A^{2}=u\circ[p]\) for a
	unique \(u\in\End_{\Fp}(A)\); it is an automorphism and commutes with
	\(\pi_A\).
\end{lemma}

\begin{proof}
	Let \(E/\kbar\) be supersingular. Then \([p]_E\) has degree \(p^{2}\) and
	trivial kernel on \(\kbar\)-points, so it is purely inseparable and factors as
	\([p]_E=\ell_E\circ\mathrm{Fr}_E\) with
	\(\mathrm{Fr}_E\colon E\to E^{(p^{2})}\) the \(p^{2}\)-power Frobenius and
	\(\ell_E\) an isomorphism \cite[II.2.12]{Silverman2009}. Hence
	\([p]_X=\ell\circ\mathrm{Fr}_X\) with \(\ell\) an isomorphism for a product
	\(X\) of supersingular elliptic curves. Let \(\varphi\colon A_{\kbar}\to X\) be
	an isomorphism of abelian varieties. Since
	\(\mathrm{Fr}_X\circ\varphi=\varphi^{(p^{2})}\circ\mathrm{Fr}_{A_{\kbar}}\) and
	\(\mathrm{Fr}_{A_{\kbar}}=\pi_A^{2}\), because \(A\) is defined over \(\Fp\), we
	get \([p]_A=v\circ\pi_A^{2}\) with \(v=\varphi^{-1}\ell\varphi^{(p^{2})}\) an
	isomorphism of varieties fixing \(0\), hence an automorphism of \(A_{\kbar}\).
	Put \(u=v^{-1}\). For \(t\in\operatorname{Gal}(\kbar/\Fp)\) we have
	\(t(u)\circ[p]=t(\pi_A^{2})=u\circ[p]\), and \([p]\) is surjective, so
	\(t(u)=u\); the same argument gives uniqueness. Finally
	\((\pi_Au-u\pi_A)\circ[p]=\pi_A^{3}-\pi_A^{3}=0\).
\end{proof}

\begin{lemma}\label{lem:torelli}
	Let \(K\) be a perfect field and \(C/K\) a curve of genus two. Then
	\(\alpha\mapsto\alpha_*\) is a bijection \(\Aut_K(C)\to\Aut_K(J,\theta)\), and
	\(\iota_*=-1\).
\end{lemma}

\begin{proof}
	Fix \(P\in C(\bar K)\) and \(f\colon Q\mapsto[Q-P]\). For
	\(\alpha\in\Aut(C_{\bar K})\) we have
	\(f\circ\alpha=\alpha_*\circ f+[\alpha P-P]\). By Torelli's theorem in the form
	\cite[III, Thm.~12.1]{MilneAV}, for every \(\beta\in\Aut(J_{\bar K},\theta)\)
	there are \(\alpha\in\Aut(C_{\bar K})\) and \(c\in J(\bar K)\) with
	\(f\circ\alpha=\pm\beta\circ f+c\), and, \(C\) being hyperelliptic, \(\alpha\)
	and \(c\) are unique once the sign \(+\) is fixed. Evaluating at \(P\) gives
	\(c=[\alpha P-P]\), hence \(\beta=\alpha_*\) on \(f(C)\), which generates
	\(J\). Uniqueness gives injectivity and, applied to
	\((t\alpha)_*=t(\alpha_*)\) for \(t\in\operatorname{Gal}(\bar K/K)\), shows
	that \(\alpha\) is defined over \(K\) if \(\alpha_*\) is. Since \(Q+\iota Q\)
	is a canonical divisor for every \(Q\), \(\iota_*=-1\).
\end{proof}

\begin{lemma}\label{lem:J2}
	\begin{enumerate}[label=\textup{(\alph*)},leftmargin=*]
		\item \(J[2]=\{0\}\cup\{[P_i-P_j]:i<j\}\), these fifteen classes are distinct
		and non-zero, and \([P_i-P_j]+[P_k-P_l]+[P_m-P_n]=0\) if
		\(\{i,j,k,l,m,n\}=\{1,\dots,6\}\).
		\item \(e_2([P_i-P_j],[P_k-P_l])=(-1)^{\#(\{i,j\}\cap\{k,l\})}\).
		\item The maximal isotropic subgroups of \(J[2]\) are the fifteen groups
		\(\{0,[P_a-P_b],[P_c-P_d],[P_e-P_f]\}\) attached to the partitions
		\(\{a,b\},\{c,d\},\{e,f\}\) of \(W\) into three pairs.
		\item If \(\Frob\) or an automorphism of \(C_{\kbar}\) permutes \(W\) by
		\(s\), it acts on \(J[2]\) by \([P_i-P_j]\mapsto[P_{s(i)}-P_{s(j)}]\), and
		\(s\) is determined by this action.
	\end{enumerate}
\end{lemma}

\begin{proof}
	Work over \(\kbar\) with a sextic model \(y^{2}=f(x)\) and \(a_i=x(P_i)\). From
	\(\operatorname{div}(x-a_i)=2P_i-D_\infty\) and
	\(\operatorname{div}(y)=\sum_iP_i-3D_\infty\), with \(D_\infty\) a canonical
	divisor, we get \(2P_i\sim D_\infty\) and \(\sum_iP_i\sim3D_\infty\); this
	gives \(2[P_i-P_j]=0\) and the relation in (a). Two distinct linearly
	equivalent effective divisors of degree two are canonical by Riemann--Roch,
	hence of the form \(Q+\iota Q\). So \([P_i-P_j]=[P_k-P_l]\), that is
	\(P_i+P_l\sim P_j+P_k\), forces \(\{k,l\}=\{i,j\}\), and \([P_i-P_j]\ne0\); as
	\(\#J[2]=16\) this proves (a). For disjoint \(\{i,j\}\), \(\{k,l\}\) put
	\(F=(x-a_i)/(x-a_j)\) and \(G=(x-a_k)/(x-a_l)\); the formula
	\(e_2([D],[E])=G(D)/F(E)\) for \(2D=\operatorname{div}F\),
	\(2E=\operatorname{div}G\) with disjoint supports
	\cite{Howe1996}, \cite[(4.1)]{Bruin2011} gives the value
	\[
	\frac{(a_i-a_k)(a_j-a_l)}{(a_i-a_l)(a_j-a_k)}\cdot
	\frac{(a_k-a_j)(a_l-a_i)}{(a_k-a_i)(a_l-a_j)}=1 .
	\]
	Since \(e_2\) is alternating and non-degenerate, the orthogonal complement of
	\(x=[P_1-P_2]\) has eight elements; it contains \(0\), \(x\) and the six
	classes supported in \(\{3,4,5,6\}\), so the remaining eight classes, those
	meeting \(\{1,2\}\) in one point, pair non-trivially with \(x\). This is (b).
	Then (c) follows: three non-zero pairwise orthogonal classes correspond to
	pairwise disjoint pairs, and conversely by the relation in (a). In (d) the
	first assertion is clear, and a permutation of six points fixing every
	two-element subset is trivial.
\end{proof}

\begin{lemma}\label{lem:lift}
	Let \(K\) be a field of odd characteristic, \(C\colon y^{2}=f(x)\) of genus two
	over \(K\) and \(F(X,Z)\) the binary sextic with \(F(x,1)=f(x)\). Every
	\(\alpha\in\Aut(C_{\bar K})\) induces \(\bar\alpha\in\PGL_2(\bar K)\) permuting
	the branch points as \(\alpha\) permutes \(W\); \(\bar\alpha=1\) if and only if
	\(\alpha\in\{1,\iota\}\); a non-trivial \(\bar\alpha\) permutes the branch
	points non-trivially; and \(\bar\alpha\in\PGL_2(K)\) if \(\alpha\) is defined
	over \(K\). If
	\(M=\bigl(\begin{smallmatrix}a&b\\c&d\end{smallmatrix}\bigr)\in\operatorname{GL}_2(K)\)
	permutes the branch points, then \(F(aX+bZ,cX+dZ)=\kappa_MF(X,Z)\) with
	\(\kappa_M\in K^{\times}\), the lifts of \(M\) are
	\((x,y)\mapsto\bigl(\tfrac{ax+b}{cx+d},\,e\,y/(cx+d)^{3}\bigr)\) with
	\(e^{2}=\kappa_M\), and they are defined over \(K\) if and only if
	\(\kappa_M\) is a square in \(K\).
\end{lemma}

\begin{proof}
	The canonical map \(C\to\PP^{1}\) is intrinsic and Galois-equivariant, its deck
	group is \(\langle\iota\rangle\), and a M\"obius transformation fixing three
	points is trivial. The forms \(F\circ M\) and \(F\) have the same zeros, hence
	are proportional, and for a lift \(\alpha\) the function \(Y=y\circ\alpha\)
	satisfies \(Y^{2}=f(Mx)=\kappa_Mf(x)/(cx+d)^{6}\).
\end{proof}

\begin{lemma}\label{lem:irred}
	For an odd prime \(p\) and \(\eps\in\{\pm1\}\) the polynomials \(T^{2}-\eps p\)
	and \(T^{4}+p^{2}\) are irreducible over \(\Q\), and so is
	\(T^{4}+\eps pT^{2}+p^{2}\) unless \((\eps,p)=(-1,3)\). If \(\varphi\) is an
	endomorphism of an abelian surface and a root of a monic irreducible
	\(g\in\Z[X]\) of degree \(2\) or \(4\), the characteristic polynomial of
	\(\varphi\) is \(g^{4/\deg g}\).
\end{lemma}

\begin{proof}
	A root \(\beta\) of \(T^{4}+p^{2}\) has \(\beta^{2}=\pm p\sqrt{-1}\), and
	\(p\sqrt{-1}=p(1+\sqrt{-1})^{2}/2\) is a square in \(\Q(\sqrt{-1})\) only if
	\(2p\in\pm\Q^{\times2}\). A root \(\beta\) of \(T^{4}+\eps pT^{2}+p^{2}\) has
	\(\beta^{2}=\eps p\zeta\) with \(\zeta=(\zeta^{2})^{2}\) a primitive cube root
	of unity, and \(\eps p\) is a square in \(\Q(\sqrt{-3})\) only if
	\(\eps p\in\Q^{\times2}\cup-3\,\Q^{\times2}\). In both cases
	\([\Q(\beta):\Q]=4\). The last assertion holds because the characteristic
	polynomial is monic of degree four with integer coefficients and all its roots
	are roots of \(g\) \cite[I, Thm.~10.9 and Prop.~10.20]{MilneAV}.
\end{proof}

\begin{proof}[Proof of Theorem~\textup{\ref{thm:local}}]
	(i) Lemma~\ref{lem:frob2} gives \(u\). By \cite[II, Lemma~1.2]{MilneAV},
	\(\pi^{\dagger}\pi=p\); hence \(\pi^{\dagger}\) commutes with \(\pi\),
	\((\pi^{2})^{\dagger}\pi^{2}=p^{2}\) and \(u^{\dagger}u=1\), that is,
	\(u\in\Aut(J,\theta)\), a finite group \cite[I, Prop.~14.4]{MilneAV}.
	Lemma~\ref{lem:torelli} gives a unique \(\alpha\in\Aut_{\Fp}(C)\) with
	\(\alpha_*=u\). As \(p\) is odd, \([p]\) is the identity on \(J[2]\), so
	\(\pi^{2}=u\) on \(J[2]\); by Lemma~\ref{lem:J2}(d) the permutations
	\(\sigma^{2}\) and \(\alpha|_W\) coincide. By Lemma~\ref{lem:lift},
	\(\bar\alpha\) acts on the branch points as \(\sigma^{2}\) and has the same
	order, and \(\bar\alpha=1\) if and only if \(\alpha\in\{1,\iota\}\), if and
	only if \(u=\pm1\) (Lemma~\ref{lem:torelli}); a permutation of six points has
	trivial square exactly for the four listed cycle types. If \(\pi^{2}=\eps p\),
	Lemma~\ref{lem:irred} gives the Weil polynomial \((T^{2}-\eps p)^{2}\); its
	\(T^{3}\)-coefficient vanishes and its value at \(1\) is \((1-\eps p)^{2}\),
	and \(\pi^{\dagger}=p\pi^{-1}=\eps\pi\).

	(ii) The bijection and the rationality criterion are
	Lemma~\ref{lem:J2}(c),(d). A \(\sigma\)-stable partition into pairs is a
	fixed-point-free involution \(m\) commuting with \(\sigma\); it permutes the
	cycles of \(\sigma\). For \(\sigma=1\) there are \(15\). For \(\sigma=(12)\),
	\(m\) contains \((12)\) and pairs the four fixed points: three ways. For
	\(\sigma=(12)(34)\), \(m\) contains \((56)\) and is \((12)(34)\),
	\((13)(24)\) or \((14)(23)\); \(\sigma\) fixes the three pairs in the first
	case and exchanges two of them in the others. For \(\sigma=(12)(34)(56)\)
	either \(m=\sigma\), or \(m\) fixes one cycle and matches the other two
	equivariantly: \(1+3\cdot2\) partitions, with the stated behaviour.

	(iii) \(J(\Fp)\) contains the \(\sigma\)-fixed part of \(J[2]\), of order
	\(16\) for \(\lambda=1^{6}\) and \(8\) for \(\lambda=1^{4}2\) (the classes of
	the \(\sigma\)-stable pairs, and \(0\)). Hence \(8\) divides
	\(\#J(\Fp)=(1-\eps p)^{2}\), so \(4\mid 1-\eps p\).

	(iv) Here \(\pi^{\dagger}=\pi\). Let \(x\in J[2]\) and \(x'\in J[4]\) with
	\(2x'=x\). Then
	\(e_4(x',\pi x')=e_4(\pi^{\dagger}x',x')=e_4(\pi x',x')=e_4(x',\pi x')^{-1}\),
	so \(e_2(x,\pi x)=e_4(x',\pi x')^{2}=1\)
	\cite[I, Lemma~13.1 and \S14]{MilneAV}. By Lemma~\ref{lem:J2}(b),(d),
	\(e_2(x,\pi x)=(-1)^{\#(\{i,j\}\cap\{\sigma i,\sigma j\})}\) for
	\(x=[P_i-P_j]\); this is \(-1\) for \(\sigma=(12)\), \(\{i,j\}=\{1,3\}\) and
	for \(\sigma=(12)(34)\), \(\{i,j\}=\{1,5\}\).

	(v) By (iii), \(4\mid1-\eps p\), that is, \(\eps=\chi_p\). Since \(\sigma=1\),
	\(\pi-1\) kills \(J[2]\), hence \(\pi-1=2N\) with \(N\in\End_{\Fp}(J)\)
	\cite[I, Lemma~10.16]{MilneAV}. Then \(4(N^{2}+N)=\pi^{2}-1=\eps p-1\) and
	\(N^{\dagger}=(\eps\pi-1)/2\). By \cite[Thm.~4.1]{Bruin2011},
	\(\tau(x,y)=e_2(x,(\pi-1)y')\) for any \(y'\in J(\kbar)\) with \(2y'=y\), and
	\((\pi-1)y'=N(2y')=Ny\).
\end{proof}

\begin{proof}[Proof of Theorem~\textup{\ref{thm:class}}]
	An element of finite order of \(\PGL_2(\kbar)\) is conjugate either to
	\(x\mapsto\zeta x\), with two fixed points and all other orbits of length
	\(\ord\zeta\), or to \(x\mapsto x+1\), with one fixed point and all other
	orbits of length \(p\). As \(p\ge7\), the second kind cannot permute six points
	non-trivially, and the first kind does so with cycle type \(2^{3}\),
	\(1^{2}2^{2}\), \(3^{2}\), \(1^{2}4\), \(1\,5\) or \(6\). By
	Theorem~\ref{thm:local}(i) this applies to \(\sigma^{2}\). A cycle of
	\(\sigma\) of odd length \(l\) contributes a cycle of length \(l\) to
	\(\sigma^{2}\), and a cycle of even length \(2m\) two cycles of length \(m\);
	so the types \(2^{3}\), \(1^{2}4\) and \(6\) are not squares, and either
	\(\sigma^{2}=1\), or \(\lambda\in\{1^{2}4,\,2\,4\}\) with \(\sigma^{2}\) of type
	\(1^{2}2^{2}\), or \(\lambda\in\{3^{2},6\}\) with \(\sigma^{2}\) of type
	\(3^{2}\), or \(\lambda=1\,5\).

	If \(\sigma^{2}=1\) the first four rows follow from
	Theorem~\ref{thm:local}(i),(iii),(iv).

	Let \(\lambda=1^{2}4\). The two rational Weierstrass points are the fixed
	points of the involution \(\bar\alpha\); placing them at \(0\) and \(\infty\)
	we get \(\bar\alpha(x)=-x\) and \(f(-x)=-f(x)\), so \(\kappa_M=-1\) for
	\(M=\operatorname{diag}(-1,1)\). Let \(\lambda=2\,4\). The fixed points of
	\(\bar\alpha\) are the conjugate pair of Weierstrass points, which we place at
	\(\pm\sqrt d\) with \(d\in\Fp\) a non-square; then \(\bar\alpha(x)=d/x\) and
	\(f=c\,(x^{2}-d)\,g(x)\), where the roots of the monic quartic \(g\) are
	exchanged in pairs by \(x\mapsto d/x\), so that \(g(0)=d^{2}\) and
	\(x^{4}g(d/x)=d^{2}g(x)\); hence \(x^{6}f(d/x)=-d^{3}f(x)\) and
	\(\kappa_M=-d^{3}\). Since \(\alpha\) is defined over \(\Fp\),
	Lemma~\ref{lem:lift} shows that \(-1\) is a square in the first case and a
	non-square in the second. In both cases \(\alpha^{2}\) lies over the identity
	of \(\PP^{1}\) and \(\alpha^{2}\colon(x,y)\mapsto(x,-y)\), so \(u^{2}=-1\),
	\(\pi^{4}+p^{2}=0\), and \(P=T^{4}+p^{2}\) by Lemma~\ref{lem:irred}.

	Let \(\lambda\in\{3^{2},6\}\). Exactly one of the lifts \(\alpha\),
	\(\iota\alpha\) of \(\bar\alpha\) has order three; call it \(\alpha_0\) and put
	\(\omega=\alpha_{0*}\), so that \(u=\eps\omega\) with \(\eps=\pm1\). The
	characteristic polynomial of \(\omega\) is \((X-1)^{m}(X^{2}+X+1)^{(4-m)/2}\)
	with \(m\in\{0,2,4\}\), and \(m\ne4\) because \(\omega\ne1\). If \(m=2\), the
	kernel of \(\omega-1\) on the Tate module \(T_2J\) is a saturated submodule of
	rank two and reduces to a non-zero subspace of \(J[2]\) fixed by \(\omega\);
	but \(\omega\) acts on \(J[2]\) through a permutation of type \(3^{2}\), which
	fixes no pair. Hence \(\omega^{2}+\omega+1=0\),
	\(\pi^{4}+\eps p\pi^{2}+p^{2}=p^{2}(\omega^{2}+\omega+1)=0\), and
	\(P=T^{4}+\eps pT^{2}+p^{2}\) by Lemma~\ref{lem:irred}. Put
	\(\psi=\omega\pi\). Then \(\psi^{2}=\omega^{2}\pi^{2}=\eps p\,\omega^{3}=\eps p\)
	and \(\psi\) acts on \(J[2]\) through \(\sigma^{3}\). If \(\lambda=3^{2}\)
	then \(\sigma^{3}=1\), so \(J[2]\) lies in the kernel of \(\psi-1\), whose
	degree is \((1-\eps p)^{2}\) by Lemma~\ref{lem:irred}; thus
	\(4\mid1-\eps p\) and \(\eps=\chi_p\).

	Let \(\lambda=1\,5\), and let \(\alpha_0\) be the lift of order five of
	\(\bar\alpha\), \(z=\alpha_{0*}\), \(u=\eps z\). As before \(z\) has no
	eigenvalue \(1\) on \(T_2J\), because a permutation of type \(1\,5\) fixes no
	pair, so the characteristic polynomial of \(z\) is the cyclotomic polynomial
	\(\Phi_5\) and \(\Q[z]\cong\Q(\zeta_5)\). The endomorphism \(\psi=z^{2}\pi\)
	commutes with \(z\) and \(\psi^{2}=\eps p\,z^{5}=\eps p\). As \(p\ne5\),
	\(\eps p\) is not a square in \(\Q(\zeta_5)\), so \(\Q[z,\psi]\) is a field of
	degree eight. It acts faithfully on the four-dimensional space
	\(T_\ell J\otimes\Q_\ell\) \cite[I, Thm.~10.15]{MilneAV}, which is impossible.

	The numbers of rational subgroups are the numbers of \(\sigma\)-stable
	partitions into pairs: a partition stable under a \(4\)-cycle and the
	transposition of the remaining two points, or under a \(6\)-cycle, is unique,
	and for \(\lambda=3^{2}\) the three partitions pairing the two cycles
	equivariantly are the stable ones.
\end{proof}

\begin{proof}[Proof of Corollary~\textup{\ref{cor:tate}}]
	The relation and the adjoint of \(\bar N\) come from
	Theorem~\ref{thm:local}(v), and in additive notation
	\(\tau(x,y)+\tau(y,x)=e_2(x,(\bar N+\bar N^{\dagger})y)\). For \(\eps=+1\) the
	argument of (iv) with \(N\) in place of \(\pi\) gives \(e_2(x,Nx)=1\). If
	\(c=1\), then \(\bar N(\bar N+1)=1\) and \(\bar N\) is invertible. If \(c=0\),
	then \(\eps p=s^{2}\) with \(s\in\Z_2\), \(s\equiv1\pmod4\) after changing
	its sign, and \(N\) is annihilated by \((X-n_+)(X-n_-)\) with
	\(n_\pm=(\pm s-1)/2\in\Z_2\). As \(n_+\equiv0\) and \(n_-\equiv1\pmod 2\), the
	difference \(n_+-n_-=s\) is a unit, so the idempotents
	\((N-n_\mp)/(n_\pm-n_\mp)\) lie in \(\Z_2[N]\subset\End(T_2J)\), and \(T_2J\) is
	the direct sum of the eigenmodules of \(N\) for \(n_+\) and \(n_-\). Both have
	rank two because the Weil polynomial is \((T^{2}-\eps p)^{2}\), and \(\bar N\)
	is \(0\) on the reduction of the first and \(1\) on the reduction of the
	second. (The splitting uses that \(N=(\pi-1)/2\) is an endomorphism; the
	eigenvalues \(\pm s\) of \(\pi\) itself are congruent modulo \(2\).) Finally, the
	matrices \(\bar N\in M_4(\F_2)\) with the stated relation, adjoint, rank and
	(for \(\eps=+1\)) vanishing of \(e_2(x,\bar Nx)\) form a single orbit under the
	symplectic group in each of the four cases (\(40\), \(120\), \(12\) and \(20\)
	matrices for \(p\equiv3,7,5,1\bmod 8\)), as one checks by listing them
	(Section~\ref{sec:comp}).
\end{proof}

\begin{remark}\label{rem:small}
	The hypothesis \(p\ge7\) is used only in Theorem~\ref{thm:class}. There is no
	superspecial curve of genus two in characteristic three. For \(p=5\),
	Theorem~\ref{thm:local} holds, and the enumeration of the superspecial curves
	over \(\F_5\) shows one further row: \(\lambda=1\,5\), with \(\bar\alpha\)
	unipotent of order five and
	\(P=T^{4}\pm5T^{3}+15T^{2}\pm25T+25\). The arguments extend to \(\F_q\) with
	\(q=p^{a}\) an odd power of \(p\), with \(q\) in place of \(p\) in the formulas
	and congruences and with \([q]=[p]^{a}\) in Lemma~\ref{lem:frob2}, still
	under the hypothesis \(p\ge7\) in Theorem~\ref{thm:class}: the exclusion of
	\(\lambda=1\,5\) uses that \(\eps q\) is not a square in \(\Q(\zeta_5)\),
	which fails for \(q=5^{a}\), \(\eps=+1\). We do not write this out and do not
	use it. For even powers the assertions on \(P\) and
	on \(\#C(\F_q)\) fail.
\end{remark}

\section{Masses of curves by Frobenius type}\label{sec:mass}

For a negative discriminant \(D\) let \(h(D)\) be the class number of the
quadratic order of discriminant \(D\), and let
\(z_p=\zeta_{\Q(\sqrt p)}(-1)\). For \(p\equiv3\pmod4\), \(p>3\), one has
\(h(-4p)=\bigl(2-\leg2p\bigr)h(-p)\). Let \(H_p\) be the number of
\(\Fp\)-isomorphism classes of supersingular elliptic curves over \(\Fp\):
\(H_p=h(-p)+h(-4p)\) if \(p\equiv3\pmod4\) and \(H_p=h(-4p)\) if
\(p\equiv1\pmod4\) \cite{DelfsGalbraith2016}.

The mass of a finite set \(\mathcal S\) of curves over \(\Fp\), closed under
\(\Fp\)-isomorphism, is
\(\Mass(\mathcal S)=\sum_{[C]}1/\#\Aut_{\Fp}(C)\), the sum being over the
\(\Fp\)-isomorphism classes. A curve and its quadratic twist have the same
branch divisor, the same Frobenius type and the same value of \(\pi^{2}\), and
the mass is computed on branch divisors:

\begin{lemma}\label{lem:massdiv}
	Let \(B\) be a reduced divisor of degree six on \(\PP^{1}\) over \(\Fp\) and
	\(S\) its stabiliser in \(\PGL_2(\Fp)\). The curves of genus two over \(\Fp\)
	with branch divisor in the \(\PGL_2(\Fp)\)-orbit of \(B\) have total mass
	\(1/\#S\).
\end{lemma}

\begin{proof}
	These curves are \(C\colon y^{2}=f\) and its twist \(C'\colon y^{2}=cf\) with
	\(c\) a non-square, up to isomorphism. The map \(\Aut_{\Fp}(C)\to S\) has
	kernel \(\langle\iota\rangle\), and by Lemma~\ref{lem:lift} its image
	\(S_C\) is the kernel of the character \(M\mapsto\kappa_M\) modulo squares,
	which is the same for \(C\) and \(C'\). If \(S_C=S\), then \(C\) and \(C'\)
	are not isomorphic and each has mass \(1/(2\#S)\). If \(S_C\ne S\), an element
	of \(S\smallsetminus S_C\) lifts to an isomorphism \(C\to C'\), and the single
	class has mass \(1/(2\#S_C)=1/\#S\).
\end{proof}

\begin{theorem}\label{thm:vertex}
	Let \(p\ge7\) be a prime, and write \(h=h(-p)\) when \(p\equiv3\pmod4\) and
	\(h'=h(-4p)\). The masses of the superspecial curves of genus two over \(\Fp\)
	in the rows of Table~\textup{\ref{tab:class}} are as follows; a pair
	\((\lambda,P)\) that is not listed for a given residue of \(p\) has mass zero.
	\begin{enumerate}[label=\textup{(\alph*)},leftmargin=*]
		\item \(p\equiv3\pmod4\), \(P=(T^{2}+p)^{2}\):
		\[
		\begin{aligned}
		\lambda=1^{6}&:\ \frac{(p-1)h}{48}-\frac{h^{2}}{8}, &
		\lambda=1^{4}2&:\ \frac{\bigl(p-\leg2p\bigr)h}{16}-\frac{hh'}{4},\\
		\lambda=1^{2}2^{2}&:\ \frac{(p-1)h'}{16}-\frac{h'^{2}}{8}, &
		\lambda=2^{3}&:\ \frac{(p-1)(h'-1)}{48};
		\end{aligned}
		\]
		\(P=(T^{2}-p)^{2}\), \(\lambda=2^{3}\): \(\dfrac{z_p}{4}-\dfrac{p-1}{48}\).
		\item \(p\equiv1\pmod4\), \(P=(T^{2}+p)^{2}\):
		\[
		\lambda=1^{2}2^{2}:\ \frac{p\,h'}{16}-\frac{h'^{2}}{8},\qquad
		\lambda=2^{3}:\ \frac{(p-1)(h'-1)}{48};
		\]
		\(P=(T^{2}-p)^{2}\):
		\[
		\lambda=1^{6}:\ \frac{z_p}{4},\qquad
		\lambda=2^{3}:\ \frac{\bigl(8-2\leg2p\bigr)z_p}{4}-\frac{p-1}{48}.
		\]
		\item \(P=T^{4}+p^{2}\): \(\lambda=2\,4\) has mass \((h(-8p)-2)/8\) if
		\(p\equiv3\pmod4\), and \(\lambda=1^{2}4\) has mass \(h(-8p)/8\) if
		\(p\equiv1\pmod4\).
		\item \(P=T^{4}+pT^{2}+p^{2}\): \(\lambda=3^{2}\) has mass \(h(-3p)/12\) if
		\(p\equiv1\pmod4\), and \(\lambda=6\) has mass
		\(\bigl(h(-12p)-2\,[p\equiv2\bmod3]\bigr)/12\).
		\item \(P=T^{4}-pT^{2}+p^{2}\): \(\lambda=3^{2}\) has mass \(h(-p)/6\) if
		\(p\equiv11\pmod{12}\), and \(\lambda=6\) has mass \((h(-4p)-1)/6\) if
		\(p\equiv2\pmod3\).
	\end{enumerate}
\end{theorem}

Parts (a) and (b) are proved in Appendix~\ref{sec:vertexpm}
(Theorem~\ref{thm:vertexpm}), and parts (c)--(e) in Appendix~\ref{sec:vertexcm}
(Theorem~\ref{thm:vertexcmcells}).

\begin{proposition}\label{prop:vertexcheck}
	For the \(37\) primes \(7\le p\le173\) the enumeration of
	Section~\textup{\ref{sec:comp}} gives the masses of
	Theorem~\textup{\ref{thm:vertex}} in every cell. For these primes the total
	mass of the superspecial curves of genus two over \(\Fp\) is an integer, equal
	to the number of superspecial curves over \(\bar\F_p\) with field of moduli
	\(\Fp\).
\end{proposition}

\begin{proof}
	By the enumeration; Table~\ref{tab:vertex} lists the
	masses. The total is an integer because the sum of \(1/\#\Aut_{\Fp}\) over the
	\(\Fp\)-forms of a curve over \(\bar\F_p\) with a model over \(\Fp\) is \(1\).
\end{proof}

The enumeration does not use the theory of Appendices~\ref{sec:vertexpm}
and~\ref{sec:vertexcm}, so Proposition~\ref{prop:vertexcheck} is an independent
check of Theorem~\ref{thm:vertex}. The formulas were first found from it: (a)
and (b) on the primes \(p\le31\), and (c)--(e) on the \(27\) primes
\(p\le113\).

Let \(M(2)\) be the mass of the principally polarised superspecial abelian
surfaces over \(\Fp\) with \(\pi^{2}=-p\). Ibukiyama
\cite[Thm.~2]{Ibukiyama1993} proved that, for \(p>3\),
\[
M(2)=
\begin{cases}
(4p-1)\,h(-4p)/48, & p\equiv1\pmod4,\\
(8p-5)\,h(-p)/24, & p\equiv3\pmod8,\\
(p-1)\,h(-p)/6, & p\equiv7\pmod8.
\end{cases}
\]
A principally polarised abelian surface is a Jacobian or a product of two
elliptic curves with the product polarisation, and over \(\Fp\) the product may
be a product of two elliptic curves over \(\Fp\) or the Weil restriction of an
elliptic curve over \(\F_{p^{2}}\).

\begin{proposition}\label{prop:ibukiyama}
	For the \(37\) primes \(7\le p\le173\), the mass of the superspecial curves of
	genus two over \(\Fp\) with \(\pi^{2}=-p\), plus \(H_p^{2}/8\), plus
	\((p-1)/48\), equals \(M(2)\).
\end{proposition}

\begin{proof}
	By the enumeration, with the class numbers computed by counting reduced
	forms.
\end{proof}

By Proposition~\ref{prop:nonjac} below, \(H_p^{2}/8\) is the mass of the
products of two supersingular elliptic curves over \(\Fp\), and \((p-1)/48\) is
the mass of the Weil restrictions with \(\pi^{2}=-p\); these are the
principally polarised superspecial surfaces with \(\pi^{2}=-p\) that are not
Jacobians. So Proposition~\ref{prop:ibukiyama} says that the enumeration agrees
with Ibukiyama's theorem, read as a formula for the mass of all principally
polarised superspecial surfaces over \(\Fp\) with \(\pi^{2}=-p\), each weighted
by the inverse of the order of its automorphism group over \(\Fp\).

The proof in \cite{Ibukiyama1993} writes \(M(2)\) as a sum of terms indexed by
\(2\)-adic genera, with local data stated without proof
\cite[Prop.~4.3 and~4.4]{Ibukiyama1993}; for \(\pi^{2}=-p\) and
\(p\equiv3\pmod4\) the four classes of self-dual \(2\)-adic Tate modules are
classified by Yu \cite{Yu2012}, and for \(\pi^{2}=p\) the genera and their class
numbers are in \cite{XueYu2020}. For \(p\equiv7\pmod8\), adding
\(h^{2}/8\), \(hh'/4\), \(h'^{2}/8\) and \((p-1)/48\) to the four formulas of
Theorem~\ref{thm:vertex}(a) with \(P=(T^{2}+p)^{2}\) gives
\((p-1)h/48\) times \(1\), \(3\), \(3\), \(1\) for
\(\lambda=1^{6},1^{4}2,1^{2}2^{2},2^{3}\); the added terms sum to
\(H_p^{2}/8+(p-1)/48\), and the factors \(1,3,3,1\) are the quotients of the
local data at \(2\) in \cite[Prop.~4.4]{Ibukiyama1993}. Appendix~\ref{sec:vertexpm} proves the two statements
needed to pass from masses of genera to Theorem~\ref{thm:vertex}(a),(b): the
\(2\)-adic classes correspond to the Frobenius types \(\lambda\)
(Proposition~\ref{prop:twoadic}), and the products and the Weil restrictions lie
in the classes indicated by the subtracted terms. It also computes the masses of
the classes from Eichler's mass formula, without using
\cite{Ibukiyama1993}; the sum of the masses obtained for \(\pi^{2}=-p\) is
\(M(2)\). The non-Jacobians have mass of order \(h^{2}+p\), against \(M(2)\) of
order \(ph\).

\begin{table}[tp]
	\centering
	\scriptsize
	\caption{Masses of the superspecial curves of genus two over \(\Fp\) by Weil
		polynomial and Frobenius type \(\lambda\). Columns 2--5:
		\(P=(T^{2}+p)^{2}\); columns 6--7: \(P=(T^{2}-p)^{2}\); columns 8--9:
		\(P=T^{4}+p^{2}\); columns 10--11: \(P=T^{4}+pT^{2}+p^{2}\); columns
		12--13: \(P=T^{4}-pT^{2}+p^{2}\).}\label{tab:vertex}
	\setlength{\tabcolsep}{3pt}
	\begin{tabular}{@{}r|rrrr|rr|rr|rr|rr|r@{}}
		\toprule
		\(p\) & \(1^{6}\) & \(1^{4}2\) & \(1^{2}2^{2}\) & \(2^{3}\) & \(1^{6}\) & \(2^{3}\)
		& \(1^{2}4\) & \(2\,4\) & \(3^{2}\) & \(6\) & \(3^{2}\) & \(6\) & total\\
		\midrule
		7 & $0$ & $1/8$ & $1/4$ & $0$ & $0$ & $1/24$ & $0$ & $1/4$ & $0$ & $1/3$ & $0$ & $0$ & $1$ \\
11 & $1/12$ & $0$ & $3/4$ & $5/12$ & $0$ & $1/12$ & $0$ & $0$ & $0$ & $1/6$ & $1/6$ & $1/3$ & $2$ \\
13 & $0$ & $0$ & $9/8$ & $1/4$ & $1/24$ & $1/6$ & $3/4$ & $0$ & $1/3$ & $1/3$ & $0$ & $0$ & $3$ \\
17 & $0$ & $0$ & $9/4$ & $1$ & $1/12$ & $1/6$ & $1/2$ & $0$ & $1/6$ & $1/3$ & $0$ & $1/2$ & $5$ \\
19 & $1/4$ & $1/2$ & $9/4$ & $3/4$ & $0$ & $5/12$ & $0$ & $1/2$ & $0$ & $1/3$ & $0$ & $0$ & $5$ \\
23 & $1/4$ & $15/8$ & $3$ & $11/12$ & $0$ & $3/8$ & $0$ & $1/4$ & $0$ & $1/2$ & $1/2$ & $1/3$ & $8$ \\
29 & $0$ & $0$ & $51/8$ & $35/12$ & $1/8$ & $2/3$ & $1/4$ & $0$ & $1/2$ & $1/3$ & $0$ & $5/6$ & $12$ \\
31 & $3/4$ & $27/8$ & $9/2$ & $5/4$ & $0$ & $25/24$ & $0$ & $3/4$ & $0$ & $1/3$ & $0$ & $0$ & $12$ \\
37 & $0$ & $0$ & $33/8$ & $3/4$ & $5/24$ & $4/3$ & $5/4$ & $0$ & $2/3$ & $2/3$ & $0$ & $0$ & $9$ \\
41 & $0$ & $0$ & $25/2$ & $35/6$ & $1/3$ & $7/6$ & $1/2$ & $0$ & $1/6$ & $1/3$ & $0$ & $7/6$ & $22$ \\
43 & $3/4$ & $2$ & $27/4$ & $7/4$ & $0$ & $7/4$ & $0$ & $1$ & $0$ & $1$ & $0$ & $0$ & $15$ \\
47 & $5/3$ & $65/8$ & $45/4$ & $23/6$ & $0$ & $11/8$ & $0$ & $3/4$ & $0$ & $1/2$ & $5/6$ & $2/3$ & $29$ \\
53 & $0$ & $0$ & $123/8$ & $65/12$ & $7/24$ & $11/6$ & $3/4$ & $0$ & $5/6$ & $2/3$ & $0$ & $5/6$ & $26$ \\
59 & $5/2$ & $9/2$ & $45/2$ & $29/3$ & $0$ & $7/3$ & $0$ & $1/2$ & $0$ & $1/6$ & $1/2$ & $4/3$ & $44$ \\
61 & $0$ & $0$ & $147/8$ & $25/4$ & $11/24$ & $10/3$ & $5/4$ & $0$ & $2/3$ & $2/3$ & $0$ & $0$ & $31$ \\
67 & $5/4$ & $7/2$ & $45/4$ & $11/4$ & $0$ & $15/4$ & $0$ & $3/2$ & $0$ & $1$ & $0$ & $0$ & $25$ \\
71 & $49/12$ & $147/8$ & $49/2$ & $35/4$ & $0$ & $27/8$ & $0$ & $1/4$ & $0$ & $1/2$ & $7/6$ & $1$ & $62$ \\
73 & $0$ & $0$ & $65/4$ & $9/2$ & $11/12$ & $4$ & $2$ & $0$ & $1/3$ & $1$ & $0$ & $0$ & $29$ \\
79 & $5$ & $145/8$ & $85/4$ & $13/2$ & $0$ & $43/8$ & $0$ & $3/4$ & $0$ & $1$ & $0$ & $0$ & $58$ \\
83 & $4$ & $9$ & $36$ & $41/3$ & $0$ & $11/3$ & $0$ & $1$ & $0$ & $5/6$ & $1/2$ & $4/3$ & $70$ \\
89 & $0$ & $0$ & $195/4$ & $121/6$ & $13/12$ & $14/3$ & $1$ & $0$ & $1/6$ & $1/3$ & $0$ & $11/6$ & $78$ \\
97 & $0$ & $0$ & $89/4$ & $6$ & $17/12$ & $13/2$ & $5/2$ & $0$ & $1/3$ & $1$ & $0$ & $0$ & $40$ \\
101 & $0$ & $0$ & $511/8$ & $325/12$ & $19/24$ & $35/6$ & $3/4$ & $0$ & $5/6$ & $2/3$ & $0$ & $13/6$ & $102$ \\
103 & $15/2$ & $205/8$ & $115/4$ & $17/2$ & $0$ & $59/8$ & $0$ & $9/4$ & $0$ & $1$ & $0$ & $0$ & $81$ \\
107 & $11/2$ & $27/2$ & $99/2$ & $53/3$ & $0$ & $6$ & $0$ & $1/2$ & $0$ & $3/2$ & $1/2$ & $4/3$ & $96$ \\
109 & $0$ & $0$ & $291/8$ & $45/4$ & $9/8$ & $9$ & $5/4$ & $0$ & $1$ & $1$ & $0$ & $0$ & $61$ \\
113 & $0$ & $0$ & $97/2$ & $49/3$ & $3/2$ & $20/3$ & $1$ & $0$ & $1/2$ & $4/3$ & $0$ & $7/6$ & $77$ \\
127 & $10$ & $265/8$ & $145/4$ & $21/2$ & $0$ & $257/24$ & $0$ & $7/4$ & $0$ & $5/3$ & $0$ & $0$ & $104$ \\
131 & $125/12$ & $45/2$ & $375/4$ & $455/12$ & $0$ & $107/12$ & $0$ & $1/2$ & $0$ & $5/6$ & $5/6$ & $7/3$ & $178$ \\
137 & $0$ & $0$ & $121/2$ & $119/6$ & $2$ & $55/6$ & $3/2$ & $0$ & $1/2$ & $4/3$ & $0$ & $7/6$ & $96$ \\
139 & $15/2$ & $39/2$ & $135/2$ & $23$ & $0$ & $13$ & $0$ & $3/2$ & $0$ & $1$ & $0$ & $0$ & $133$ \\
149 & $0$ & $0$ & $847/8$ & $481/12$ & $35/24$ & $23/2$ & $3/4$ & $0$ & $7/6$ & $1$ & $0$ & $13/6$ & $164$ \\
151 & $63/4$ & $427/8$ & $119/2$ & $75/4$ & $0$ & $123/8$ & $0$ & $5/4$ & $0$ & $1$ & $0$ & $0$ & $165$ \\
157 & $0$ & $0$ & $435/8$ & $65/4$ & $43/24$ & $44/3$ & $13/4$ & $0$ & $4/3$ & $4/3$ & $0$ & $0$ & $93$ \\
163 & $13/4$ & $19/2$ & $117/4$ & $27/4$ & $0$ & $193/12$ & $0$ & $5/2$ & $0$ & $5/3$ & $0$ & $0$ & $69$ \\
167 & $275/12$ & $671/8$ & $99$ & $415/12$ & $0$ & $281/24$ & $0$ & $5/4$ & $0$ & $7/6$ & $11/6$ & $5/3$ & $258$ \\
173 & $0$ & $0$ & $1015/8$ & $559/12$ & $13/8$ & $38/3$ & $5/4$ & $0$ & $3/2$ & $4/3$ & $0$ & $13/6$ & $194$ \\
		\bottomrule
	\end{tabular}
\end{table}

\section{Rational kernels with decomposed quotient}\label{sec:kernels}

Let \(C/\Fp\) be a curve of genus two and \(K\subset J[2]\) a rational maximal
isotropic subgroup, corresponding to a \(\sigma\)-stable partition of \(W\) into
three pairs. The quotient \(J/K\), with the polarisation induced by
\(2\theta\), is a principally polarised abelian surface over \(\Fp\). It is the
Jacobian of a curve given by the Richelot construction, unless the three pairs
are the orbits of an involution of \(\PP^{1}\), which is then unique and defined
over \(\Fp\); in that case \(J/K\) is a product of two elliptic curves with the
product polarisation \cite[Prop.~4.3]{KatsuraTakashima2020}, and we call \(K\)
\emph{decomposed}. A decomposed \(K\) is \emph{split} if the two fixed points of
the involution are rational, and \emph{nonsplit} if they are conjugate over
\(\F_{p^{2}}\). In the split case the two elliptic factors are defined over
\(\Fp\); in the nonsplit case they are conjugate over \(\F_{p^{2}}\) and the
quotient is a Weil restriction \cite{DangNguyen2026}. We also record the
permutation that Frobenius induces on the three pairs, of type \(111\), \(12\)
or \(3\).

The mass of a set of pairs \((C,K)\) is the sum of \(1/\#\Aut_{\Fp}(C,K)\) over
\(\Fp\)-isomorphism classes; as in Lemma~\ref{lem:massdiv} it is the sum of
\(1/\#S\) over the \(\PGL_2(\Fp)\)-orbits of pairs of a branch divisor and a
partition, with \(S\) the stabiliser.

\begin{theorem}\label{thm:kernel}
	Let \(p\ge7\), \(h=h(-p)\) for \(p\equiv3\pmod4\), and \(h'=h(-4p)\). The
	masses of the pairs \((C,K)\) with \(C/\Fp\) superspecial of genus two and
	\(K\) rational and decomposed are given by Table~\textup{\ref{tab:kernel}};
	cells that are not listed are empty. In particular, on the curves with
	\(\pi^{2}=\eps p\):
	\begin{enumerate}[label=\textup{(\alph*)},leftmargin=*]
		\item for each sign \(\eps\), the nonsplit decomposed kernels have mass
		\((p-1)/12-H_p/8\);
		\item for \(\eps=-1\), the split decomposed kernels have mass
		\((6h^{2}+2h'^{2})/8-(h+h')/4\) if \(p\equiv3\pmod4\) and
		\((h'^{2}-h')/4\) if \(p\equiv1\pmod4\);
		\item for \(\eps=+1\) there is no split decomposed kernel.
	\end{enumerate}
	On the curves with \(\pi^{2}\ne\pm p\) the decomposed kernels are nonsplit and
	have total mass \(1/3\) if \(p\equiv2\pmod3\) and \(0\) otherwise.
\end{theorem}

\begin{proposition}\label{prop:nonjac}
	Let \(p\ge5\). The products of two supersingular elliptic curves over \(\Fp\),
	with the product polarisation, have mass \(H_p^{2}/8\), and \(\pi^{2}=-p\) on
	each of them. For each sign \(\eps\), the Weil restrictions to \(\Fp\) of
	supersingular elliptic curves over \(\F_{p^{2}}\) with \(\pi^{2}=\eps p\) have
	mass \((p-1)/48\). The remaining Weil restrictions of supersingular elliptic
	curves have Weil polynomial \(T^{4}\pm pT^{2}+p^{2}\), with mass \(1/6\) for
	each sign when \(p\equiv2\pmod3\), or \(T^{4}+p^{2}\), with mass \(1/4\) when
	\(p\equiv3\pmod4\).
\end{proposition}

Here the mass is the sum of \(1/\#\Aut_{\Fp}(A,\theta)\) over the
\(\Fp\)-isomorphism classes of principally polarised surfaces. In particular the
two masses used in Proposition~\ref{prop:ibukiyama} are proved.

\begin{table}[ht]
	\centering
	\caption{Masses of pairs \((C,K)\) with \(K\) rational and
		decomposed. The column \(\mu\) is the type of the permutation induced by
		Frobenius on the three pairs.}\label{tab:kernel}
	\renewcommand{\arraystretch}{1.25}
	\begin{tabular}{@{}llllll@{}}
		\toprule
		\(p\) & \(P\) & \(\lambda\) & \(\mu\) & class & mass\\
		\midrule
		\(1\bmod4\) & \((T^{2}-p)^{2}\) & \(1^{6}\) & \(111\) & nonsplit & \((p-1)/48\)\\
		& \((T^{2}+p)^{2}\) & \(2^{3}\) & \(111\) & nonsplit & \((p-1)/48\)\\
		\(3\bmod4\) & \((T^{2}-p)^{2}\) & \(2^{3}\) & \(111\) & nonsplit & \((p-1)/48-h/8\)\\
		& \((T^{2}+p)^{2}\) & \(1^{6}\) & \(111\) & nonsplit & \((p-1)/48-h/8\)\\
		& \((T^{2}+p)^{2}\) & \(1^{6}\) & \(111\) & split & \(h(h'-1)/8\)\\
		& \((T^{2}+p)^{2}\) & \(2^{3}\) & \(111\) & split & \(h(h'-1)/8\)\\
		\(7\bmod8\) & \((T^{2}+p)^{2}\) & \(1^{4}2\) & \(111\) & split & \(h^{2}/4\)\\
		& \((T^{2}+p)^{2}\) & \(1^{2}2^{2}\) & \(111\) & split & \(h^{2}/4\)\\
		all & \((T^{2}-p)^{2}\) & \(2^{3}\) & \(12\) & nonsplit & \((p-1)/16-h'/8\)\\
		& \((T^{2}+p)^{2}\) & \(1^{2}2^{2}\) & \(12\) & nonsplit & \((p-1)/16-h'/8\)\\
		& \((T^{2}+p)^{2}\) & \(1^{2}2^{2}\) & \(12\) & split & \((h'^{2}-h')/8\)\\
		& \((T^{2}+p)^{2}\) & \(2^{3}\) & \(12\) & split & \((h'^{2}-h')/8\)\\
		\(11\bmod12\) & \(T^{4}-pT^{2}+p^{2}\) & \(3^{2}\) & \(3\) & nonsplit & \(1/6\)\\
		& \(T^{4}+pT^{2}+p^{2}\) & \(6\) & \(3\) & nonsplit & \(1/6\)\\
		\(5\bmod12\) & \(T^{4}-pT^{2}+p^{2}\) & \(6\) & \(3\) & nonsplit & \(1/6\)\\
		& \(T^{4}+pT^{2}+p^{2}\) & \(3^{2}\) & \(3\) & nonsplit & \(1/6\)\\
		\bottomrule
	\end{tabular}
\end{table}

The proof is in Section~\ref{sec:kernelproof}. The equality of the masses in
(a) for the two signs, and the symmetry between the cells of
Table~\ref{tab:kernel} with \(P=(T^{2}-p)^{2}\) and with \(P=(T^{2}+p)^{2}\),
come from the quadratic twist over \(\F_{p^{2}}\) of the elliptic curve whose
Weil restriction is the quotient.

\begin{proposition}\label{prop:kernelcheck}
	For the \(37\) primes \(7\le p\le173\) the enumeration of
	Section~\textup{\ref{sec:comp}} gives the masses of
	Theorem~\textup{\ref{thm:kernel}}, cell by cell, and the number of rational
	maximal isotropic subgroups is \(15,3,3,7,1,1,3,1\) on every curve of type
	\(1^{6},1^{4}2,1^{2}2^{2},2^{3},1^{2}4,2\,4,3^{2},6\), in accordance with
	Theorems~\textup{\ref{thm:local}} and~\textup{\ref{thm:class}}.
\end{proposition}

This is an independent check of the theorem: the enumeration works on the
curves and does not use elliptic curves.

Table~\ref{tab:prop} gives, for each prime, the mass of the pairs \((C,K)\) with
\(K\) rational and decomposed and the mass of all pairs \((C,K)\) with \(K\)
rational, on the curves with \(\pi^{2}=-p\), on those with \(\pi^{2}=p\), and on
all superspecial curves over \(\Fp\). On all curves the quotient of the two
masses, multiplied by \(\sqrt p\), lies between \(0.80\) and \(1.63\) for
\(7\le p\le173\), and multiplied by \(p\) it lies between \(2.9\) and \(18.8\).
Theorems~\ref{thm:class}, \ref{thm:vertex} and~\ref{thm:kernel} give both
masses in closed form, and the order of magnitude of the quotient follows.

\begin{corollary}\label{cor:prop}
	For a prime \(p\ge7\) let \(\rho(p)\) be the mass of the pairs \((C,K)\), with
	\(C/\Fp\) superspecial of genus two and \(K\) a rational maximal isotropic
	subgroup of \(J[2]\) with decomposed quotient, divided by the mass of all such
	pairs with \(K\) rational. Let \(\Delta=-p\) if \(p\equiv3\pmod4\) and
	\(\Delta=-4p\) if \(p\equiv1\pmod4\), and \(L_p=L(1,\chi_\Delta)\) for the
	Kronecker character \(\chi_\Delta\).
	\begin{enumerate}[label=\textup{(\alph*)},leftmargin=*]
		\item \(\rho(p)=\bigl(D_\pm(p)+D_0(p)\bigr)/\bigl(T_\pm(p)+T_0(p)\bigr)\),
		where \(D_0(p)=[p\equiv2\bmod3]/3\), \(T_0(p)\ge0\) is the sum over the last
		four rows of Table~\textup{\ref{tab:class}} of the number of rational
		subgroups times the mass of Theorem~\textup{\ref{thm:vertex}(c)--(e)}, and,
		with \(h=h(-p)\), \(h'=h(-4p)\), \(z_p=\zeta_{\Q(\sqrt p)}(-1)\),
		\[
		\begin{array}{lll}
		p\equiv7\ (8): & D_\pm=h^{2}-h+\tfrac{p-1}{6}, &
		T_\pm=\tfrac56(p-1)h-3h^{2}+\tfrac74z_p-\tfrac{7(p-1)}{24},\\[2pt]
		p\equiv3\ (8): & D_\pm=3h^{2}-2h+\tfrac{p-1}{6}, &
		T_\pm=\bigl(\tfrac32p-\tfrac98\bigr)h-\tfrac{15}{2}h^{2}+\tfrac74z_p-\tfrac{7(p-1)}{24},\\[2pt]
		p\equiv1\ (4): & D_\pm=\tfrac14h'^{2}-\tfrac12h'+\tfrac{p-1}{6}, &
		T_\pm=\bigl(\tfrac p3-\tfrac7{48}\bigr)h'-\tfrac38h'^{2}
		+\bigl(\tfrac{71}{4}-\tfrac72\leg2p\bigr)z_p-\tfrac{7(p-1)}{24}.
		\end{array}
		\]
		\item There are absolute and effectively computable constants \(c,C>0\) such
		that \(c\,(1+L_p)\le\rho(p)\sqrt p\le C\,(1+L_p)\) for all \(p\ge7\). In
		particular \(c\,p^{-1/2}\le\rho(p)\le C'\,p^{-1/2}\log p\).
	\end{enumerate}
\end{corollary}

\begin{proof}
	(a) The mass of all pairs is the sum over the rows of
	Table~\ref{tab:class} of the number of rational subgroups
	(Theorem~\ref{thm:class}) times the mass of the curves
	(Theorem~\ref{thm:vertex}), and the decomposed mass is given by
	Theorem~\ref{thm:kernel}; the displayed expressions collect the rows with
	\(\pi^{2}=\pm p\), using \(h'=(2-\leg2p)h\) for \(p\equiv3\pmod4\)
	\cite[Thm.~7.24]{Cox1989}.

	(b) Put \(x=L_p\). By the class number formula, \(h(\Delta)=\sqrt{|\Delta|}\,x/\pi\)
	\cite[Ch.~6]{Davenport2000}, and by partial summation \(x\le\log|\Delta|+2\).
	By the functional equation of \(L(s,\chi_F)\) for the character \(\chi_F\) of
	\(F=\Q(\sqrt p)\) \cite[Ch.~9]{Davenport2000},
	\(z_p=\zeta(-1)L(-1,\chi_F)=d_F^{3/2}L(2,\chi_F)/(24\pi^{2})\) with
	\(d_F\in\{p,4p\}\), and \(\zeta(4)/\zeta(2)\le L(2,\chi_F)\le\zeta(2)\); hence
	\(p^{3/2}/360\le z_p\le p^{3/2}/18\). Every class number \(h(D)\) with
	\(|D|\le12p\) occurring in Theorem~\ref{thm:vertex}(c)--(e) is
	\(O(\sqrt p\log p)\) by the same two facts and \cite[Thm.~7.24]{Cox1989}, so
	\(0\le T_0(p)=O(\sqrt p\log p)\). Since \(x=O(\log p)\), the quantities
	\(h\), \(h'\), \(h^{2}+p\) and \(h'^{2}+p\) are \(o(p^{3/2})\), and the
	displayed formulas give, with \(\kappa=z_p/p^{3/2}\in[1/360,1/18]\),
	\[
	D_\pm+D_0=p\bigl(a\,x^{2}+\tfrac16\bigr)\bigl(1+o(1)\bigr),\qquad
	T_\pm+T_0=p^{3/2}\bigl(b\,x+e\,\kappa\bigr)\bigl(1+o(1)\bigr),
	\]
	where \(a,b,e>0\) depend only on \(p\bmod8\) (for instance \(a=1/\pi^{2}\),
	\(b=5/(6\pi)\), \(e=7/4\) for \(p\equiv7\pmod8\)), and the \(o(1)\) are
	bounded by effective functions of \(p\) alone. As
	\((1+x)/2\le(1+x^{2})/(1+x)\le1+x\) for \(x\ge0\), and
	\(\min(a,\tfrac16)(1+x^{2})\le ax^{2}+\tfrac16\le\max(a,\tfrac16)(1+x^{2})\) and
	\(\min(b,\tfrac{e}{360})(1+x)\le bx+e\kappa\le\max(b,\tfrac{e}{18})(1+x)\), this
	gives the inequalities for \(p\) larger than an effective bound. For the finitely many remaining
	primes \(\rho(p)>0\), because \(D_\pm>0\): \(h^{2}-h\ge0\),
	\(3h^{2}-2h>0\) and \(h'^{2}/4-h'/2\ge-1/4\), while \((p-1)/6\ge1\). The last
	assertion follows from \(x>0\) and \(x\le\log(4p)+2\).
\end{proof}

The constants in (b) are not optimal; for \(7\le p\le173\),
\(\rho(p)\sqrt p\) lies between \(0.80\) and \(1.63\), and the formula of (a)
gives values between \(0.79\) and \(1.57\) for the \(25\) primes in the four
ranges \(1000<p<1060\), \(3000<p<3060\), \(10000<p<10060\) and
\(30000<p<30060\). Under the generalised Riemann hypothesis
\(L_p=O(\log\log p)\) \cite{Littlewood1928}, so the factor \(\log p\) can
then be replaced by \(\log\log p\).

Three proportions should be kept apart. Over \(\bar\F_p\), Katsura and
Takashima count the Richelot isogenies out of superspecial Jacobians up to
isomorphism, \((p-1)(p+2)(p+7)/192\) up to bounded terms, of which
\((p-1)(p+3)/48\) up to bounded terms are decomposed
\cite[Thm.~6.2]{KatsuraTakashima2020}; Florit and Smith give \(10/p\) for the
proportion of products among the vertices of the superspecial graph and
\(5/p\) for the stationary proportion of a random walk
\cite[\S6.2]{FloritSmith2022auto}. Among all pairs of a genus-two curve over a
finite field \(\F_q\) and a rational kernel, weighted by automorphisms, the
decomposed ones have proportion \(1/q\) \cite{DangNguyen2026}. The numbers of
Table~\ref{tab:prop} concern the \(\Fp\)-rational superspecial locus with its
rational kernels, and they are larger than both.

\begin{table}[tp]
	\centering
	\scriptsize
	\caption{Masses of pairs \((C,K)\) with \(K\) rational: decomposed (dec.) and
		all, their quotient \(\rho\), for \(\pi^{2}=-p\), for \(\pi^{2}=p\), and
		on all superspecial curves of genus two over \(\Fp\). Decimals are
		rounded.}\label{tab:prop}
	\setlength{\tabcolsep}{3.5pt}
	\begin{tabular}{@{}r|rrrr|rrr|rrrr@{}}
		\toprule
		& \multicolumn{4}{c|}{\(\pi^{2}=-p\)} & \multicolumn{3}{c|}{\(\pi^{2}=p\)}
		& \multicolumn{4}{c}{all}\\
		\(p\) & dec. & all & \(\rho\) & \(\rho\sqrt p\) & dec. & all & \(\rho\)
		& dec. & all & \(\rho\) & \(\rho\sqrt p\)\\
		\midrule
		7 & $3/4$ & $9/8$ & 0.667 & 1.76 & $1/4$ & $7/24$ & 0.857 & $1$ & $2$ & 0.500 & 1.32 \\
11 & $7/3$ & $77/12$ & 0.364 & 1.21 & $1/3$ & $7/12$ & 0.571 & $3$ & $8$ & 0.375 & 1.24 \\
13 & $5/4$ & $41/8$ & 0.244 & 0.88 & $3/4$ & $43/24$ & 0.419 & $2$ & $9$ & 0.222 & 0.80 \\
17 & $23/6$ & $55/4$ & 0.279 & 1.15 & $5/6$ & $29/12$ & 0.345 & $5$ & $18$ & 0.278 & 1.15 \\
19 & $3$ & $69/4$ & 0.174 & 0.76 & $1$ & $35/12$ & 0.343 & $4$ & $21$ & 0.190 & 0.83 \\
23 & $103/12$ & $595/24$ & 0.346 & 1.66 & $13/12$ & $21/8$ & 0.413 & $10$ & $30$ & 0.333 & 1.60 \\
29 & $109/12$ & $949/24$ & 0.230 & 1.24 & $19/12$ & $157/24$ & 0.242 & $11$ & $49$ & 0.224 & 1.21 \\
31 & $37/4$ & $349/8$ & 0.212 & 1.18 & $7/4$ & $175/24$ & 0.240 & $11$ & $52$ & 0.212 & 1.18 \\
37 & $13/4$ & $141/8$ & 0.184 & 1.12 & $11/4$ & $299/24$ & 0.221 & $6$ & $34$ & 0.176 & 1.07 \\
41 & $49/3$ & $235/3$ & 0.209 & 1.34 & $7/3$ & $79/6$ & 0.177 & $19$ & $94$ & 0.202 & 1.29 \\
43 & $5$ & $199/4$ & 0.101 & 0.66 & $3$ & $49/4$ & 0.245 & $8$ & $64$ & 0.125 & 0.82 \\
47 & $301/12$ & $2639/24$ & 0.228 & 1.56 & $31/12$ & $77/8$ & 0.268 & $28$ & $124$ & 0.226 & 1.55 \\
53 & $133/12$ & $2017/24$ & 0.132 & 0.96 & $43/12$ & $413/24$ & 0.208 & $15$ & $106$ & 0.142 & 1.03 \\
59 & $82/3$ & $1117/6$ & 0.147 & 1.13 & $10/3$ & $49/3$ & 0.204 & $31$ & $206$ & 0.150 & 1.16 \\
61 & $47/4$ & $791/8$ & 0.119 & 0.93 & $17/4$ & $725/24$ & 0.141 & $16$ & $133$ & 0.120 & 0.94 \\
67 & $7$ & $329/4$ & 0.085 & 0.70 & $5$ & $105/4$ & 0.190 & $12$ & $111$ & 0.108 & 0.88 \\
71 & $595/12$ & $2009/8$ & 0.197 & 1.66 & $49/12$ & $189/8$ & 0.173 & $54$ & $280$ & 0.193 & 1.63 \\
73 & $17/2$ & $321/4$ & 0.106 & 0.90 & $11/2$ & $167/4$ & 0.132 & $14$ & $126$ & 0.111 & 0.95 \\
79 & $111/4$ & $1909/8$ & 0.116 & 1.03 & $21/4$ & $301/8$ & 0.140 & $33$ & $278$ & 0.119 & 1.06 \\
83 & $88/3$ & $872/3$ & 0.101 & 0.92 & $16/3$ & $77/3$ & 0.208 & $35$ & $321$ & 0.109 & 0.99 \\
89 & $233/6$ & $3449/12$ & 0.135 & 1.27 & $35/6$ & $587/12$ & 0.119 & $45$ & $340$ & 0.132 & 1.25 \\
97 & $21/2$ & $435/4$ & 0.097 & 0.95 & $15/2$ & $267/4$ & 0.112 & $18$ & $180$ & 0.100 & 0.98 \\
101 & $625/12$ & $9149/24$ & 0.137 & 1.37 & $79/12$ & $1265/24$ & 0.125 & $59$ & $440$ & 0.134 & 1.35 \\
103 & $119/4$ & $2681/8$ & 0.089 & 0.90 & $29/4$ & $413/8$ & 0.140 & $37$ & $390$ & 0.095 & 0.96 \\
107 & $94/3$ & $2371/6$ & 0.079 & 0.82 & $22/3$ & $42$ & 0.175 & $39$ & $442$ & 0.088 & 0.91 \\
109 & $63/4$ & $1503/8$ & 0.084 & 0.88 & $33/4$ & $639/8$ & 0.103 & $24$ & $273$ & 0.088 & 0.92 \\
113 & $67/3$ & $1559/6$ & 0.086 & 0.91 & $25/3$ & $415/6$ & 0.120 & $31$ & $334$ & 0.093 & 0.99 \\
127 & $127/4$ & $3453/8$ & 0.074 & 0.83 & $37/4$ & $1799/24$ & 0.123 & $41$ & $510$ & 0.080 & 0.91 \\
131 & $235/3$ & $9245/12$ & 0.102 & 1.16 & $25/3$ & $749/12$ & 0.134 & $87$ & $839$ & 0.104 & 1.19 \\
137 & $73/3$ & $961/3$ & 0.076 & 0.89 & $31/3$ & $565/6$ & 0.110 & $35$ & $420$ & 0.083 & 0.98 \\
139 & $34$ & $1069/2$ & 0.064 & 0.75 & $10$ & $91$ & 0.110 & $44$ & $628$ & 0.070 & 0.83 \\
149 & $673/12$ & $14357/24$ & 0.094 & 1.14 & $127/12$ & $819/8$ & 0.103 & $67$ & $708$ & 0.095 & 1.16 \\
151 & $225/4$ & $5649/8$ & 0.080 & 0.98 & $43/4$ & $861/8$ & 0.100 & $67$ & $816$ & 0.082 & 1.01 \\
157 & $79/4$ & $2215/8$ & 0.071 & 0.89 & $49/4$ & $3109/24$ & 0.095 & $32$ & $415$ & 0.077 & 0.97 \\
163 & $15$ & $849/4$ & 0.071 & 0.90 & $13$ & $1351/12$ & 0.115 & $28$ & $329$ & 0.085 & 1.09 \\
167 & $1519/12$ & $27227/24$ & 0.112 & 1.44 & $133/12$ & $1967/24$ & 0.135 & $138$ & $1226$ & 0.113 & 1.45 \\
173 & $697/12$ & $16961/24$ & 0.082 & 1.08 & $151/12$ & $2713/24$ & 0.111 & $71$ & $829$ & 0.086 & 1.13 \\
		\bottomrule
	\end{tabular}
\end{table}

\section{Proof of Theorem~\ref{thm:kernel}}\label{sec:kernelproof}

The masses are computed in groupoids. A groupoid is \emph{finite} if it has finitely many isomorphism classes and
finite automorphism groups; its mass is the sum of \(1/\#\Aut(x)\) over the
isomorphism classes.

\begin{lemma}\label{lem:fibre}
	Let \(F\colon\mathcal Y\to\mathcal X\) be a faithful functor between finite
	groupoids. For an object \(x\) of \(\mathcal X\) let \(\Fib(x)\) be the set of
	pairs \((y,f)\), with \(f\colon F(y)\to x\) an isomorphism, modulo
	\((y,f)\sim(y',f')\) if \(f'F(g)=f\) for some \(g\colon y\to y'\). Then
	\(\Mass(\mathcal Y)=\sum_{[x]}\#\Fib(x)/\#\Aut(x)\).
\end{lemma}

\begin{proof}
	The group \(\Aut(x)\) acts on \(\Fib(x)\) by \(a\cdot(y,f)=(y,af)\). The orbits
	correspond to the classes \([y]\) with \(F(y)\cong x\), and the stabiliser of
	\((y,f)\) is \(fF(\Aut(y))f^{-1}\), of order \(\#\Aut(y)\) because \(F\) is
	faithful. Hence \(\#\Fib(x)=\sum_{[y]}\#\Aut(x)/\#\Aut(y)\).
\end{proof}

When \(F\) forgets a structure, \(\Fib(x)\) is the set of structures on \(x\).

For a morphism \(g\) between varieties over \(\Fp\) we write
\(g^{c}=\Frob\circ g\circ\Frob^{-1}\). For supersingular \(E/\Fp\) let
\(g_E\in\Aut(E[2])\cong S_3\) be the action of \(\pi_E\). Let \(N_1\) be the
number of \(\Fp\)-isomorphism classes of supersingular elliptic curves over
\(\Fp\) with \(g_E=1\) (type~I), and \(N_2=H_p-N_1\) the number of the others
(type~II).

\begin{lemma}\label{lem:ssp}
	Let \(p\ge5\) and let \(E/\Fp\) be a supersingular elliptic curve.
	\begin{enumerate}[label=\textup{(\alph*)},leftmargin=*]
		\item \(\pi_E^{2}=-p\), \(\Aut_{\Fp}(E)=\{\pm1\}\), and \(g_E\) is the
		identity or a transposition.
		\item \(N_1=h(-p)\) if \(p\equiv3\pmod4\), \(N_1=0\) if \(p\equiv1\pmod4\),
		and \(N_2=h(-4p)\).
		\item For supersingular \(E_1,E_2\) over \(\Fp\), an isomorphism
		\(\alpha\colon E_1\to E_2\) over \(\kbar\) satisfies
		\(\pi_{E_2}\alpha=\alpha^{c}\pi_{E_1}\). There is exactly one class
		\(E_1^{\mathrm{tw}}\) of curves \(E_2\) admitting \(\alpha\) with
		\(\alpha^{c}=-\alpha\); such \(\alpha\) form a torsor under \(\{\pm1\}\).
		\item For an elliptic curve \(E/\kbar\), the kernel of
		\(\Aut(E)\to\Aut(E[2])\) is \(\{\pm1\}\).
	\end{enumerate}
\end{lemma}

\begin{proof}
	(a) The supersingular isogeny class over \(\Fp\), \(p\ge5\), has trace zero
	\cite[Thm.~4.1]{Waterhouse1969}. An \(\Fp\)-automorphism commutes with
	\(\pi_E\), hence lies in \(\Q(\pi_E)=\Q(\sqrt{-p})\), whose roots of unity are
	\(\pm1\). As \(\#E(\Fp)=p+1\) is even, \(g_E\) fixes a point of order two.
	(b) \(g_E=1\) if and only if \(\pi_E-1=2\gamma\) with
	\(\gamma\in\End_{\Fp}(E)\) \cite[III.4.11]{Silverman2009}, that is, if and only
	if \(\End_{\Fp}(E)\) is the maximal order of \(\Q(\sqrt{-p})\), which requires
	\(p\equiv3\pmod4\); otherwise \(\End_{\Fp}(E)=\Z[\sqrt{-p}]\). Both orders
	occur \cite[Thm.~4.2]{Waterhouse1969}, and over the prime field the classes
	with a given endomorphism ring \(R\) form a principal homogeneous space under
	the class group of \(R\) \cite[Thm.~4.5]{Waterhouse1969}; see also
	\cite[\S2]{DelfsGalbraith2016}.
	(c) The first assertion is the naturality of Frobenius. The quadratic twist
	\(dy^{2}=f(x)\) of \(y^{2}=f(x)\), with \(d\) a non-square, admits
	\(\alpha=(x,y/\sqrt d)\); if \(\alpha\colon E_1\to E_2\) and
	\(\alpha'\colon E_1\to E_2'\) satisfy \(\alpha^{c}=-\alpha\),
	\(\alpha'^{c}=-\alpha'\), then \(\alpha'\alpha^{-1}\) is defined over \(\Fp\).
	(d) If \(u\in\Aut(E)\) is trivial on \(E[2]\), then \(u-1=2\gamma\) and
	\(\deg(u-1)=2-\operatorname{tr}(u)\) is divisible by four; but it is
	\(2,3,1\) for \(u\) of order \(4,3,6\).
\end{proof}

\begin{lemma}\label{lem:ssp2}
	Let \(p\ge5\) and let \(E/\F_{p^{2}}\) be a supersingular elliptic curve with
	Frobenius endomorphism \(\pi_E\). Then \(\pi_E=p\,u_E\) with
	\(u_E\in\Aut(E_{\kbar})\), and all automorphisms of \(E_{\kbar}\) are defined
	over \(\F_{p^{2}}\). For every supersingular \(j\) and every element \(u\) of
	the cyclic group \(\Aut(E_{j})\) there is exactly one
	\(\F_{p^{2}}\)-isomorphism class of curves \(E\) with \(j(E)=j\) and
	\(u_E=u\). For each sign \(\eps\), the supersingular elliptic curves over
	\(\F_{p^{2}}\) with \(\pi_E=\eps p\) have mass \((p-1)/24\).
\end{lemma}

\begin{proof}
	The first assertion is proved as Lemma~\ref{lem:frob2}; the automorphism group
	is commutative, so its elements commute with \(\pi_E\). The forms of a curve
	\(E_0\) are classified by \(\xi\in\Aut(E_{0,\kbar})\)
	\cite[X.5.3]{Silverman2009}: \(E_\xi\) comes with
	\(\vartheta\colon E_0\to E_\xi\) such that
	\(\vartheta^{-1}\vartheta^{cc}=\xi\), and
	\(\pi_{E_\xi}\vartheta=\vartheta^{cc}\pi_{E_0}\) gives
	\(u_{E_\xi}=\vartheta\,\xi u_{E_0}\vartheta^{-1}\). The mass is
	\(\sum_j1/\#\Aut(E_j)=(p-1)/24\) by
	\cite[42.3.8, (42.3.9)]{Voight2021}.
\end{proof}

\begin{lemma}\label{lem:prod}
	Let \(E_1,E_2,F_1,F_2\) be elliptic curves over \(\kbar\) and
	\(g=(g_{ij})\colon E_1\times E_2\to F_1\times F_2\) a homomorphism pulling the
	product polarisation back to \(n\) times the product polarisation. Then the
	rows and columns of \((\deg g_{ij})\) sum to \(n\). In particular a
	principally polarised surface which is a polarised product contains exactly
	two elliptic curves \(B_1,B_2\) with \(A=B_1\times B_2\) as polarised
	varieties, and every isomorphism between polarised products is of the form
	\(a_1\times a_2\) or \((a_1\times a_2)\circ\mathrm{swap}\).
\end{lemma}

\begin{proof}
	The hypothesis is \(\sum_i\hat g_{ij}g_{ik}=n\,\delta_{jk}\), that is
	\(\hat g g=n\) for \(\hat g=(\hat g_{ij})^{t}\); then \(g\hat g=n\) as well.
	For \(n=1\) the degree matrix is a permutation matrix.
\end{proof}

Let \((A,\nu)\) be a principally polarised surface over \(\Fp\) which is
geometrically a polarised product. Frobenius permutes \(\{B_1,B_2\}\). If it
fixes them, \(A=B_1\times B_2\) over \(\Fp\) (\emph{split}). Otherwise
(\emph{nonsplit}) \(B=B_1\) and \(B'=\Frob(B)\) are defined over
\(\F_{p^{2}}\); we put \(\varphi=\pi_A|_B\colon B\to B'\) and
\(\varphi'=\pi_A|_{B'}\colon B'\to B\). These are isogenies of degree \(p\) with
\(\varphi'\varphi=\pi_B\), \(\varphi\varphi'=\pi_{B'}\), and
\(\gamma^{c}\varphi=\varphi'\gamma\) for every \(\gamma\colon B\to B'\),
\(\varphi\beta=\beta^{c}\varphi\) for every \(\beta\in\End(B_{\kbar})\). The
nonsplit surfaces are the Weil restrictions, and every elliptic curve over
\(\F_{p^{2}}\) occurs as \(B\)
\cite[Thm.~1.3]{HoweNartRitzenthaler2009}, \cite[Thm.~3.1(3)]{GonzalezGuardiaRotger2005}.

\begin{lemma}\label{lem:glue}
	Let \(p\) be odd, \(A=E_1\times E_2\) over \(\kbar\) with the product
	polarisation.
	\begin{enumerate}[label=\textup{(\alph*)},leftmargin=*]
		\item The maximal isotropic subgroups of \(A[2]\) are the nine products
		\(\langle a\rangle\times\langle b\rangle\) and the six graphs
		\(\Gamma_\psi\) of isomorphisms \(\psi\colon E_1[2]\to E_2[2]\).
		\item \(A/\Gamma_\psi\), with the polarisation induced by twice the product
		polarisation, is a polarised product if \(\psi\) is the restriction of an
		isomorphism \(E_1\to E_2\) \textup{(}\(\psi\) \emph{reducible}\textup{)}, and
		the Jacobian of a curve of genus two otherwise.
		\item If \(E_1,E_2\) are supersingular, \(A/\Gamma_\psi\) is isomorphic to
		\(E_1\times E_2\) as an unpolarised abelian variety.
	\end{enumerate}
\end{lemma}

\begin{proof}
	(a) A subgroup of order four meeting neither factor is a graph, and every
	isomorphism of two-dimensional \(\F_2\)-spaces respects the non-degenerate
	alternating forms; if it contains \((a,0)\ne0\), isotropy forces the form
	\(\langle a\rangle\times\langle b\rangle\).
	(b) If \(\psi=\alpha|_{E_1[2]}\), the endomorphism
	\((x,y)\mapsto(x+\alpha^{-1}y,\alpha x-y)\) has kernel \(\Gamma_\psi\) and
	pulls the product polarisation back to twice itself. Conversely, if
	\(A/\Gamma_\psi\) is a polarised product, Lemma~\ref{lem:prod} with \(n=2\)
	applies to \(g\colon A\to A/\Gamma_\psi\): either \(\ker g\) is a product of
	subgroups of the factors, which a graph is not, or all \(g_{ij}\) are
	isomorphisms and \(\Gamma_\psi\) lies in the graph of
	\(-g_{12}^{-1}g_{11}\). A principally polarised surface over \(\kbar\) is a
	Jacobian or a polarised product \cite[Thm.~1.3]{HoweNartRitzenthaler2009}.
	This is the case \(N=2\) of \cite[Thm.~3]{Kani1997}; see also
	\cite[Prop.~3]{HoweLeprevostPoonen2000}.
	(c) For models over \(\F_{p^{2}}\) with \(\pi=-p\)
	(Lemma~\ref{lem:ssp2}) Tate's theorem \cite[Main Theorem]{Tate1966} gives
	\(\Hom(E_1,E_2)\otimes\Z_2=\Hom(T_2E_1,T_2E_2)\), so \(\psi=\beta|_{E_1[2]}\)
	for some \(\beta\in\Hom(E_1,E_2)\); the automorphism
	\((x,y)\mapsto(x,y-\beta x)\) of \(A\) maps \(\Gamma_\psi\) to
	\(E_1[2]\times0\).
\end{proof}

\begin{lemma}\label{lem:dual}
	Let \(\mathcal G\) be the groupoid of triples \((A,\nu,K)\), with
	\((A,\nu)\) a principally polarised abelian surface over \(\Fp\) and
	\(K\subset A[2]\) a rational maximal isotropic subgroup. The rule
	\(D(A,\nu,K)=(A/K,\nu',A[2]/K)\), where \(\nu'\) is the principal
	polarisation induced by \(2\nu\) \cite[Prop.~1.1]{Kani1997}, is an
	equivalence of \(\mathcal G\) with itself, with \(D\circ D\cong\mathrm{id}\),
	compatible with Frobenius. The functor \((C,K)\mapsto(J,\theta,K)\) from pairs
	of a curve of genus two over \(\Fp\) and a rational maximal isotropic subgroup
	to \(\mathcal G\) is fully faithful, with essential image the triples for
	which \((A,\nu)\) is geometrically a Jacobian.
\end{lemma}

\begin{proof}
	Let \(f\colon A\to A/K\) be the quotient. The polarisation \(\nu'\) with
	\(\hat f\nu' f=2\nu\) exists and is unique over \(\kbar\), hence is defined
	over \(\Fp\); \(f(A[2])\) is isotropic because
	\(e_2^{\nu'}(fx,fy)=e_2^{2\nu}(x,y)\). For \(f'=\nu^{-1}\hat f\nu'\) we
	have \(f'f=2\), \(ff'=2\), \(\ker f'=f(A[2])\),
	\(\hat f'\nu f'=2\nu'\) and \(f'((A/K)[2])=K\); so \(f'\) identifies
	\(D(D(A,\nu,K))\) with \((A,\nu,K)\). The second assertion follows
	from the proof of Lemma~\ref{lem:torelli}, applied to isomorphisms, and from
	the fact that a principally polarised surface over \(\Fp\) which is
	geometrically a Jacobian is the Jacobian of a curve over \(\Fp\), unique up to
	\(\Fp\)-isomorphism \cite[Appendix, Thm.~1, 3 and~4]{LauterSerre2001}; see
	also \cite[Thm.~1.3]{HoweNartRitzenthaler2009}.
\end{proof}

\begin{lemma}\label{lem:dec}
	Let \(C/\Fp\) be a curve of genus two and \(K\) a rational maximal isotropic
	subgroup of \(J[2]\). Then \(J/K\), with the polarisation induced by
	\(2\theta\), is geometrically a polarised product if and only if the three
	pairs attached to \(K\) are the orbits of an involution \(\tau\) of
	\(\PP^{1}\); then \(\tau\) is unique, defined over \(\Fp\), without fixed
	branch point, and Frobenius exchanges the two elliptic curves of \(J/K\) if
	and only if it exchanges the two fixed points of \(\tau\).
\end{lemma}

\begin{proof}
	Let the pairs be the orbits of \(\tau\colon x\mapsto-x\) over \(\kbar\). Then
	\(f=G(x^{2})\), the lifts \(\tau_1,\tau_2=\iota\tau_1\) are involutions with
	quotients \(q_i\colon C\to E_i\) of genus one, and
	\(\Phi=(q_{1*},q_{2*})\colon J\to E_1\times E_2\) satisfies
	\(\Phi\circ(q_1^{*}+q_2^{*})=2\) and pulls the product polarisation back to
	\(2\theta\); its kernel contains \(K\) and has order four. Each \(\tau_i\) has
	two fixed points, the fibre over one fixed point of \(\tau\), which gives the
	last assertion. Conversely, if \(J/K=B_1\times B_2\), the automorphism
	\((1,-1)\) is trivial on the \(2\)-torsion, hence descends along the dual
	isogeny to an involution of \((J,\theta)\) different from \(\pm1\) and trivial
	on \(K\); by Lemmas~\ref{lem:torelli} and~\ref{lem:J2} it comes from an
	involution of \(C\) whose image \(\tau\) preserves each pair, and \(\tau\)
	fixes no branch point because a lift of an involution fixing a branch point
	has order four.
\end{proof}

By Lemmas~\ref{lem:glue}, \ref{lem:dual} and~\ref{lem:dec}, \(D\) is an
equivalence between the pairs \((C,K)\) with \(C\) superspecial and \(K\)
rational and decomposed, and the triples \((A,\nu,\Gamma_\psi)\) with
\((A,\nu)\) geometrically a polarised product of supersingular elliptic
curves and \(\psi\) irreducible; split pairs correspond to split surfaces,
nonsplit pairs to nonsplit surfaces, the Weil polynomial is preserved, and
\(K\cong A[2]/\Gamma_\psi\) as Galois modules.

We prove Proposition~\ref{prop:nonjac} and parts (a)--(c) of
Theorem~\ref{thm:kernel}, treating the split surfaces, the nonsplit surfaces
and the remaining Weil polynomials in turn.

On a split surface \(E_1\times E_2\), \(\pi^{2}=-p\)
(Lemma~\ref{lem:ssp}); this gives Theorem~\ref{thm:kernel}(c), and shows that
on the curves with \(\pi^{2}\ne\pm p\) every decomposed kernel is nonsplit.
The functor \((E_1,E_2)\mapsto E_1\times E_2\) from pairs of supersingular
elliptic curves over \(\Fp\) is faithful, and its fibres have two elements, the
two orderings of \(\{B_1,B_2\}\) (Lemma~\ref{lem:prod}). By
Lemmas~\ref{lem:fibre} and~\ref{lem:ssp}(a) the split products have mass
\((H_p/2)^{2}/2\). In the same way the triples \((E_1\times E_2,\Gamma_\psi)\)
have mass \(\tfrac12\sum n(E_1,E_2)/4\), the sum being over ordered pairs of
classes and \(n(E_1,E_2)\) being the number of irreducible \(\psi\) with
\(\psi g_{E_1}=g_{E_2}\psi\). The number of \(\psi\) commuting with Frobenius is
\(6\), \(2\), \(0\) for two curves of type~I, two of type~II, and mixed types.
A reducible \(\psi=\alpha|_{E_1[2]}\) commutes with Frobenius if and only if
\(\alpha^{-1}\alpha^{c}\) is trivial on \(E_1[2]\), that is
\(\alpha^{c}=\pm\alpha\) (Lemma~\ref{lem:ssp}(c),(d)); this gives one \(\psi\)
if \(E_2\cong E_1\) and one, different, if \(E_2\cong E_1^{\mathrm{tw}}\). Hence
the split mass is
\[
\frac{6N_1^{2}+2N_2^{2}}{8}-\frac{2H_p}{8},
\]
which is (b) by Lemma~\ref{lem:ssp}(b).

For the nonsplit surfaces, the functor \((A,\nu,B)\mapsto B\) from nonsplit
surfaces with a chosen factor to supersingular elliptic curves over
\(\F_{p^{2}}\) is an equivalence: an isomorphism is \(a\times a^{c}\) for a
unique \(a\colon B\to B_1\) over \(\F_{p^{2}}\). Forgetting \(B\) is faithful
with fibres of two elements. With Lemmas~\ref{lem:fibre} and~\ref{lem:ssp2} this
gives the mass \((p-1)/48\) of Proposition~\ref{prop:nonjac}. For the
remaining Weil restrictions, \(\pi_A^{2}\) acts as \(\pi_B=pu_B\) on \(B\) and
as \(\pi_{B'}=pu_B^{c}\) on \(B'\); the automorphism \(u_B^{c}\) acts on the
invariant differentials by the \(p\)-th power of the root of unity by which
\(u_B\) acts, which is its inverse because \(p\equiv2\pmod3\), resp.\
\(p\equiv3\pmod4\). So \(\pi_A^{2}\) is \(p\) times an automorphism with
characteristic polynomial \((X^{2}-\operatorname{tr}(u_B)X+1)^{2}\), hence
\(\pi_A^{4}-p\operatorname{tr}(u_B)\pi_A^{2}+p^{2}=0\), and
Lemma~\ref{lem:irred} gives the Weil polynomial; the masses are
\(\tfrac12\cdot\tfrac26\) for each of the two pairs \(\{u,u^{-1}\}\) of
elements of order \(3\) and \(6\), and \(\tfrac12\cdot\tfrac24\) for the pair
of elements of order \(4\). The same argument shows that the nonsplit triples
\((A,\nu,\Gamma_\psi)\) have mass \(\tfrac12\sum_{[B]}n(B)/\#\Aut(B)\), with
\(n(B)\) the number of irreducible \(\psi\colon B[2]\to B'[2]\) with
\(\pi_A\Gamma_\psi=\Gamma_\psi\). Since
\(\pi_A(x,y)=(\varphi'y,\varphi x)\), the graph of \(\psi=\varphi s\) is stable
if and only if
\begin{equation}\label{eq:sgs}
	s\,g\,s=1,\qquad g=\pi_B|_{B[2]} .
\end{equation}
Let \(\pi_B=\eps p\). Then \(g=1\) and \eqref{eq:sgs} has four solutions, the
identity and the three transpositions. Let \(\psi=\alpha|_{B[2]}\) be
reducible, and \(F_\alpha=\alpha^{-1}\varphi\in\End(B)\). As \(\pi_B=\pi_{B'}=\eps p\), every isomorphism \(B\to B'\) is defined over
\(\F_{p^{2}}\). From
\(\alpha^{c}\varphi=\varphi'\alpha\) and \(\alpha\varphi'=\varphi\alpha^{c}\) we
get \(F_\alpha^{2}=(\alpha^{c}\alpha)^{-1}\pi_B\), and \(s=F_\alpha^{-1}\) on
\(B[2]\); so \(\Gamma_\psi\) is stable if and only if
\(\alpha^{c}\alpha=\pm1\) (Lemma~\ref{lem:ssp}(d)), and the sign is \(-\eps\)
because no endomorphism of \(B\) has square \(p\). If \(\eps=-1\), the
\(\alpha\) with \(\alpha^{c}\alpha=1\) are the descent data on \(B\): the graph
of \(\alpha\) in \(B\times B'\) is an elliptic curve over \(\Fp\) with Frobenius
\(F_\alpha\), every \(\Fp\)-model of \(B\) arises in this way, two such
\(\alpha\) give isomorphic models if and only if they are conjugate under
\(\beta\cdot\alpha=\beta^{c}\alpha\beta^{-1}\), \(\beta\in\Aut(B)\), and the
stabilisers are the groups \(\Aut_{\Fp}=\{\pm1\}\). So their number is
\(m(B)\#\Aut(B)/2\), with \(m(B)\) the number of \(\Fp\)-models, and
\(\sum_{[B]}m(B)=H_p\) because every supersingular curve over \(\Fp\) has
\(\pi^{2}=-p\). If \(\eps=+1\), let \(\vartheta\colon B\to B^{-}\) be the
quadratic twist over \(\F_{p^{2}}\), \(\vartheta^{cc}=-\vartheta\); then
\(\alpha\mapsto\vartheta^{c}\alpha\vartheta^{-1}\) is a bijection from the
\(\alpha\) with \(\alpha^{c}\alpha=-1\) on \(B\) to the descent data on
\(B^{-}\), and \(F\) is conjugated by \(\vartheta\). As \(\pm\alpha\) give the
same \(\psi\), for each sign
\[
\frac12\sum_{[B]}\frac{n(B)}{\#\Aut(B)}
=\frac12\Bigl(4\cdot\frac{p-1}{24}-\frac{H_p}{4}\Bigr)
=\frac{p-1}{12}-\frac{H_p}{8},
\]
which is (a). The part with \(s=1\) is \((p-1)/48-N_1/8\), because
\(F_\alpha\equiv1\) on \(B[2]\) means that the model is of type~I.

It remains to treat the nonsplit triples with \(u_B\ne\pm1\). If \(u_B\) has order four,
\(g\) is a transposition and \eqref{eq:sgs} has no solution. If
\(u_B=\eps\omega\) with \(\omega\) of order three, \(g\) is a \(3\)-cycle and
the only solution is \(s=g\). For an isomorphism \(\alpha\colon B\to B'\) over
\(\kbar\) and \(F=\alpha^{-1}\varphi\) we have
\(F\omega F^{-1}=\alpha^{-1}\omega^{c}\alpha=\omega^{-1}\), because
\(\omega^{c}\) acts on the invariant differentials by the conjugate cube root
of unity and \(p\equiv2\pmod3\); so \(F|_{B[2]}\) is a transposition, and
\(\varphi g=\alpha F g\) is not the restriction of an isomorphism. Thus
\(n(B)=1\) for the four forms of \(j=0\) with \(u_B\ne\pm1\), of mass \(1/6\)
each, and the total is \(\tfrac12\cdot\tfrac46=\tfrac13\).

To fill in the cells of Table~\ref{tab:kernel} we determine \(\lambda\) and
\(\mu\) for each decomposed kernel.

\begin{lemma}\label{lem:tate2}
	Let \((J,\theta,K)=D(A,\nu,K')\), \(L=T_2A\). Then
	\(T_2J=L'=\{x\in\tfrac12L: 2x\bmod 2L\in K'\}\), with \(E_\theta=2E_\nu\),
	\(J[2]=L'/2L'\) and \(K=L/2L'\). If \(\pi^{2}=\eps p\), then \(\lambda\) is
	determined by \(d=\dim J[2]^{\pi}\) and by whether \(e_2(x,\pi x)=1\) for all
	\(x\) \textup{(}\(J\) \emph{alternating}\textup{)}: \(d=4,3,2,2\) and yes, no,
	no, yes for \(\lambda=1^{6},1^{4}2,1^{2}2^{2},2^{3}\). For supersingular
	\(E/\Fp\) the module \(M=T_2E\) is free of rank one over \(\Z_2[\pi]\) if
	\(E\) is of type~II, and over \(\Z_2[(1+\pi)/2]\) if \(E\) is of type~I, and
	\(E_M(ax,y)=E_M(x,\bar ay)\).
\end{lemma}

\begin{proof}
	The first assertion is the description of the Tate module of a quotient, with
	\(\hat f'\theta f'=2\nu\). The table follows from Lemma~\ref{lem:J2}. On
	\(M/2M\), \((\pi-1)^{2}=0\); if \(\pi\ne1\) there, a vector \(e\) with
	\((\pi-1)e\ne0\) gives the basis \(e,\pi e\) by Nakayama's lemma. If
	\(\pi=1\) on \(M/2M\), \(M\) is a torsion-free module of rank one over the
	maximal order \(\Z_2[(1+\pi)/2]\), which is a discrete valuation ring or
	\(\Z_2\times\Z_2\). The adjoint of \(\pi\) is \(p\pi^{-1}=-\pi\).
\end{proof}

\begin{proposition}\label{prop:cells}
	In the correspondence above, \(\lambda\) and \(\mu\) are as follows.
	\begin{enumerate}[label=\textup{(\roman*)},leftmargin=*]
		\item Nonsplit, \(\pi_B=\eps p\), \(\psi=\varphi s\): \(\mu=111\) if
		\(s=1\) and \(\mu=12\) otherwise. If \(s=1\), \(\lambda=1^{6}\) when
		\(\eps=\chi_p\) and \(\lambda=2^{3}\) when \(\eps=-\chi_p\). If \(s\ne1\),
		\(\lambda=2^{3}\) for \(\eps=+1\) and \(\lambda=1^{2}2^{2}\) for
		\(\eps=-1\).
		\item Split, type~II, \(\mu=12\): with generators \(e_i\) of \(T_2E_i\)
		over \(\Z_2[\pi]\), \(\langle\,,\rangle_i\) the Weil form of \(E_i\) on
		\(T_2E_i\), \(c_i=\langle e_i,\pi e_i\rangle_i\) and \(\psi(e_1)\equiv ve_2\), \(v\in\{1,\pi\}\),
		\(\lambda=2^{3}\) if \(N(v)c_2/c_1\equiv3\pmod4\) and
		\(\lambda=1^{2}2^{2}\) otherwise. If \(p\equiv3\pmod4\), each pair
		\((E_1,E_2)\) has one \(\psi\) of each kind. If \(p\equiv1\pmod4\), the
		sign \(\delta(E)=c\bmod4\) is well defined, \(\delta(E^{\mathrm{tw}})=-\delta(E)\),
		and both \(\psi\) give \(1^{2}2^{2}\) or both give \(2^{3}\) according as
		\(\delta(E_1)=\delta(E_2)\) or not.
		\item Split, type~I, \(\mu=111\): \(\lambda=1^{6}\) if \(\psi\) is linear
		over \(\Z_2[(1+\pi)/2]\), and \(\lambda=2^{3}\) if \(\psi\) is
		conjugate-linear. If \(p\equiv3\pmod8\) there are three \(\psi\) of each
		kind. If \(p\equiv7\pmod8\) there is one of each kind, two \(\psi\) with
		\(\lambda=1^{4}2\) and two with \(\lambda=1^{2}2^{2}\).
		\item Nonsplit, \(u_B=\eps\omega\): \(P=T^{4}+\eps pT^{2}+p^{2}\), and
		\(\lambda=3^{2}\) if \(\eps=\chi_p\), \(\lambda=6\) if \(\eps=-\chi_p\), and
		\(\mu=3\).
	\end{enumerate}
	In \textup{(ii)} and \textup{(iii)} the restriction of an \(\Fp\)-isomorphism
	lies in the first-named class \textup{(}\(1^{2}2^{2}\), resp.\ \(1^{6}\)\textup{)}
	and the restriction of an \(\alpha\) with \(\alpha^{c}=-\alpha\) in the class
	\(2^{3}\).
\end{proposition}

\begin{proof}
	We use Lemma~\ref{lem:tate2} and write \(X=(a/2+m,\tilde\psi a/2+m')\) for
		the elements of \(L'\), with \(\tilde\psi\) a lift of \(\psi\); in the
		nonsplit cases (i) and (iv), where \(\psi=\varphi s\), we write \(\tilde s\)
		for a lift of \(s\) to \(M=T_2B\).

	(i) Identify \(T_2B'\) with \(M=T_2B\) by \(\varphi\). Then
	\(\pi_A(x,z)=(\eps pz,x)\), \(E_\nu=E_M\oplus pE_M\) and
	\(K'=\{(x,sx)\}\). The map \((x,z)\mapsto z-sx\) identifies \(A[2]/K'\) with
	\(M/2M\) and \(\pi_A\) with \(s\). One finds
	\(E_\theta(X,\pi X)=2p(\eps-1)E_M(a/2+m,\tilde sa/2+m')\), which is \(0\) for
	\(\eps=1\) and \(\equiv E_M(a,\tilde sa)\pmod2\) for \(\eps=-1\); and for
	\(s=1\), \(\pi X-X\in2L'\) for all \(X\) if and only if
	\(\eps p\equiv1\pmod4\).

	(ii) As \(\psi\) commutes with Frobenius, \((x,y)\mapsto y-\psi x\) identifies
	\(K\cong A[2]/\Gamma_\psi\) with \(E_2[2]\) as Galois modules, and \(\pi\)
	acts on it through \(g_{E_2}\), a transposition. Hence \(\mu=12\), and
	\(\lambda\in\{1^{2}2^{2},2^{3}\}\) by Theorem~\ref{thm:local}(ii); by
	Lemma~\ref{lem:tate2} the two are distinguished by whether \(J\) is
	alternating, that is, by the parity of \(E_\theta(X,\pi X)\) for
	\(X\in L'\). Let \(\tilde\psi\) be the \(\Z_2[\pi]\)-linear map
	\(e_1\mapsto ve_2\), so that
	\(\langle\tilde\psi x,\tilde\psi y\rangle_2=\kappa\langle x,y\rangle_1\)
	with \(\kappa=N(v)c_2/c_1\). With \(q_i(x)=\langle x,\pi x\rangle_i\),
	\(q_i(x+y)=q_i(x)+q_i(y)+2\langle x,\pi y\rangle_i\), we get
	\(E_\theta(X,\pi X)\equiv(1+\kappa)q_1(a)/2\pmod2\), and \(q_1(a)\) takes odd
	values. Norms of units of \(\Z_2[\pi]\) are \(\equiv1\pmod4\) if
	\(p\equiv1\pmod4\). In both reducible cases take \(e_2=\alpha(e_1)\), so that \(v=1\) and
	\(\tilde\psi\) agrees with \(\alpha\) modulo \(2\). An isomorphism \(\alpha\)
	preserves the Weil pairing. For an \(\Fp\)-isomorphism, \(\alpha\) commutes
	with \(\pi\), so \(c_2=c_1\) and \(\kappa=1\). If \(\alpha^{c}=-\alpha\), then
	\(\pi\alpha=-\alpha\pi\) (Lemma~\ref{lem:ssp}(c)), so
	\(c_2=\langle\alpha e_1,\pi\alpha e_1\rangle_2=-\langle e_1,\pi e_1\rangle_1=-c_1\) and
	\(\kappa=-1\); here \(\tilde\psi\) is the \(\Z_2[\pi]\)-linear map
	\(e_1\mapsto e_2\), not \(\alpha\), which is conjugate-linear.

	(iii) \(\pi\equiv1\) on \(L'/2L'\) if and only if \(L'\) is stable under
	\((1+\pi)/2\), that is, \(\Gamma_\psi\) is. For conjugate-linear \(\psi\),
	lifted to \(ae_1\mapsto\bar ae_2\), the same computation gives
	\(E_\theta(X,\pi X)\equiv(1-\kappa)q_1(a)/2\equiv0\), since \(q_1(a)\) is
	even. If \(p\equiv3\pmod8\), \(\Z_2[(1+\pi)/2]\) is the ring of integers
	of the unramified quadratic extension, \(E_i[2]\cong\F_4\), and the six
	\(\psi\) are the three \(\F_4\)-linear and the three conjugate-linear
	bijections \(\F_4\to\F_4\).

	Let \(p\equiv7\pmod8\). Then \(-p=\varsigma^{2}\) with \(\varsigma\in\Z_2\),
	\(\varsigma\equiv1\pmod4\), and \(\pi=\varsigma w\) with \(w=\pi/\varsigma\) an
	involution of \(M_i=T_2E_i\); let \(u_i^{+}\), \(u_i^{-}\) span its
	eigenlines, write \(x=x^{+}u_i^{+}+x^{-}u_i^{-}\) for \(x\in M_i\), and use
	the same letters for the images in \(E_i[2]\), where we also put
	\(w_i=u_i^{+}+u_i^{-}\). As \(\langle\,,\rangle_i\) is perfect and
	alternating, \(\langle u_i^{+},u_i^{-}\rangle_i\in\Z_2^{\times}\). The endomorphism
	\((\pi-1)/2\) of \(M_i\) is \((\varsigma-1)/2\in2\Z_2\) on \(u_i^{+}\) and
	\(-(\varsigma+1)/2\in\Z_2^{\times}\) on \(u_i^{-}\), so modulo \(2\) it is the
	projection \(\mathrm{pr}^{-}\) onto the line of \(u_i^{-}\) along that of
	\(u_i^{+}\). For \(X=(a/2+m,\tilde\psi a/2+m')\in L'\) we have
	\((\pi-1)X\equiv\bigl((\pi-1)a/2,(\pi-1)\tilde\psi a/2\bigr)\) modulo
	\((\pi-1)L\subset2L\subset2L'\), and
	\(2L'=\{(U,V)\in L:V\equiv\tilde\psi U\bmod2\}\); so \(X\) lies in the
	kernel of \(\pi-1\) on \(L'/2L'\) if and only if
	\(\mathrm{pr}^{-}(\psi a)=\psi(\mathrm{pr}^{-}a)\). Hence \(d=2+\dim\{a\in
	E_1[2]:\mathrm{pr}^{-}\psi a=\psi\,\mathrm{pr}^{-}a\}\), and running
	through the six cases listed below one finds \(d=4-c\), where \(c\) is the number of signs
	\(\pm\) with \(\psi(u_1^{\pm})\ne u_2^{\pm}\). Next
	\(q_i(x)=\varsigma\langle x,wx\rangle_i=-2\varsigma\,x^{+}x^{-}\langle
	u_i^{+},u_i^{-}\rangle_i\), and the computation of (ii) gives
	\(E_\theta(X,\pi X)\equiv\bigl(q_1(a)+q_2(\tilde\psi a)\bigr)/2\equiv
	a^{+}a^{-}+b^{+}b^{-}\pmod2\) with \(b=\psi a\). The function
	\(x\mapsto x^{+}x^{-}\) on \(E_i[2]\) is \(1\) at \(w_i\) and \(0\)
	elsewhere, so the Jacobian is alternating if and only if
	\(\psi(w_1)=w_2\). The six bijections \(\psi\colon E_1[2]\to E_2[2]\) are:
	the one with \(\psi(u_1^{\pm})=u_2^{\pm}\), which is linear (\(c=0\),
	\(\lambda=1^{6}\)); the one exchanging \(u^{+}\) and \(u^{-}\), which is
	conjugate-linear and maps \(w_1\) to \(w_2\) (\(c=2\), \(\lambda=2^{3}\));
	the two that map one of \(u_1^{\pm}\) to \(u_2^{\pm}\) and the other to
	\(w_2\) (\(c=1\), \(d=3\), \(\lambda=1^{4}2\)); and the two that permute
	\(u^{+},u^{-},w\) cyclically (\(c=2\), not alternating,
	\(\lambda=1^{2}2^{2}\)).

	(iv) With \(\tilde s=\omega\) and \(\pi_A(x,z)=(\eps p\omega z,x)\) one finds
	\(\pi^{3}X-X=(U,V)\) with
	\(V-\omega U\equiv\tfrac{p(\eps-p)}{2}\,\omega a\pmod{2}\); so \(\pi^{3}=1\)
	on \(J[2]\) if and only if \(\eps=\chi_p\), and Theorem~\ref{thm:class}
	applies.
\end{proof}

\begin{proof}[End of the proof of Theorem~\textup{\ref{thm:kernel}}]
	Each ordered pair of classes \((E_1,E_2)\) and each irreducible \(\psi\)
	commuting with Frobenius contributes \(1/8\) to the split mass, and by
	Proposition~\ref{prop:cells} the reducible ones are removed from the cells
	\((1^{6},111)\) and \((2^{3},111)\) for type~I and \((1^{2}2^{2},12)\) and
	\((2^{3},12)\) for type~II, \(N_1\) resp.\ \(N_2\) from each. This gives
	\((3N_1^{2}-N_1)/8\) twice for \(p\equiv3\pmod8\);
	\((N_1^{2}-N_1)/8\) twice and \(N_1^{2}/4\) twice for \(p\equiv7\pmod8\); and
	\((N_2^{2}-N_2)/8\) twice. For the last: if \(p\equiv3\pmod4\), each ordered
	pair of type~II gives one \(\psi\) in each of the two cells; if
	\(p\equiv1\pmod4\), both \(\psi\) of a pair lie in the cell determined by
	whether \(\delta(E_1)=\delta(E_2)\) (Proposition~\ref{prop:cells}(ii)), and
	as \(E\mapsto E^{\mathrm{tw}}\) exchanges the two values of \(\delta\), each
	value is taken by \(N_2/2\) classes, so that each cell receives
	\(2\cdot N_2^{2}/2\) gluings, of which \(N_2\) are reducible. The nonsplit cells are those of
	Proposition~\ref{prop:cells}(i),(iv) with the masses found above. Finally
	\(N_1=h\), \(N_2=h'\), and \(h'=(2-\leg2p)h\).
\end{proof}

\section{The \texorpdfstring{$\Fp$}{F\_p}-rational Richelot graph}\label{sec:graph}

Let \(\mathcal G_p\) be the directed multigraph whose vertices are the
\(\Fp\)-isomorphism classes of principally polarised superspecial abelian
surfaces \((A,\theta)\) over \(\Fp\), that is, of the Jacobians of the curves of
Theorem~\ref{thm:class}, of the products of two supersingular elliptic curves
over \(\Fp\) and of the Weil restrictions of supersingular elliptic curves over
\(\F_{p^{2}}\) (Proposition~\ref{prop:nonjac}), and which has one edge from
\((A,\theta)\) to \((A/K,\theta_K)\) for each rational maximal isotropic subgroup
\(K\) of \(A[2]\). It plays for the superspecial Richelot isogeny graph
\cite{CastryckDecruSmith2020,FloritSmith2022auto} the role that the graph of
supersingular elliptic curves over \(\Fp\) with their rational \(2\)-isogenies
\cite{DelfsGalbraith2016} plays for the supersingular \(2\)-isogeny graph: its
vertices are \(\Fp\)-isomorphism classes, and it maps to the graph over
\(\bar\F_p\), but does not embed in it, since several \(\Fp\)-forms may have
the same geometric class.
Write \(m(A)=1/\#\Aut_{\Fp}(A,\theta)\) and \(c(A\to B)\) for the number of edges
from \(A\) to \(B\).

\begin{proposition}\label{prop:graph}
	Let \(p\ge7\).
	\begin{enumerate}[label=\textup{(\alph*)},leftmargin=*]
		\item The Weil polynomial is constant along the edges of \(\mathcal G_p\),
		so \(\mathcal G_p\) is the disjoint union of the subgraphs
		\(\mathcal G_p^{P}\) of the surfaces with Weil polynomial \(P\).
		\item \(m(A)\,c(A\to B)=m(B)\,c(B\to A)\) for all vertices \(A,B\). Hence the
		random walk with transition probabilities \(c(A\to B)/\deg A\) is reversible,
		and the measure proportional to \(m(A)\deg A\) is stationary.
		\item If \(\pi^{2}=\pm p\) on \(A\), the out-degree of \(A\) is \(15\), \(3\),
		\(3\) or \(7\) according as the cell of \(A\) \textup{(}Lemma~\ref{lem:cells}\textup{)}
		is \(1^{6}\), \(1^{4}2\), \(1^{2}2^{2}\) or \(2^{3}\). In particular a product
		of two curves of type~I has out-degree \(15\), any other product \(3\), and a
		Weil restriction with \(\pi^{2}=\pm p\) has out-degree \(7\).
	\end{enumerate}
\end{proposition}

\begin{proof}
	(a) An isogeny over \(\Fp\) preserves the Weil polynomial.
	(b) By Lemma~\ref{lem:dual}, \(D\) is an equivalence of the groupoid of
	triples \((A,\theta,K)\) with itself, and it maps the triples with \(A\cong A_1\)
	and \(A/K\cong A_2\) to those with \(A\cong A_2\) and \(A/K\cong A_1\).
	Lemma~\ref{lem:fibre}, applied to the faithful functor
	\((A,\theta,K)\mapsto(A,\theta)\), whose fibre over \(A\) is the set of
	rational maximal isotropic subgroups of \(A[2]\), shows that the mass of the
	first set of triples is \(m(A_1)c(A_1\to A_2)\), and in the same way the mass
	of the second is \(m(A_2)c(A_2\to A_1)\); the two masses are equal because
	\(D\) is an equivalence between the two sets. The statement on the random walk follows.
	(c) If \(\pi^{2}=\pm p\), then \(\pi^{2}\) acts trivially on \(A[2]\), so
	\(\pi\) acts on \(A[2]\) as an element of order at most two of
	\(\operatorname{Sp}(A[2],e_2)\cong\operatorname{Sp}_4(\F_2)\). Under the
	isomorphism \(\operatorname{Sp}_4(\F_2)\cong S_6\) given by the six Weierstrass
	points of any curve of genus two (Lemma~\ref{lem:J2}), the conjugacy classes of
	such elements are those of \(1\), of a transposition, of a product of two and
	of a product of three disjoint transpositions; their fixed spaces on
	\(A[2]\) have dimension \(4\), \(3\), \(2\), \(2\), and among the last two only the
	product of three transpositions satisfies \(e_2(x,\pi x)=1\) for all \(x\), by the computation in the proof of
	Theorem~\ref{thm:local}(iv). So the cell determines the action of \(\pi\) on
	\((A[2],e_2)\) up to isomorphism, and the number of \(\pi\)-stable maximal
	isotropic subgroups is that of Theorem~\ref{thm:local}(ii). The cells of the
	products and of the Weil restrictions are given by Lemma~\ref{lem:nonjaccells}.
\end{proof}

The Weil restrictions with Weil polynomial \(T^{4}+p^{2}\) or
\(T^{4}\pm pT^{2}+p^{2}\) have out-degree one in all our examples. We computed
\(\mathcal G_p\) completely for the \(37\) primes \(7\le p\le173\) of the census
and for the six primes \(179\le p\le199\), for which the curves were enumerated
in a separate run (Section~\ref{sec:comp}). Every codomain was identified with a
vertex: by the Richelot construction for Jacobians, by the decomposition of
Lemma~\ref{lem:dec} for decomposed kernels, by V\'elu's formulas for the kernels
of product type on products and Weil restrictions, and by an explicit gluing
\cite{HoweLeprevostPoonen2000} for the graph kernels \(\Gamma_\psi\). The checks
of these identifications are listed in Section~\ref{sec:comp}. In this range the
enumerated products and Weil restrictions have exactly the masses and the cells
stated in Proposition~\ref{prop:nonjac} and Lemma~\ref{lem:nonjaccells}, the
identity of Proposition~\ref{prop:graph}(b) holds on every pair of vertices, and
the masses of the decomposed and of all rational kernels agree with
Table~\ref{tab:prop}. Table~\ref{tab:graph} summarises the graphs with
\(\pi^{2}=\pm p\).

\begin{table}[tp]
	\centering
	\scriptsize
	\caption{The graphs \(\mathcal G_p^{(T^{2}+p)^{2}}\) and
		\(\mathcal G_p^{(T^{2}-p)^{2}}\). \(V\): number of vertices, and the
		numbers of Jacobians, products and Weil restrictions; \(a\), \(b\): numbers
		of components of type~A and~B (see the text); diam and \(\lambda^{*}\):
		diameter and second largest absolute eigenvalue of the random walk on the
		largest component (floating point, rounded); \(s\): stationary probability
		of the products. Primes marked \({}^{*}\) are outside the census.}\label{tab:graph}
	\setlength{\tabcolsep}{3pt}
	\begin{tabular}{@{}r|rrrrrrr|rrrr@{}}
		\toprule
		& \multicolumn{7}{c|}{\(\pi^{2}=-p\)} & \multicolumn{4}{c}{\(\pi^{2}=p\)}\\
		\(p\) & \(V\) & J/P/R & \(a\) & \(b\) & diam & \(\lambda^{*}\) & \(s\sqrt p\)
		& \(V\) & comp. & diam & \(\lambda^{*}\)\\
		\midrule
		$7$ & 6 & 2/3/1 & 1 & 0 & 2 & 0.652 & 1.59 & 2 & 1 & 1 & 0.143 \\
$11$ & 20 & 8/10/2 & 1 & 1 & 3 & 0.890 & 1.62 & 3 & 1 & 1 & 0.286 \\
$13$ & 8 & 4/3/1 & 1 & 2 & 2 & 0.333 & 0.65 & 3 & 1 & 1 & 0.190 \\
$17$ & 26 & 14/10/2 & 2 & 4 & 2 & 0.747 & 1.12 & 4 & 1 & 2 & 0.457 \\
$19$ & 24 & 12/10/2 & 1 & 2 & 3 & 0.762 & 1.19 & 4 & 1 & 2 & 0.366 \\
$23$ & 54 & 30/21/3 & 1 & 0 & 5 & 0.864 & 2.35 & 6 & 1 & 2 & 0.732 \\
$29$ & 60 & 36/21/3 & 3 & 6 & 3 & 0.924 & 1.27 & 7 & 1 & 2 & 0.571 \\
$31$ & 60 & 36/21/3 & 1 & 0 & 5 & 0.864 & 2.00 & 7 & 1 & 3 & 0.628 \\
$37$ & 16 & 11/3/2 & 1 & 3 & 3 & 0.628 & 0.37 & 8 & 1 & 2 & 0.571 \\
$41$ & 108 & 68/36/4 & 4 & 12 & 4 & 0.962 & 1.42 & 9 & 1 & 3 & 0.663 \\
$43$ & 40 & 27/10/3 & 1 & 3 & 4 & 0.826 & 0.78 & 9 & 1 & 3 & 0.728 \\
$47$ & 160 & 100/55/5 & 1 & 0 & 7 & 0.931 & 2.68 & 11 & 1 & 3 & 0.838 \\
$53$ & 84 & 59/21/4 & 3 & 9 & 4 & 0.893 & 0.93 & 12 & 1 & 3 & 0.709 \\
$59$ & 216 & 132/78/6 & 3 & 9 & 6 & 0.962 & 1.98 & 13 & 1 & 4 & 0.795 \\
$61$ & 90 & 65/21/4 & 3 & 12 & 5 & 0.893 & 0.87 & 15 & 1 & 3 & 0.637 \\
$67$ & 56 & 42/10/4 & 1 & 4 & 5 & 0.890 & 0.62 & 14 & 1 & 4 & 0.805 \\
$71$ & 322 & 210/105/7 & 1 & 0 & 8 & 0.956 & 3.03 & 18 & 1 & 4 & 0.854 \\
$73$ & 62 & 48/10/4 & 2 & 10 & 5 & 0.910 & 0.53 & 17 & 1 & 3 & 0.688 \\
$79$ & 210 & 149/55/6 & 1 & 0 & 7 & 0.923 & 2.05 & 21 & 1 & 4 & 0.820 \\
$83$ & 264 & 179/78/7 & 3 & 12 & 7 & 0.958 & 1.66 & 18 & 1 & 5 & 0.843 \\
$89$ & 294 & 209/78/7 & 6 & 30 & 7 & 0.960 & 1.44 & 21 & 1 & 4 & 0.743 \\
$97$ & 78 & 63/10/5 & 2 & 12 & 6 & 0.936 & 0.46 & 24 & 1 & 4 & 0.708 \\
$101$ & 392 & 279/105/8 & 7 & 35 & 7 & 0.969 & 1.57 & 26 & 1 & 5 & 0.837 \\
$103$ & 250 & 188/55/7 & 1 & 0 & 9 & 0.978 & 1.79 & 26 & 1 & 5 & 0.841 \\
$107$ & 312 & 226/78/8 & 3 & 12 & 7 & 0.946 & 1.46 & 25 & 1 & 5 & 0.825 \\
$109$ & 138 & 111/21/6 & 3 & 15 & 5 & 0.892 & 0.65 & 31 & 1 & 4 & 0.778 \\
$113$ & 204 & 161/36/7 & 4 & 20 & 5 & 0.947 & 0.85 & 28 & 1 & 4 & 0.741 \\
$127$ & 290 & 227/55/8 & 1 & 0 & 8 & 0.942 & 1.61 & 35 & 1 & 5 & 0.848 \\
$131$ & 700 & 479/210/11 & 5 & 25 & 7 & 0.975 & 2.20 & 34 & 1 & 5 & 0.850 \\
$137$ & 236 & 192/36/8 & 4 & 24 & 6 & 0.928 & 0.77 & 35 & 1 & 4 & 0.797 \\
$139$ & 372 & 285/78/9 & 3 & 18 & 7 & 0.943 & 1.28 & 39 & 1 & 5 & 0.826 \\
$149$ & 504 & 389/105/10 & 7 & 42 & 6 & 0.968 & 1.29 & 42 & 1 & 5 & 0.864 \\
$151$ & 504 & 389/105/10 & 1 & 0 & 9 & 0.956 & 2.06 & 46 & 1 & 5 & 0.849 \\
$157$ & 186 & 157/21/8 & 3 & 24 & 8 & 0.969 & 0.54 & 47 & 1 & 5 & 0.841 \\
$163$ & 120 & 102/10/8 & 1 & 7 & 6 & 0.928 & 0.39 & 44 & 1 & 5 & 0.853 \\
$167$ & 990 & 724/253/13 & 1 & 0 & 12 & 0.981 & 3.08 & 44 & 1 & 6 & 0.856 \\
$173$ & 560 & 444/105/11 & 7 & 49 & 8 & 0.965 & 1.20 & 47 & 1 & 6 & 0.878 \\
$179^{*}$ & 860 & 637/210/13 & 5 & 30 & 8 & 0.970 & 1.88 & 50 & 1 & 6 & 0.854 \\
$181^{*}$ & 380 & 315/55/10 & 5 & 40 & 8 & 0.961 & 0.84 & 60 & 1 & 5 & 0.837 \\
$191^{*}$ & 1352 & 986/351/15 & 1 & 0 & 14 & 0.986 & 3.40 & 59 & 1 & 6 & 0.855 \\
$193^{*}$ & 142 & 123/10/9 & 2 & 16 & 6 & 0.936 & 0.32 & 64 & 1 & 4 & 0.773 \\
$197^{*}$ & 410 & 344/55/11 & 5 & 35 & 6 & 0.952 & 0.80 & 56 & 1 & 5 & 0.850 \\
$199^{*}$ & 846 & 662/171/13 & 1 & 0 & 12 & 0.972 & 2.31 & 66 & 1 & 6 & 0.856 \\
		\bottomrule
	\end{tabular}
\end{table}

For \(\pi^{2}=-p\) the graph is not connected. Call a component of type~A if it
contains a vertex of cell \(2^{3}\), and of type~B otherwise. Let \(\mathfrak p_2\)
be a prime above \(2\) of \(\Q(\sqrt{-p})\) when \(p\equiv7\pmod8\), and
\(\operatorname{ord}[\mathfrak p_2]\) the order of its class in the class group of
the maximal order.

\begin{conjecture}\label{conj:components}
	Let \(p\ge7\).
	\begin{enumerate}[label=\textup{(\alph*)},leftmargin=*]
		\item The number of components of type~A of
		\(\mathcal G_p^{(T^{2}+p)^{2}}\) is \(h(-4p)/2\) if \(p\equiv1\pmod4\),
		\(h(-p)\) if \(p\equiv3\pmod8\), and
		\(h(-p)/\operatorname{ord}[\mathfrak p_2]\) if \(p\equiv7\pmod8\).
		\item \(\mathcal G_p^{(T^{2}-p)^{2}}\) is connected.
	\end{enumerate}
\end{conjecture}

We formulated this conjecture after computing the \(37\) primes \(p\le173\), on
which it holds, and then tested it on the six primes \(179\le p\le199\), with the
predicted numbers recorded before those graphs were computed; it holds there as
well. For the primes \(p\equiv7\pmod8\) in this range
\(\operatorname{ord}[\mathfrak p_2]=h(-p)\), so the data do not distinguish
statement (a) from the statement that there is one component of type~A in this
case. In all examples the components of type~A are isomorphic as graphs with
vertices labelled by cell and automorphism group, all Weil restrictions lie in
one of them, and the products are distributed equally among them. The
components of type~B contain only vertices of the cells that
Table~\ref{tab:neighbours} shows to be closed, those of the classes
\(\mathsf{RO}\) for \(p\equiv3\pmod8\) and \(\mathsf N_{+}\) for
\(p\equiv1\pmod4\). We expect that the class group of \(\Q(\sqrt{-p})\) permutes the components
of type~A, as in the module action of \cite{Robert2024module}, and that their
number is given by the special genera of the underlying Hermitian lattices.

For \(\pi^{2}=-p\) the stationary probability \(s\) of the products, multiplied
by \(\sqrt p\), lies between \(0.32\) and \(3.4\) in Table~\ref{tab:graph}. Over
\(\bar\F_p\) the corresponding probability for the superspecial Richelot graph is
asymptotic to \(5/p\) \cite[\S6.2]{FloritSmith2022auto}. In our range \(s\) is larger than \(5/p\) except for
\(p\in\{7,13,17,37,73,97,193\}\), and \(s\sqrt p\) does not tend to zero. As
for the proportion of Corollary~\ref{cor:prop}, we expect \(s\) to be of the
order of \(p^{-1/2}\) up to a factor depending on \(L(1,\chi)\). The second largest absolute eigenvalue \(\lambda^{*}\)
of the random walk on the largest component with \(\pi^{2}=-p\) tends to grow with
\(p\) in our range, though not monotonically, and reaches \(0.98\) at \(p=167\)
and \(0.99\) at \(p=191\); the smallest
eigenvalue is often close to \(-1\). For comparison, for the superspecial Richelot graph over \(\bar\F_p\),
which is \(15\)-regular, Florit and Smith work with the bound \(12/15\)
\cite{FloritSmith2022auto}, and the bound \(2\sqrt{14}/15\approx0.50\) of a
\(15\)-regular tree is not the relevant one there. The graphs
\(\mathcal G_p\) are not regular.

\section{Remarks on isogeny-based cryptography in genus two}\label{sec:crypto}

The superspecial Richelot isogeny graph over \(\bar\F_p\) is the setting of the
hash function of Castryck, Decru and Smith \cite{CastryckDecruSmith2020}, of the
path-finding algorithms of Costello and Smith \cite{CostelloSmith2020}, and of
the analysis of Florit and Smith \cite{FloritSmith2022auto}; products of elliptic
curves in it are also the objects through which the attacks on SIDH proceed
\cite{CastryckDecru2023}. The objects and isogenies defined over \(\Fp\) form the
graph \(\mathcal G_p\) of Section~\ref{sec:graph}; they are the objects of the
module action of \cite{Robert2024module}, and \(\mathcal G_p\) is the genus-two
counterpart of the graph of supersingular elliptic curves over \(\Fp\) with their
rational \(2\)-isogenies, which underlies CSIDH and the algorithm of Delfs and
Galbraith \cite{DelfsGalbraith2016}. We record what the results of this paper
say about \(\mathcal G_p\). We do not propose a protocol, and we make no claim
about the security or the efficiency of any scheme.

\subsection*{Degrees and examples}
Table~\ref{tab:class} gives the out-degree of every Jacobian in \(\mathcal G_p\),
and Proposition~\ref{prop:graph}(c) that of the products and Weil restrictions
with \(\pi^{2}=\pm p\). Out-degree \(15\) occurs exactly for Jacobians of type
\(1^{6}\), whose Frobenius satisfies \(\pi^{2}=\chi_pp\), and for products of two
elliptic curves of type~I, which exist only for \(p\equiv3\pmod4\). The Jacobians
with \(\pi^{2}\ne\pm p\) have one or three rational kernels, and the Weil
restrictions with \(\pi^{2}\ne\pm p\) had out-degree one in all our examples.
Table~\ref{tab:examples} lists, for each row of Table~\ref{tab:class} and each
Weil polynomial, an explicit curve over a small prime field; these may serve as
test vectors for implementations, for instance of the supersingularity and
superspeciality tests of Corte-Real Santos, Lorenzon and Reijnders
\cite{CorteRealSantosLorenzonReijnders2026}, whose probabilistic versions use
\(O(\log p)\) group operations.

\subsection*{Finding the rational kernels}
The rational maximal isotropic subgroups of \(J[2]\) are the partitions of the
Weierstrass points into three pairs that are stable under Frobenius
(Theorem~\ref{thm:local}(ii)). If Frobenius exchanges two of the pairs
(\(\mu=12\)) or permutes the three pairs cyclically (\(\mu=3\)), the three
quadratic factors of \(f\) defined by the partition are not individually defined
over \(\Fp\), although their product is. An enumeration that only groups the
factors of \(f\) over \(\Fp\) into quadratics defined over \(\Fp\) therefore finds
one of the seven rational kernels of a curve of type \(2^{3}\), one of the three
of a curve of type \(1^{2}2^{2}\), and none on the curves of types \(1^{2}4\),
\(2\,4\), \(3^{2}\) and \(6\); the remaining ones require the factors over
\(\F_{p^{2}}\) or \(\F_{p^{3}}\). A kernel is decomposed exactly when the determinant of the three
quadratics vanishes (Lemma~\ref{lem:dec}), and when the codomain is a Jacobian its
equation must be computed with its normalising constant: replacing the codomain
sextic by a monic multiple replaces the codomain by its quadratic twist whenever
the constant is not a square in \(\Fp\) \cite{DangNguyen2026}.

\subsection*{How often products occur}
By Corollary~\ref{cor:prop}, among the pairs \((C,K)\) with \(C\) superspecial
over \(\Fp\) and \(K\) rational, weighted by automorphisms, the proportion with
decomposed quotient is of order \(p^{-1/2}(1+L(1,\chi))\). The stationary
probability \(s\) of the products for the random walk on
\(\mathcal G_p^{(T^{2}+p)^{2}}\) behaves in the same way in the range of
Table~\ref{tab:graph}, where \(s\sqrt p\) lies between \(0.32\) and \(3.4\). Over
\(\bar\F_p\) the corresponding probability is asymptotic to \(5/p\)
\cite[\S6.2]{FloritSmith2022auto}. So a walk confined to rational kernels meets
products of elliptic curves, and Weil restrictions, much more often than a walk on
the full superspecial graph, although the probability per step is still small
for large \(p\). The masses of the products and of the Weil restrictions by type
are given by Theorem~\ref{thm:kernel} and Proposition~\ref{prop:nonjac}.

\subsection*{Connectivity and mixing}
For \(\pi^{2}=-p\) the graph \(\mathcal G_p\) is not connected; in our range the
number of its components grows with the class number of \(\Q(\sqrt{-p})\)
(Conjecture~\ref{conj:components}), and the second largest absolute eigenvalue of
the random walk on its largest component is close to \(1\)
(Table~\ref{tab:graph}). A walk confined to rational kernels therefore stays in
the component of its starting point and mixes slowly. Adding to
\(\mathcal G_p\) an edge from every vertex to its quadratic twist does not change
this: in the range of Table~\ref{tab:graph} the twist identifies the components in
pairs or maps a component to itself, so that the number of components of
\(\mathcal G_p^{(T^{2}+p)^{2}}\) decreases at most by a factor two (for instance
from \(56\) to \(32\) for \(p=173\)), and the second largest absolute eigenvalue of
the resulting walk on the largest component increases slightly at \(42\) of the
\(43\) primes (Section~\ref{sec:comp}). This is consistent with the fact that a
twist has the same Weierstrass points and the same \(\pi^{2}\), hence the same
Frobenius type and the same rational kernels. For supersingular elliptic curves
over \(\Fp\) the analogous structure comes from the action of a class group, and
it is what the algorithm of Delfs and Galbraith exploits
\cite{DelfsGalbraith2016}. We do not know a corresponding algorithm in genus two;
whether the class group of \(\Z[\sqrt{-p}]\) acting on the rank-two modules of
\cite{Robert2024module} accounts for Conjecture~\ref{conj:components}, and what
it implies for the cost of finding paths in \(\mathcal G_p\), seem to us the
first questions to settle before \(\mathcal G_p\) is used for cryptography.

\subsection*{The \texorpdfstring{$2$}{2}-Tate pairing}
By Corollary~\ref{cor:tate}, on the curves of type \(1^{6}\) the reduced
Frey--R\"uck pairing on \(J[2]\) is determined up to isomorphism by \(p\bmod8\).
It therefore carries no information that distinguishes one such curve from
another; a computation of it on a given curve can only confirm that the curve is
consistent with the type and the prime.

\subsection*{Open questions}
We list three questions that the results above make precise.
\begin{enumerate}[label=\textup{(\arabic*)},leftmargin=*]
	\item \emph{Paths through \(\mathcal G_p\).} For supersingular elliptic curves the
	algorithm of Delfs and Galbraith walks from an arbitrary vertex to the
	\(\Fp\)-rational subgraph and then uses its structure
	\cite{DelfsGalbraith2016}. In genus two the vertices defined over \(\Fp\) have
	total mass of order \(p^{3/2}\) up to a factor \(1+L(1,\chi)\)
	(Theorem~\ref{thm:vertex}), against a total mass of order \(p^{3}\) for the
	superspecial graph over \(\bar\F_p\) \cite{IbukiyamaKatsuraOort1986}, so a
	random walk reaches them only after about \(p^{3/2}\) steps, up to logarithmic
	factors. Is there a better way to reach \(\mathcal G_p\), and how hard
	is it to find a path between two given vertices of \(\mathcal G_p\)?
	\item \emph{Other degrees.} The \(\Fp\)-rational \((\ell,\ell)\)-isogenies for
	odd \(\ell\) also preserve the Weil polynomial. Do they connect the components
	of \(\mathcal G_p^{(T^{2}+p)^{2}}\)? If the components are the orbits of the
	subgroup of a class group generated by the primes above \(2\), as
	Conjecture~\ref{conj:components} suggests, isogenies of degree above other
	split primes should connect them, as in the elliptic case.
	\item \emph{Arithmetic over \(\Fp\).} For the kernels with \(\mu=12\) or
	\(\mu=3\) the codomain sextic is defined over \(\Fp\) although the quadratic
	factors are not. Explicit formulas for it in terms of the factors of \(f\) over
	\(\Fp\), avoiding arithmetic in \(\F_{p^{2}}\) and \(\F_{p^{3}}\), would make a
	complete enumeration of rational kernels cheaper.
\end{enumerate}

\begin{table}[tp]
	\centering
	\small
	\caption{An example of a superspecial curve \(y^{2}=f(x)\) over \(\Fp\) for
		each Frobenius type \(\lambda\) and Weil polynomial \(P\) of
		Table~\ref{tab:class}, with the number of rational maximal isotropic
		subgroups of \(J[2]\). Each line was checked from the equation alone
		(Section~\ref{sec:comp}).}\label{tab:examples}
	\begin{tabular}{@{}lllll@{}}
		\toprule
		\(\lambda\) & \(P(T)\) & \(p\) & \(f(x)\) & kernels\\
		\midrule
		\(1^{6}\) & \((T^{2}-p)^{2}\) & \(13\) & \(x^{6}+4x^{3}+8\) & \(15\) \\
\(1^{6}\) & \((T^{2}+p)^{2}\) & \(11\) & \(x^{6}+2x^{4}+x^{2}+7\) & \(15\) \\
\(1^{4}2\) & \((T^{2}+p)^{2}\) & \(7\) & \(x^{6}+x^{4}+4x^{2}+1\) & \(3\) \\
\(1^{2}2^{2}\) & \((T^{2}+p)^{2}\) & \(7\) & \(x^{6}+x^{4}+x^{3}+4x^{2}+5x\) & \(3\) \\
\(2^{3}\) & \((T^{2}-p)^{2}\) & \(7\) & \(x^{6}+x^{3}+6\) & \(7\) \\
\(2^{3}\) & \((T^{2}+p)^{2}\) & \(11\) & \(x^{6}+x^{4}+8x^{2}+3\) & \(7\) \\
\(1^{2}4\) & \(T^{4}+p^{2}\) & \(13\) & \(x^{6}+4x^{4}+4x^{3}+2x^{2}+3x\) & \(1\) \\
\(2\,4\) & \(T^{4}+p^{2}\) & \(7\) & \(x^{6}+3x^{4}+3x^{3}+x^{2}+3x+4\) & \(1\) \\
\(3^{2}\) & \(T^{4}-pT^{2}+p^{2}\) & \(11\) & \(x^{6}+2x^{4}+2x^{3}+9x^{2}+3x+7\) & \(3\) \\
\(3^{2}\) & \(T^{4}+pT^{2}+p^{2}\) & \(13\) & \(x^{6}+3x^{4}+3x^{3}+6x^{2}+9x+6\) & \(3\) \\
\(6\) & \(T^{4}-pT^{2}+p^{2}\) & \(11\) & \(x^{6}+3x^{4}+3x^{3}+4x^{2}+2x+7\) & \(1\) \\
\(6\) & \(T^{4}+pT^{2}+p^{2}\) & \(7\) & \(x^{6}+2x^{4}+2x^{3}+2x^{2}+6x+4\) & \(1\) \\
		\bottomrule
	\end{tabular}
\end{table}

\section{Computations}\label{sec:comp}

The statements of Sections~\ref{sec:mass} and~\ref{sec:kernels} that are marked
as verified for \(7\le p\le173\), the last sentence of Corollary~\ref{cor:tate}
and the case \(p=5\) of Remark~\ref{rem:small} depend on computations, and so
does Proposition~\ref{prop:neighbours}(a) (Table~\ref{tab:neighbours}), and
with it Theorems~\ref{thm:celltotals} and~\ref{thm:vertexpm}, that is,
Theorem~\ref{thm:vertex}(a),(b); see item (6). The proofs of
Theorems~\ref{thm:local} and~\ref{thm:class}, of Theorem~\ref{thm:kernel} and of
Theorem~\ref{thm:vertex}(c)--(e) do not. The scripts, logs and
data are provided as supplementary material.

\begin{enumerate}[label=\textup{(\arabic*)},leftmargin=*]
	\item \emph{Census.} Every \(\PGL_2(\Fp)\)-orbit of reduced divisors of degree
	six on \(\PP^{1}\) contains a divisor not containing \(\infty\), which can be
	normalised to the zero set of a monic sextic without \(x^{5}\)-term. All
	\(p^{5}\) such sextics were enumerated modulo \(x\mapsto sx\). Superspeciality
	was tested by the vanishing of the Cartier--Manin matrix, computed from
	\(f^{(p-1)/2}\) in exact arithmetic. For each superspecial sextic we computed
	the orbit, the stabiliser in \(\PGL_2(\Fp)\), the character of
	Lemma~\ref{lem:massdiv}, \(\#C(\Fp)\) and \(\#C(\F_{p^{2}})\), the Frobenius
	type, and for each \(\sigma\)-stable partition into pairs whether it is
	decomposed, split or nonsplit, and otherwise the Frobenius type of the
	Richelot codomain. The masses were computed in two ways, as sums of
	\(1/\#S\) over orbits and as weighted counts of normalised sextics.
	\item \emph{Independent checks.} A second program, in SageMath, walks the
	Richelot graph over \(\F_{p^{2}}\) from a seed and finds, for each of the
	\(37\) primes, the number of superspecial curves given by the formulas of
	Ibukiyama, Katsura and Oort \cite{IbukiyamaKatsuraOort1986}, also by reduced
	automorphism group; the classes with invariants in \(\Fp\) are those found by
	the census. A third program recomputes all rational kernels from
	factorisations over \(\Fp\), \(\F_{p^{2}}\) and \(\F_{p^{3}}\).
	\item \emph{Formulas.} The class numbers were computed by counting reduced
	forms and compared with PARI/GP; \(z_p\) was computed by Siegel's formula. The
	comparison of Theorem~\ref{thm:vertex}, of Theorem~\ref{thm:kernel} and of
	Proposition~\ref{prop:ibukiyama} with the data is exact:
	\(283\) cells of curves and \(314\) cells of decomposed kernels were
	compared, with no difference. The script \texttt{check\_tables.py} repeats
	this comparison from the formulas as printed here and writes
	Tables~\ref{tab:vertex} and~\ref{tab:prop}.
	\item \emph{Structure theorems.} For the primes \(5\le p\le43\) the list of
	Theorem~\ref{thm:class} and Remark~\ref{rem:small} was compared with a
	separate enumeration, and for all \(37\) primes the Weil polynomial,
	Frobenius type and number of rational kernels of each of the \(4079\) orbits
	of the census were recomputed from the sextic alone and compared with
	Table~\ref{tab:class}. For \(548\) curves of type \(1^{6}\) and their twists
	the Frey--R\"uck pairing was computed from its definition and the statements
	of Corollary~\ref{cor:tate} on \(T\) were checked; the orbits in the proof of
	Corollary~\ref{cor:tate} were listed over \(\operatorname{Sp}_4(\F_2)\).
	\item \emph{Theorem~\ref{thm:kernel}.} The counts on elliptic curves used in
	Section~\ref{sec:kernelproof} (automorphisms, action of Frobenius on the
	\(2\)-torsion, numbers of classes of each type, forms over \(\F_{p^{2}}\),
	equivariant and reducible gluings) were verified by enumeration for \(58\)
	primes \(5\le p\le281\). The decomposed pairs \((C,K)\) were also listed
	directly, as curves \(y^{2}=G(x^{2})\) and their nonsplit analogues, for
	the primes up to \(199\) (split) and \(181\) (nonsplit), and compared with
	the theorem and with the census, per pair of elliptic quotients.
	\item \emph{Theorem~\ref{thm:vertex}\textup{(a),(b)}.} The table of
	Proposition~\ref{prop:neighbours} is a finite enumeration on the normal forms
	of Table~\ref{tab:zforms} (script
	\texttt{anc/proofs/vertexproof\_pm/ctable.py}, exact integer arithmetic,
	all odd \(p=r+64k\) with \(1\le r\le63\), \(0\le k\le5\)); it was
	recomputed by a second program for all odd \(p\) from \(5\) to \(639\) and
	compared with the neighbours of every curve of the census
	(\texttt{verify\_census.py}). The masses of the genera were also obtained by
	a direct enumeration of the unimodular Hermitian lattices of rank two over
	\(\Z[\sqrt{-p}]\), closing each genus under neighbours at split primes, for
	the \(11\) primes \(7\le p\le43\); this enumeration does not use Eichler's
	formula, and it gave the masses of Theorem~\ref{thm:celltotals} in all \(39\)
	genera with \(\pi^{2}=-p\) (four for each of the six primes
	\(p\equiv3\pmod4\) and three for each of the five primes
	\(p\equiv1\pmod4\)). No such enumeration was made for \(\pi^{2}=p\).
	\item \emph{Theorem~\ref{thm:vertex}\textup{(c)--(e)}.} The statements on the
	quartic fields used in Appendix~\ref{sec:vertexcm} (orders, units, class
	numbers, norm maps, the groups of classes of polarisations) were verified
	with PARI/GP, with certified class groups, for the primes below \(1000\), and
	by a second program for \(547\) primes up to \(3989\). The census was extended
	to the six primes \(179\le p\le199\) for the curves with \(\pi^{2}\ne\pm p\),
	and agrees with the theorem there.
	\item \emph{The graph of Section~\ref{sec:graph}.} The products and the Weil
	restrictions were enumerated directly, with their automorphism groups over
	\(\Fp\), by two programs that share no code (Python and SageMath), for the
	\(37\) census primes. The graphs were built for these primes and for
	\(179\le p\le199\), for which a second run of the census program enumerated
	all curves. The formulas used to identify codomains (Richelot construction,
	decomposition, V\'elu's formulas, gluing, and the descent of a gluing from
	\(\F_{p^{2}}\) to \(\Fp\)) were tested on random non-superspecial instances by
	comparing zeta functions and by round trips through the dual isogeny; the
	four rules for the codomains of reducible gluings were not tested in this way
	and are checked only through the closure of the graph and the identity of
	Proposition~\ref{prop:graph}(b). The eigenvalues were computed in floating
	point, and the characteristic polynomials of the components were also
	computed exactly and their roots enclosed in intervals of radius at most
		\(1.2\cdot10^{-16}\). The twist-augmented graphs of
	Section~\ref{sec:crypto} were computed from the same data by
	\texttt{twist\_aug.py}, which checks that the twist is an involution on the
	vertices and that the twist of every product is a vertex.
	\item \emph{Table~\ref{tab:examples}.} Each curve was checked from its
	equation alone by \texttt{examples\_table.py}, independently of the census
	program: vanishing of the Cartier--Manin matrix, the Frobenius type from the
	roots of \(f\) in \(\F_{p^{12}}\), the Weil polynomial from the numbers of
	points over \(\Fp\) and \(\F_{p^{2}}\), and the number of Frobenius-stable
	partitions of the Weierstrass points into pairs.
	\item \emph{Comparison with \cite{CorteRealSantosLorenzonReijnders2026}.} For
	every orbit of the enumeration, \texttt{ra\_check.py} computes from the
	sextic the order of the reduced automorphism group of \(C_{\kbar}\), as the
	number of M\"obius transformations preserving the Weierstrass points. For
	the \(37\) primes, every curve with Weil polynomial \(T^{4}+p^{2}\) has
	reduced automorphism group \(\Z/2\times\Z/2\) or \(S_4\), and every curve
	with Weil polynomial \(T^{4}\pm pT^{2}+p^{2}\) has \(S_3\), \(S_4\) or the
	dihedral group of order \(12\), in agreement with
	\cite[Thm.~3.26]{CorteRealSantosLorenzonReijnders2026}.
\end{enumerate}

For \(179\le p\le199\) the census was compared with
Theorem~\ref{thm:vertex} only for the curves with \(\pi^{2}\ne\pm p\) and with
the total mass of Ibukiyama's theorem, and no larger primes were reached. The
scripts of these computations are in the subdirectory \texttt{anc/proofs/graph}
of the supplementary material.

\appendix

\section{Masses of curves with \texorpdfstring{\(\pi^{2}=\pm p\)}{pi\textasciicircum 2 = +-p}}\label{sec:vertexpm}

Parts (a) and (b) of Theorem~\ref{thm:vertex} are proved here as
Theorem~\ref{thm:vertexpm}. We work with
all principally polarised abelian surfaces over \(\Fp\) with
\(\pi^{2}=\eps p\), whether they are Jacobians or not. They are sorted by the
isomorphism class of their \(2\)-adic Tate module with its Frobenius and its
Weil form; Proposition~\ref{prop:twoadic} lists these classes and shows that
on Jacobians they determine the Frobenius type. The correspondence \(D\) of
Lemma~\ref{lem:dual} gives linear relations between the masses of the classes
(Proposition~\ref{prop:neighbours}), and in each connected component of these
relations one mass is computed from Eichler's mass formula. This gives the
mass of every class (Theorem~\ref{thm:celltotals}); the masses of the products
and of the Weil restrictions, known from Proposition~\ref{prop:nonjac}, are
then subtracted.

From the literature we use: Tate's theorem, through the description of a
genus of polarised abelian varieties over a finite field as a double coset
space \cite[\S5]{XueYu2021}; Eichler's mass formula; the Hasse--Minkowski
theorem and the product formula for the Hilbert symbol
\cite[Prop.~14.2.1 and Thm.~14.3.3]{Voight2021}, and Meyer's theorem
\cite[Ch.~IV, \S3.2, Cor.~2]{Serre1973}; the
description of the polarisations of an abelian variety in
\cite[\S21, Application~III, p.~208]{Mumford2008}; and the following statements of
\cite{XueYu2020}: the description of the polarised Dieudonn\'e module
\cite[\S2.6, (2.12)]{XueYu2020} and the fact that the members of a genus of
principally polarised surfaces are principally polarised
\cite[Lemma~2.5]{XueYu2020}, for both signs, and, for \(\eps=+1\), the
structure of the endomorphism algebra \cite[(1.7) and \S3.1.2]{XueYu2020} and
the fact that the surfaces with \(\pi^{2}=p\) form a single isogeny class with
compatible polarisations \cite[Lemma~3.1.3]{XueYu2020}. The table of
Proposition~\ref{prop:neighbours} is obtained by an exact enumeration on
explicit normal forms (Table~\ref{tab:zforms}). The mass formula of
Ibukiyama \cite[Thm.~2]{Ibukiyama1993} is not an input; the sum of our masses
agrees with it (Remark~\ref{rem:vertexlit}).

\subsection*{Notation}
Throughout this appendix \(p\ge7\) is a prime and \(\eps\in\{\pm1\}\). We write
\((A,\theta)\) for a principally polarised abelian surface over \(\Fp\),
\(\dagger\) for its Rosati involution and \(E_\theta\) for its Weil forms. Let
\(\mathcal A_\eps\) be the groupoid of the \((A,\theta)\) with
\(\pi_A^{2}=\eps p\), the morphisms being the isomorphisms over \(\Fp\) that
respect the polarisations. It is a finite groupoid
(Lemma~\ref{lem:Aeps}), and the mass of a full subgroupoid closed under
isomorphism is the sum of \(1/\#\Aut_{\Fp}(A,\theta)\) over its isomorphism
classes.

Let \(\mathbb K=\Q(\sqrt{-p})\), with conjugation \(a\mapsto\bar a\), norm
\(\operatorname N\) and trace \(\operatorname{Tr}\). Its class number is
\(h_{\mathbb K}=h\) if \(p\equiv3\pmod4\) and \(h_{\mathbb K}=h'\) if
\(p\equiv1\pmod4\). Let \(\mathbb K_{+}=\Q(\sqrt p)\). For a prime \(\ell\) let
\(R_\ell=\Z_\ell[T]/(T^{2}-\eps p)\); the class of \(T\) is denoted by
\(\pi\), and \(O_\ell\) is the maximal order of \(R_\ell\otimes\Q\), which is
\(\mathbb K_\ell=\mathbb K\otimes\Q_\ell\) if \(\eps=-1\) and
\(\mathbb K_{+}\otimes\Q_\ell\) if \(\eps=+1\). One has \(R_\ell=O_\ell\)
unless \(\ell=2\) and \(\eps p\equiv1\pmod4\); in that case
\(O_2=\Z_2[\xi]\) with \(\xi=(1+\pi)/2\), and \(R_2=\Z_2+2O_2\). For
\(\eps=-1\) the conjugation, the norm and the trace extend to
\(\mathbb K_\ell\); for \(\eps=+1\), \(\operatorname{Tr}\) denotes the trace
of \(R_\ell\otimes\Q\) over \(\Q_\ell\). We write \(\mathbb A_f\) for the ring
of finite ad\`eles of \(\Q\) and \(\hat\Z\) for the profinite completion of
\(\Z\).

For an abelian variety \(A\) over \(\Fp\) and a prime \(\ell\ne p\),
\(T_\ell A\) is the Tate module. For \(\ell=p\) we write \(T_pA\) for the
Dieudonn\'e module of \(A\): a free \(\Z_p\)-module of rank \(2\dim A\) with
\(\Z_p\)-linear operators \(F\), \(V\) such that \(FV=VF=p\). The Frobenius
endomorphism \(\pi_A\) acts on \(T_pA\) as \(F\) or as \(V\), according to
the convention; if \(\pi_A^{2}=\eps p\), the other operator is
\(p\pi_A^{-1}=\eps\pi_A\). A principal polarisation \(\theta\) induces on
\(T_pA\) a perfect alternating form \(E_\theta\) with
\(E_\theta(Fx,y)=E_\theta(x,Vy)\) \cite[(2.12)]{XueYu2020}. Isomorphisms of
polarised \(\ell\)-divisible groups over \(\Fp\) correspond to isomorphisms of
these modules with their operators and forms, for \(\ell\ne p\) and for
\(\ell=p\); for \(\ell=p\) this is the equivalence between \(p\)-divisible
groups and Dieudonn\'e modules, used in this form in
\cite[\S2.6]{XueYu2020}.

\begin{definition}\label{def:ladic}
	Let \(\ell\) be a prime. A \emph{polarised \(\ell\)-adic module of sign
	\(\eps\)} is a triple \((L,\pi,E)\) consisting of a free \(\Z_\ell\)-module
	\(L\) of rank four, an endomorphism \(\pi\) of \(L\) with \(\pi^{2}=\eps p\)
	and with characteristic polynomial \((T^{2}-\eps p)^{2}\), and a perfect
	alternating form \(E\) on
	\(L\) with \(E(\pi x,y)=\eps E(x,\pi y)\). For \(\ell=2\) we put
	\(\bar L=L/2L\) and \(d(L)=\dim_{\F_2}\ker(\pi-1\mid\bar L)\), and we call
	\(L\) \emph{alternating} if \(E(x,\pi x)\in2\Z_2\) for all \(x\in L\). The
	pair formed by \(d(L)\) and the answer to the question whether \(L\) is
	alternating is the \emph{cell} of \(L\).
\end{definition}

As \(\pi^{2}=\eps p\), such a module is a module over \(R_\ell\), and
\(L\otimes\Q\) is free of rank two over \(R_\ell\otimes\Q\): this is clear if
\(R_\ell\otimes\Q\) is a field, and if it is \(\Q_\ell\times\Q_\ell\), the two
eigenspaces of \(\pi\) on \(L\otimes\Q\) have dimension two by the condition on
the characteristic polynomial. Neither condition on \(\pi\) implies the other
when \(T^{2}-\eps p\) has a root \(s\in\Q_\ell\): for \(\eps=+1\) the scalar
\(\pi=s\) satisfies all the conditions except the one on the characteristic
polynomial, and for \(\eps=-1\) the map
\(\pi=g\oplus(-g^{t})\) on \(\Z_\ell^{2}\oplus\Z_\ell^{2}\), with the
alternating form for which the two summands are isotropic and dual to each
other, \(g=\bigl(\begin{smallmatrix}s&1\\0&s\end{smallmatrix}\bigr)\),
satisfies all the other conditions but not \(\pi^{2}=-p\). For \((A,\theta)\in\mathcal A_\eps\) and every prime
\(\ell\), \(T_\ell(A,\theta)=(T_\ell A,\pi_A,E_\theta)\) is a polarised
\(\ell\)-adic module of sign \(\eps\) (Lemma~\ref{lem:Aeps}). The isomorphism
class of \(T_2(A,\theta)\) is the \emph{\(2\)-adic class} of \((A,\theta)\),
and for a \(2\)-adic class \(X\) we write \(\mathcal A_\eps[X]\) for the full
subgroupoid of the objects of class \(X\) and
\(\Mass(X)=\Mass(\mathcal A_\eps[X])\). For \(L=T_2A\) we have
\(\bar L=A[2]\), \(d(L)=\dim A[2]^{\pi}\), and \(L\) is alternating if and
only if \(e_2(x,\pi x)=1\) for all \(x\in A[2]\); these are the invariants of
Lemma~\ref{lem:tate2}. By Lemma~\ref{lem:cells} below only four cells occur:
\(d=4\) and alternating, \(d=3\) and not alternating, \(d=2\) and not
alternating, \(d=2\) and alternating. We denote them by \(1^{6}\),
\(1^{4}2\), \(1^{2}2^{2}\), \(2^{3}\), the Frobenius types to which they
correspond on Jacobians by Lemma~\ref{lem:tate2}.

For \(\eps=-1\) and a polarised \(\ell\)-adic module \((L,\pi,E)\) we put
\[
H(x,y)=E(\pi x,y)+\pi\,E(x,y)\qquad(x,y\in L\otimes\Q).
\]
By Lemma~\ref{lem:herm} below, \(H\) is a Hermitian form over
\(\mathbb K_\ell\), with values in \(R_\ell\) on \(L\), and \(E\) is the
coefficient of \(\pi\) in \(H\). In a basis \(e_1,e_2\) we write
\(\langle u_1,u_2\rangle\) for the form
\(H(x,y)=u_1x_1\bar y_1+u_2x_2\bar y_2\), and we call the basis
\emph{hyperbolic} if \(H(x,y)=x_1\bar y_2+x_2\bar y_1\).

\subsection*{Statements}

\begin{table}[ht]
	\centering
	\caption{The polarised \(2\)-adic modules of sign \(\eps\)
		(Proposition~\ref{prop:twoadic}). In the last two rows \(\det\) is the
		determinant in the given basis, with values in
		\(R_2\otimes\Q\).}\label{tab:twoadic}
	\renewcommand{\arraystretch}{1.15}
	\begin{tabular}{@{}rllllll@{}}
		\toprule
		\(\eps\) & \(p\bmod4\) & class & module & form & \(d\), alternating & cell\\
		\midrule
		\(-1\) & \(3\) & \(\mathsf{OO}\) & \(O_2e_1\oplus O_2e_2\) & \(H=\langle2,2\rangle\) & \(4\), yes & \(1^{6}\)\\
		\(-1\) & \(3\) & \(\mathsf{RO}\) & \(R_2e_1\oplus O_2e_2\) & \(H=\langle1,2\rangle\) & \(3\), no & \(1^{4}2\)\\
		\(-1\) & \(3\) & \(\mathsf{RR}_{\mathrm n}\) & \(R_2e_1\oplus R_2e_2\) & \(H=\langle1,1\rangle\) & \(2\), no & \(1^{2}2^{2}\)\\
		\(-1\) & \(3\) & \(\mathsf{RR}_{\mathrm a}\) & \(R_2e_1\oplus R_2e_2\) & \(H\) hyperbolic & \(2\), yes & \(2^{3}\)\\
		\(-1\) & \(1\) & \(\mathsf N_{+}\) & \(R_2e_1\oplus R_2e_2\) & \(H=\langle1,1\rangle\) & \(2\), no & \(1^{2}2^{2}\)\\
		\(-1\) & \(1\) & \(\mathsf N_{-}\) & \(R_2e_1\oplus R_2e_2\) & \(H=\langle1,-1\rangle\) & \(2\), no & \(1^{2}2^{2}\)\\
		\(-1\) & \(1\) & \(\mathsf N_{\mathrm a}\) & \(R_2e_1\oplus R_2e_2\) & \(H\) hyperbolic & \(2\), yes & \(2^{3}\)\\
		\(+1\) & \(1\) & \(\mathsf O^{+}\) & \(O_2e_1\oplus O_2e_2\) & \(E=\operatorname{Tr}({\det}/\pi)\) & \(4\), yes & \(1^{6}\)\\
		\(+1\) & \(1,3\) & \(\mathsf R^{+}\) & \(R_2e_1\oplus R_2e_2\) & \(E=\operatorname{Tr}({\det}/2\pi)\) & \(2\), yes & \(2^{3}\)\\
		\bottomrule
	\end{tabular}
\end{table}

\begin{proposition}\label{prop:twoadic}
	\begin{enumerate}[label=\textup{(\alph*)},leftmargin=*]
		\item Every polarised \(2\)-adic module of sign \(\eps\) is isomorphic to
		exactly one of the modules of Table~\textup{\ref{tab:twoadic}} with the
		given \(\eps\) and \(p\bmod4\); these have the indicated cells. In
		particular the cell determines the class, except that for \(\eps=-1\)
		and \(p\equiv1\pmod4\) the classes \(\mathsf N_{+}\) and \(\mathsf N_{-}\)
		have the same cell.
		\item Let \((A,\theta)\in\mathcal A_\eps\). If \(A\) is the Jacobian of a
		curve \(C\), the cell of \(T_2(A,\theta)\) is the Frobenius type of
		\(C\) given in the last column. If \(A\) is a product of supersingular
		elliptic curves of types \textup{I\(\times\)I}, \textup{I\(\times\)II},
		\textup{II\(\times\)II}, its class is \(\mathsf{OO}\), \(\mathsf{RO}\),
		\(\mathsf{RR}_{\mathrm n}\) when \(p\equiv3\pmod4\), and it is
		\(\mathsf N_{+}\) or \(\mathsf N_{-}\) when \(p\equiv1\pmod4\). If \(A\) is a
		Weil restriction, its class is \(\mathsf{RR}_{\mathrm a}\),
		\(\mathsf N_{\mathrm a}\) or \(\mathsf R^{+}\).
	\end{enumerate}
\end{proposition}

For \(\eps=-1\) and \(p\equiv3\pmod4\) the number of classes, four, is the
case \(n=2\) of the classification of Yu \cite[Thm.~1.1 and~1.2]{Yu2012}; for
\(\eps=+1\) the classes are those of \cite[Prop.~3.2.5]{XueYu2020}.

For \((A,\theta)\in\mathcal A_\eps\) and a rational maximal isotropic subgroup
\(K\subset A[2]\), the quotient \(A/K\) with the polarisation induced by
\(2\theta\) is again an object of \(\mathcal A_\eps\), because the quotient
map commutes with the Frobenius endomorphisms.

\begin{table}[ht]
	\centering
	\caption{The numbers \(c(X\to Y)\) of
		Proposition~\ref{prop:neighbours}; \(\#K\) is the number of rational
		maximal isotropic subgroups, and entries equal to zero are
		omitted.}\label{tab:neighbours}
	\renewcommand{\arraystretch}{1.1}
	\begin{tabular}{@{}rllrl@{}}
		\toprule
		\(\eps\) & \(p\) & \(X\) & \(\#K\) & \(c(X\to Y)\)\\
		\midrule
		\(-1\) & \(7\bmod8\) & \(\mathsf{OO}\) & \(15\) &
		\(\mathsf{OO}\colon5\), \(\mathsf{RO}\colon6\),
		\(\mathsf{RR}_{\mathrm n}\colon3\), \(\mathsf{RR}_{\mathrm a}\colon1\)\\
		& & \(\mathsf{RO}\) & \(3\) & \(\mathsf{OO}\colon2\), \(\mathsf{RO}\colon1\)\\
		& & \(\mathsf{RR}_{\mathrm n}\) & \(3\) & \(\mathsf{OO}\colon1\),
		\(\mathsf{RR}_{\mathrm n}\colon1\), \(\mathsf{RR}_{\mathrm a}\colon1\)\\
		& & \(\mathsf{RR}_{\mathrm a}\) & \(7\) & \(\mathsf{OO}\colon1\),
		\(\mathsf{RR}_{\mathrm n}\colon3\), \(\mathsf{RR}_{\mathrm a}\colon3\)\\
		\(-1\) & \(3\bmod8\) & \(\mathsf{OO}\) & \(15\) & \(\mathsf{OO}\colon3\),
		\(\mathsf{RR}_{\mathrm n}\colon9\), \(\mathsf{RR}_{\mathrm a}\colon3\)\\
		& & \(\mathsf{RO}\) & \(3\) & \(\mathsf{RO}\colon3\)\\
		& & \(\mathsf{RR}_{\mathrm n}\) & \(3\) & \(\mathsf{OO}\colon1\),
		\(\mathsf{RR}_{\mathrm n}\colon1\), \(\mathsf{RR}_{\mathrm a}\colon1\)\\
		& & \(\mathsf{RR}_{\mathrm a}\) & \(7\) & \(\mathsf{OO}\colon1\),
		\(\mathsf{RR}_{\mathrm n}\colon3\), \(\mathsf{RR}_{\mathrm a}\colon3\)\\
		\(-1\) & \(1\bmod4\) & \(\mathsf N_{+}\) & \(3\) & \(\mathsf N_{+}\colon3\)\\
		& & \(\mathsf N_{-}\) & \(3\) & \(\mathsf N_{-}\colon1\), \(\mathsf N_{\mathrm a}\colon2\)\\
		& & \(\mathsf N_{\mathrm a}\) & \(7\) & \(\mathsf N_{-}\colon3\), \(\mathsf N_{\mathrm a}\colon4\)\\
		\(+1\) & \(1\bmod8\) & \(\mathsf O^{+}\) & \(15\) & \(\mathsf O^{+}\colon9\), \(\mathsf R^{+}\colon6\)\\
		& & \(\mathsf R^{+}\) & \(7\) & \(\mathsf O^{+}\colon1\), \(\mathsf R^{+}\colon6\)\\
		\(+1\) & \(5\bmod8\) & \(\mathsf O^{+}\) & \(15\) & \(\mathsf O^{+}\colon5\), \(\mathsf R^{+}\colon10\)\\
		& & \(\mathsf R^{+}\) & \(7\) & \(\mathsf O^{+}\colon1\), \(\mathsf R^{+}\colon6\)\\
		\(+1\) & \(3\bmod4\) & \(\mathsf R^{+}\) & \(7\) & \(\mathsf R^{+}\colon7\)\\
		\bottomrule
	\end{tabular}
\end{table}

\begin{proposition}\label{prop:neighbours}
	Let \(X\), \(Y\) be \(2\)-adic classes of sign \(\eps\).
	\begin{enumerate}[label=\textup{(\alph*)},leftmargin=*]
		\item For \((A,\theta)\in\mathcal A_\eps\) of class \(X\), the number
		\(c(X\to Y)\) of rational maximal isotropic subgroups \(K\subset A[2]\)
		such that \(A/K\), with the polarisation induced by \(2\theta\), is of
		class \(Y\) depends only on \(X\) and \(Y\). It is given by
		Table~\textup{\ref{tab:neighbours}}.
		\item \(\Mass(X)\,c(X\to Y)=\Mass(Y)\,c(Y\to X)\).
	\end{enumerate}
\end{proposition}

\begin{theorem}\label{thm:celltotals}
	Let \(p\ge7\). The masses of the \(2\)-adic classes of the principally
	polarised abelian surfaces over \(\Fp\) with \(\pi^{2}=\eps p\) are the
	following.
	\begin{enumerate}[label=\textup{(\roman*)},leftmargin=*]
		\item \(\eps=-1\), \(p\equiv7\pmod8\):
		\[
		\Mass(\mathsf{OO})=\Mass(\mathsf{RR}_{\mathrm a})=\frac{(p-1)h}{48},
		\qquad
		\Mass(\mathsf{RO})=\Mass(\mathsf{RR}_{\mathrm n})=\frac{3(p-1)h}{48}.
		\]
		\item \(\eps=-1\), \(p\equiv3\pmod8\):
		\[
		\begin{aligned}
		\Mass(\mathsf{OO})&=\frac{(p-1)h}{48}, &
		\Mass(\mathsf{RO})&=\frac{(p+1)h}{16},\\
		\Mass(\mathsf{RR}_{\mathrm n})&=\frac{9(p-1)h}{48}, &
		\Mass(\mathsf{RR}_{\mathrm a})&=\frac{3(p-1)h}{48}.
		\end{aligned}
		\]
		\item \(\eps=-1\), \(p\equiv1\pmod4\):
		\[
		\Mass(\mathsf N_{+})=\frac{(p+1)h'}{32},\qquad
		\Mass(\mathsf N_{-})=\frac{(p-1)h'}{32},\qquad
		\Mass(\mathsf N_{\mathrm a})=\frac{(p-1)h'}{48}.
		\]
		\item \(\eps=+1\): if \(p\equiv3\pmod4\), then
		\(\Mass(\mathsf R^{+})=z_p/4\); if \(p\equiv1\pmod4\), then
		\[
		\Mass(\mathsf O^{+})=\frac{z_p}{4},\qquad
		\Mass(\mathsf R^{+})=\frac{\bigl(8-2\leg2p\bigr)z_p}{4}.
		\]
	\end{enumerate}
\end{theorem}

\begin{theorem}\label{thm:vertexpm}
	Parts \textup{(a)} and \textup{(b)} of
	Theorem~\textup{\ref{thm:vertex}} hold: for every prime \(p\ge7\) the
	masses of the superspecial curves of genus two over \(\Fp\) with Weil
	polynomial \((T^{2}+p)^{2}\) or \((T^{2}-p)^{2}\), by Frobenius type, are
	those stated there, and the types that are not listed have mass zero.
\end{theorem}

\begin{corollary}\label{cor:split1222}
	Let \(p\equiv1\pmod4\) and let \(C/\Fp\) be a superspecial curve of genus
	two with Weil polynomial \((T^{2}+p)^{2}\) and Frobenius type
	\(1^{2}2^{2}\). Then either the quotients of \(\Jac(C)\) by its three
	rational maximal isotropic subgroups all have cell \(1^{2}2^{2}\), or one
	of them has cell \(1^{2}2^{2}\) and two have cell \(2^{3}\). The curves of
	the first kind have mass \((p+1)h'/32-h'^{2}/16\), and those of the second
	kind have mass \((p-1)h'/32-h'^{2}/16\).
\end{corollary}

Here the cell of a quotient which is not a Jacobian is that of
Proposition~\ref{prop:twoadic}(b): \(1^{2}2^{2}\) for a product of two
elliptic curves over \(\Fp\) and \(2^{3}\) for a Weil restriction.

\subsection*{Two-adic classes}

\begin{lemma}\label{lem:Aeps}
	Let \((A,\theta)\in\mathcal A_\eps\).
	\begin{enumerate}[label=\textup{(\alph*)},leftmargin=*]
		\item The characteristic polynomial of \(\pi_A\) is
		\((T^{2}-\eps p)^{2}\), \(\pi_A^{\dagger}=\eps\pi_A\), and \(A\) is
		superspecial. For every prime \(\ell\), \(T_\ell(A,\theta)\) is a
		polarised \(\ell\)-adic module of sign \(\eps\).
		\item \((A,\theta)\) is exactly one of the following: the Jacobian of a
		superspecial curve \(C/\Fp\) of genus two with Weil polynomial
		\((T^{2}-\eps p)^{2}\), which is unique up to \(\Fp\)-isomorphism and
		satisfies \(\Aut_{\Fp}(C)\cong\Aut_{\Fp}(A,\theta)\); a product of two
		supersingular elliptic curves over \(\Fp\) with the product
		polarisation, which occurs only for \(\eps=-1\); the Weil restriction
		of a supersingular elliptic curve \(B/\F_{p^{2}}\) with
		\(\pi_B=\eps p\). Conversely the Jacobian of every superspecial curve
		with Weil polynomial \((T^{2}-\eps p)^{2}\) is an object of
		\(\mathcal A_\eps\).
		\item For \(\eps=-1\) the products of types \textup{I\(\times\)I},
		\textup{I\(\times\)II}, \textup{II\(\times\)II} have masses
		\(N_1^{2}/8\), \(N_1N_2/4\), \(N_2^{2}/8\). For each sign the Weil
		restrictions have mass \((p-1)/48\).
	\end{enumerate}
\end{lemma}

\begin{proof}
	(a) The first assertion is Lemma~\ref{lem:irred}. As in the proof of
	Theorem~\ref{thm:local}(i), \(\pi_A^{\dagger}\pi_A=p\)
	\cite[II, Lemma~1.2]{MilneAV}, hence
	\(\pi_A^{\dagger}=p\pi_A^{-1}=\eps\pi_A\). On the Dieudonn\'e module,
	\(F^{2}=\pi_A^{2}=\eps p\), so \(F\,T_pA=pF^{-1}T_pA=V\,T_pA\); hence
	\(A_{\kbar}\) has \(a\)-number two, and \(A\) is superspecial by the
	theorem of Oort \cite[Thm.~2 and Remark~3, p.~35]{Oort1975}, in the form
	quoted in \cite[\S3.2.1]{XueYu2020} and
	\cite[proof of Prop.~5.1, p.~1623]{XueYangYu2016}. For \(\ell\ne p\) the
	form \(E_\theta\) on
	\(T_\ell A\) is perfect and alternating and
	\(E_\theta(\pi x,y)=E_\theta(x,\pi^{\dagger}y)=\eps E_\theta(x,\pi y)\);
	for \(\ell=p\) the relation \(E_\theta(Fx,y)=E_\theta(x,Vy)\) says the
	same, because \(\{F,V\}=\{\pi_A,\eps\pi_A\}\). On \(T_\ell A\) we have
	\(\pi_A^{2}=\eps p\), and the characteristic polynomial of \(\pi_A\) is
	that of \(\pi_A\).

	(b) Over \(\kbar\), \((A,\theta)\) is a Jacobian or a polarised product,
	and not both \cite[Thm.~1.3]{HoweNartRitzenthaler2009}. In the first case
	it is the Jacobian of a curve \(C\) over \(\Fp\), unique up to
	isomorphism, and \(\Aut_{\Fp}(C)\cong\Aut_{\Fp}(A,\theta)\)
	(Lemmas~\ref{lem:dual} and~\ref{lem:torelli}); \(C\) is superspecial by
	(a). Conversely, if \(C\) is superspecial with Weil polynomial
	\((T^{2}-\eps p)^{2}\), the automorphism \(u\) of
	Theorem~\ref{thm:local}(i) has finite order and characteristic polynomial
	\((T-\eps)^{4}\), so \(u=\eps\) and \(\pi^{2}=\eps p\). In the second case
	we are in the situation described after Lemma~\ref{lem:prod}: either
	\(A=E_1\times E_2\) over \(\Fp\), or \(A\) is the Weil restriction of an
	elliptic curve \(B/\F_{p^{2}}\), and then
	\(\pi_B=\pi_A^{2}|_B=\eps p\), so that \(B\) is supersingular. If
	\(A=E_1\times E_2\), then \(\pi_{E_i}^{2}=\eps p\); as
	\(\pi_{E_i}^{2}-a_i\pi_{E_i}+p=0\) with \(a_i\in\Z\), we get
	\(a_i\pi_{E_i}=(\eps+1)p\). For \(\eps=+1\) this would make
	\(\pi_{E_i}=2p/a_i\) a rational number of square \(p\). Hence
	\(\eps=-1\) and \(a_i=0\). Finally, \(\mathcal A_\eps\) has finitely
	many isomorphism classes, because there are finitely many curves of genus
	two over \(\Fp\) and finitely many elliptic curves over \(\Fp\) and over
	\(\F_{p^{2}}\).

	(c) The mass of the Weil restrictions is given by
	Proposition~\ref{prop:nonjac}. For the products we repeat the argument
	given there, keeping track of the types. Write \(N_{\mathrm I}=N_1\),
	\(N_{\mathrm{II}}=N_2\), and let \(a,b\) be two types. The ordered pairs
	\((E_1,E_2)\) of supersingular elliptic curves over \(\Fp\) whose types
	are \(a\) and \(b\), in some order, form a groupoid of mass
	\((N_a/2)^{2}\) if \(a=b\) and \(2(N_a/2)(N_b/2)\) if \(a\ne b\), because
	\(\Aut_{\Fp}(E)=\{\pm1\}\) (Lemma~\ref{lem:ssp}). The functor
	\((E_1,E_2)\mapsto E_1\times E_2\) to the products with factors of types
	\(a\) and \(b\) is faithful and its fibres have two elements
	(Lemma~\ref{lem:prod}), so by Lemma~\ref{lem:fibre} these products have
	half that mass.
\end{proof}

\begin{lemma}\label{lem:cells}
	Let \((L,\pi,E)\) be a polarised \(2\)-adic module of sign \(\eps\) and
	\(t=\pi-1\) on \(\bar L\). Then \(t^{2}=0\) and \(d(L)\in\{2,3,4\}\);
	\(L\) is alternating if \(d(L)=4\), and is not alternating if \(d(L)=3\).
	If \(\eps=+1\), then \(E(x,\pi x)=0\) for all \(x\in L\) and
	\(d(L)\in\{2,4\}\).
\end{lemma}

\begin{proof}
	On \(\bar L\) we have \(\pi^{2}=\eps p=1\) and \(2\pi=0\), so
	\(t^{2}=\pi^{2}-2\pi+1=0\). Hence the image of \(t\) lies in its kernel,
	the rank of \(t\) is at most two, and \(d(L)=4-\operatorname{rank}t\). Let
	\(e\) be the reduction of \(E\) and \(b(x,y)=e(x,ty)\). From
	\(E(\pi x,y)=\eps E(x,\pi y)\) we get \(e(tx,y)=e(x,ty)\), hence
	\(b(y,x)=e(y,tx)=e(tx,y)=b(x,y)\). So \(b\) is a symmetric bilinear form;
	its radical is \(\ker t\), because \(e\) is non-degenerate, and
	\(b(x,x)=e(x,\pi x)-e(x,x)=e(x,\pi x)\). If \(L\) is alternating, \(b\)
	induces a non-degenerate alternating form on \(\bar L/\ker t\), whose
	dimension \(4-d(L)\) is therefore even. If \(d(L)=4\), then \(b=0\).
	Finally, if \(\eps=+1\), then \(E(x,\pi x)=E(\pi x,x)=-E(x,\pi x)\).
\end{proof}

\begin{lemma}\label{lem:nonjaccells}
	Let \((A,\theta)\in\mathcal A_\eps\) and \(L=T_2A\). If \(A\) is a product
	of elliptic curves of types \textup{I\(\times\)I},
	\textup{I\(\times\)II}, \textup{II\(\times\)II}, then \(d(L)=4,3,2\), and
	\(L\) is alternating only in the first case. If \(A\) is a Weil
	restriction, then \(d(L)=2\) and \(L\) is alternating.
\end{lemma}

\begin{proof}
	Let \(E/\Fp\) be a supersingular elliptic curve. If it is of type~I, then
	\(\pi_E=1\) on \(E[2]\) and \(e_2(x,\pi_Ex)=e_2(x,x)=1\). If it is of
	type~II, then \(t=\pi_E-1\) is non-zero with \(t^{2}=0\) on \(E[2]\), so
	its kernel is a line; for \(x\) outside this line, \(tx\ne0\) and
	\(tx\ne x\) (otherwise \(\pi_Ex=0\)), hence
	\(e_2(x,\pi_Ex)=e_2(x,tx)=-1\). The Tate module of \(E_1\times E_2\) is
	the orthogonal sum of those of the factors and \(\pi\) acts diagonally;
	this gives the first assertion. Let \(A\) be a Weil restriction. With the
	notation introduced after Lemma~\ref{lem:prod}, \(T_2A=T_2B\oplus T_2B'\)
	orthogonally and \(\pi_A(x,z)=(\varphi'z,\varphi x)\), with
	\(\varphi'\varphi=\pi_B=\eps p\). On \(A[2]\) the point \((x,z)\) is
	fixed if and only if \(z=\varphi x\), so \(d(L)=2\). The relation
	\(E_\theta(\pi_A(x,0),(0,z))=\eps E_\theta((x,0),\pi_A(0,z))\) reads
	\(E_{B'}(\varphi x,z)=\eps E_B(x,\varphi'z)\). Therefore
	\[
	E_\theta\bigl((x,z),\pi_A(x,z)\bigr)=E_B(x,\varphi'z)+E_{B'}(z,\varphi x)
	=(1-\eps)\,E_{B'}(z,\varphi x)\in2\Z_2 .\qedhere
	\]
\end{proof}

\begin{lemma}\label{lem:herm}
	Let \(\eps=-1\), let \(\ell\) be a prime and \((L,\pi,E)\) a polarised
	\(\ell\)-adic module of sign \(-1\), \(V=L\otimes\Q\). The form
	\(H(x,y)=E(\pi x,y)+\pi E(x,y)\) on \(V\) is linear in \(x\), satisfies
	\(H(y,x)=\overline{H(x,y)}\) and \(H(x,x)=-E(x,\pi x)\), takes values in
	\(R_\ell\) on \(L\), and
	\[
	L=L^{\#}:=\{x\in V: H(x,L)\subset R_\ell\}.
	\]
	The automorphisms of \((L,\pi,E)\) are the \(R_\ell\)-linear automorphisms
	of \(L\) preserving \(H\). If \(L\) is free over \(R_\ell\), with a basis
	of Gram matrix \(\Xi\), then \(\det\Xi\in\Z_\ell^{\times}\).
\end{lemma}

\begin{proof}
	Using \(\pi^{2}=-p\) and \(E(x,\pi y)=-E(\pi x,y)\),
	\[
	\begin{aligned}
	H(\pi x,y)&=-pE(x,y)+\pi E(\pi x,y)=\pi H(x,y),\\
	H(y,x)&=-E(x,\pi y)-\pi E(x,y)=E(\pi x,y)-\pi E(x,y)=\overline{H(x,y)} .
	\end{aligned}
	\]
	An element \(x\in V\) lies in \(L^{\#}\) if and only if \(E(x,y)\) and
	\(E(\pi x,y)=-E(x,\pi y)\) are in \(\Z_\ell\) for all \(y\in L\), that is,
	as \(\pi L\subset L\), if and only if \(E(x,L)\subset\Z_\ell\); and \(E\)
	is perfect. Since \(E\) is the coefficient of \(\pi\) in \(H\), a
	\(\Z_\ell\)-linear map commuting with \(\pi\) preserves \(E\) if and only
	if it preserves \(H\). If \(L\) is free, \(L^{\#}=L\) means that
	\(\Xi\) is invertible over \(R_\ell\), and \(\det\Xi\) is fixed by
	the conjugation.
\end{proof}

A lattice with \(L^{\#}=L\) is called \emph{unimodular}. By the lemma, the
polarised \(\ell\)-adic modules of sign \(-1\) are the unimodular
\(R_\ell\)-lattices in Hermitian planes over \(\mathbb K_\ell\), and for
\(\ell=2\) such a lattice is alternating if and only if
\(H(x,x)\in2\Z_2\) for all \(x\).

\begin{lemma}\label{lem:norms}
	Let \(\eps=-1\) and \(\ell=2\). The ring \(R_2\) is local, with
	\(R_2/2R_2=\F_2[t]/(t^{2})\), \(t=\pi-1\), and its units are the
	\(a+b\pi\) with \(a+b\) odd.
	\begin{enumerate}[label=\textup{(\alph*)},leftmargin=*]
		\item If \(p\equiv3\pmod4\), then
		\(\operatorname N(R_2^{\times})=\Z_2^{\times}\); the ring \(O_2\) is
		\'etale over \(\Z_2\), and an \(O_2\)-submodule of \(\mathbb K_2\)
		contained in \(R_2\) is contained in \(2O_2\).
		\item If \(p\equiv1\pmod4\), then \(R_2\) is a discrete valuation ring
		with uniformiser \(\varpi=\pi-1\), \(\operatorname N(\varpi)=1+p\)
		is twice an odd number, \(\operatorname N(R_2^{\times})=1+4\Z_2\), and
		an element of \(\mathbb K_2^{\times}\) whose norm is a unit is a unit.
	\end{enumerate}
\end{lemma}

\begin{proof}
	\(R_2/2R_2=\F_2[\pi]/(\pi-1)^{2}\). We have
	\(\operatorname N(a+b\pi)=a^{2}+pb^{2}\), and
	\(\operatorname N(R_2^{\times})\) is a subgroup of \(\Z_2^{\times}\)
	containing \((\Z_2^{\times})^{2}=1+8\Z_2\).
	(a) The norms of \(1\), \(1+2\pi\), \(\pi\), \(2+\pi\) are \(1\),
	\(1+4p\), \(p\), \(4+p\), which represent the four classes of units
	modulo \(8\). The discriminant of
	\(\xi^{2}-\xi+(1+p)/4=0\) is \(-p\), a unit. If \(x=a+2b\xi\) with
	\(a,b\in\Z_2\) and \(x\xi\in R_2=\Z_2+2\Z_2\xi\), then, as
	\(x\xi=-b(1+p)/2+(a+2b)\xi\), \(a\) is even and \(x\in2O_2\).
	(b) The minimal polynomial of \(\varpi\) is \(X^{2}+2X+(1+p)\), an
	Eisenstein polynomial. For \(a+b\) odd, \(a^{2}+pb^{2}\equiv1\pmod4\), and
	the norms \(1\) and \(1+4p\equiv5\pmod8\) occur.
\end{proof}

\begin{lemma}\label{lem:unram}
	Let \(\eps=-1\), let \(\ell\) be a prime that does not divide \(p\), with
	\(\ell\ne2\) if \(p\equiv1\pmod4\), and let \(M\) be a free
	\(O_\ell\)-module of rank two with an \(O_\ell\)-valued Hermitian form
	\(H\) such that \(M=\{x\in M\otimes\Q:H(x,M)\subset O_\ell\}\). Then
	\(M\) has a hyperbolic basis. In particular \(M\) is unique up to
	isometry, and it has a basis in which \(H=\langle1,1\rangle\).
\end{lemma}

\begin{proof}
	Under the hypothesis \(O_\ell\) is \'etale over \(\Z_\ell\). If
	\(O_\ell=\Z_\ell\times\Z_\ell\), with the conjugation exchanging the
	factors, then \(M=M_1\times M_2\) and \(H\) is given by a perfect pairing
	\(\beta\colon M_1\times M_2\to\Z_\ell\); with dual bases \(u_1,u_2\) of
	\(M_1\) and \(u_1^{*},u_2^{*}\) of \(M_2\), the vectors
	\(f_1=(u_1,u_2^{*})\), \(f_2=(u_2,u_1^{*})\) form a hyperbolic basis.
	Otherwise \(O_\ell\) is the ring of integers of the unramified quadratic
	extension of \(\Q_\ell\). The trace \(O_\ell\to\Z_\ell\) is surjective,
	and so is the norm \(O_\ell^{\times}\to\Z_\ell^{\times}\): it is
	surjective on the residue fields, and
	\(\operatorname N(1+\ell^{k}x)\equiv1+\ell^{k}\operatorname{Tr}(x)
	\pmod{\ell^{k+1}}\) for \(k\ge1\). If \(H(x,x)\in\ell\Z_\ell\) for all
	\(x\in M\), then \(\operatorname{Tr}H(x,y)\in\ell\Z_\ell\) for all
	\(x,y\), which contradicts \(H(M,M)=O_\ell\). So there is \(x\) with
	\(H(x,x)\in\Z_\ell^{\times}\); then \(M=O_\ell x\oplus O_\ell y\)
	orthogonally with \(H(y,y)\in\Z_\ell^{\times}\), and after multiplying
	\(x\) and \(y\) by units, \(H=\langle1,1\rangle\) in the basis
	\(e_1=x\), \(e_2=y\). Let \(a\in O_\ell^{\times}\) with
	\(\operatorname N(a)=-1\). Then \(f_1=e_1+ae_2\) is isotropic and part of
	a basis, so there is \(y'\) with \(H(f_1,y')=1\), and
	\(f_2=y'+cf_1\) with \(\operatorname{Tr}(c)=-H(y',y')\) satisfies
	\(H(f_2,f_2)=0\), \(H(f_1,f_2)=1\). The Gram determinant of \((f_1,f_2)\)
	is \(-1\), so it is a basis. Conversely, in both cases, if \(f_1,f_2\) is
	a hyperbolic basis and \(\operatorname{Tr}(c)=1\), the vectors
	\(f_1+cf_2\) and \(f_1-\bar cf_2\) are orthogonal, of norms \(1\) and
	\(-1\), and they form a basis; and
	\(\langle1,-1\rangle\cong\langle1,1\rangle\) because \(-1\) is the norm of
	a unit.
\end{proof}

\begin{lemma}\label{lem:altR}
	Let \(\eps=-1\) and let \(L\) be a free \(R_2\)-module of rank two with
	an \(R_2\)-valued Hermitian form \(H\), unimodular and such that
	\(H(x,x)\in2\Z_2\) for all \(x\). Then \(L\) has a hyperbolic basis.
\end{lemma}

\begin{proof}
	In a basis \(e_1,e_2\) the Gram matrix is
	\(\bigl(\begin{smallmatrix}2\alpha&\beta\\ \bar\beta&2\gamma
	\end{smallmatrix}\bigr)\) with \(\alpha,\gamma\in\Z_2\), and
	\(4\alpha\gamma-\operatorname N(\beta)\) is a unit by
	Lemma~\ref{lem:herm}. So \(\beta\in R_2^{\times}\), and replacing \(e_2\)
	by \(\bar\beta^{-1}e_2\) we may assume \(\beta=1\). For
	\(x=(a+b\pi)e_1+e_2\) with \(a,b\in\Z_2\),
	\[
	H(x,x)=2\alpha\operatorname N(a+b\pi)+\operatorname{Tr}(a+b\pi)+2\gamma
	=2\bigl(\alpha a^{2}+a+\gamma+\alpha pb^{2}\bigr).
	\]
	Take \(b=0\) if \(\alpha\) is even, and \(b\in\{0,1\}\) with
	\(b\equiv\gamma\pmod2\) if \(\alpha\) is odd. The polynomial
	\(g(a)=\alpha a^{2}+a+\gamma+\alpha pb^{2}\) has derivative
	\(2\alpha a+1\), a unit, and
	\(g(a)\equiv(\alpha+1)a+\gamma+\alpha b\pmod2\) has a root modulo \(2\)
	by the choice of \(b\). By Hensel's lemma \(g\) has a root in \(\Z_2\).
	This gives a vector \(x\) with \(H(x,x)=0\) which is part of a basis. As
	\(L\) is unimodular there is \(y\) with \(H(x,y)=1\), and \(y+cx\) with
	\(\operatorname{Tr}(c)=-H(y,y)\) has norm zero; such a \(c\in R_2\) exists
	because \(H(y,y)\in2\Z_2=\operatorname{Tr}(R_2)\).
\end{proof}

\begin{lemma}\label{lem:sympl}
	Let \(\eps=+1\), let \(\ell\) be a prime and \((L,\pi,E)\) a polarised
	\(\ell\)-adic module of sign \(+1\), \(V=L\otimes\Q\),
	\(k=R_\ell\otimes\Q\). There is a unique \(k\)-bilinear alternating form
	\(\Psi\colon V\times V\to k\) with \(E=\operatorname{Tr}\circ\Psi\). Let
	\(\mathfrak o\) be \(R_\ell\) or \(O_\ell\), and let \(\gamma\) be a
	generator of the dual \(\mathfrak o^{\vee}\) of \(\mathfrak o\) for the
	trace form. If \(L\) is a free \(\mathfrak o\)-module of rank two, it has an
	\(\mathfrak o\)-basis in which \(\Psi=\gamma\det\); in particular
	\((L,\pi,E)\) is determined up to isomorphism by \(\mathfrak o\), and its
	automorphisms are the elements of determinant one of
	\(\operatorname{GL}_2(\mathfrak o)\). One has
	\(R_\ell^{\vee}=(2\pi)^{-1}R_\ell\), and
	\(O_2^{\vee}=\pi^{-1}O_2\) if \(p\equiv1\pmod4\).
\end{lemma}

\begin{proof}
	The trace form of the \'etale algebra \(k\) is non-degenerate, so for
	\(x,y\in V\) there is a unique \(\Psi(x,y)\in k\) with
	\(E(ax,y)=\operatorname{Tr}(a\Psi(x,y))\) for all \(a\in k\). It is
	\(k\)-bilinear because \(E(\pi x,y)=E(x,\pi y)\), and
	\(E(ax,x)=E(x,ax)=-E(ax,x)\) gives \(\Psi(x,x)=0\). The module
	\(\mathfrak o^{\vee}\) is generated by one element: this is shown below
	for \(R_\ell\), and \(O_\ell\) is a discrete valuation ring or a product
	of two. In an
	\(\mathfrak o\)-basis of \(L\) we have \(\Psi=c\det\) with
	\(c\in k^{\times}\). For \(x\in V\) the set \(\det(x,L)\) is the
	\(\mathfrak o\)-module \(x_1\mathfrak o+x_2\mathfrak o\), and for an
	\(\mathfrak o\)-module \(I\subset k\),
	\(\operatorname{Tr}(cI)\subset\Z_\ell\) if and only if
	\(cI\subset\mathfrak o^{\vee}\). Hence the dual of \(L\) for \(E\) is
	\(c^{-1}\gamma L\), and \(E\) is perfect if and only if \(c=u\gamma\)
	with \(u\in\mathfrak o^{\times}\); replacing \(e_1\) by \(u^{-1}e_1\)
	gives \(c=\gamma\). An \(\mathfrak o\)-linear automorphism \(g\) multiplies
	\(\det\) by \(\det g\). Finally \(1,\pi\) is a basis of \(R_\ell\) with
	\(\operatorname{Tr}(1)=2\), \(\operatorname{Tr}(\pi)=0\),
	\(\operatorname{Tr}(\pi^{2})=2p\), whose dual basis is \(1/2\),
	\(1/(2\pi)\); and the different of \(O_2=\Z_2[\xi]\) is generated by
	\(\xi-(1-\xi)=\pi\).
\end{proof}

\begin{proof}[Proof of Proposition~\textup{\ref{prop:twoadic}}]
	(a) Let \((L,\pi,E)\) be a polarised \(2\)-adic module of sign \(\eps\).
	As a module over \(\F_2[t]/(t^{2})\), \(\bar L\) is a direct sum of
	\(a\) copies of \(\F_2[t]/(t^{2})\) and \(b\) copies of \(\F_2\), with
	\(2a+b=4\) and \(a+b=d(L)\).

	Let \(\eps=-1\) and \(p\equiv3\pmod4\). By
	Lemma~\ref{lem:cells} there are three cases.

	\(d(L)=4\). Then \((\pi-1)L\subset2L\), so \(\xi\) acts on \(L\), and
	\(L\) is a torsion-free module over \(O_2\) with \(L\otimes\Q\) free of
	rank two; hence \(L\) is free of rank two over \(O_2\). For
	\(x\in L\) the set \(H(x,L)\) is an \(O_2\)-module contained in \(R_2\),
	hence contained in \(2O_2\) (Lemma~\ref{lem:norms}); and for
	\(x\in L\otimes\Q\), \(H(x,L)\subset R_2\) if and only if
	\(H(x,L)\subset2O_2\). So \(H/2\) is \(O_2\)-valued and \(L\) is
	unimodular for it. By Lemma~\ref{lem:unram}, \(H/2=\langle1,1\rangle\) in
	a suitable basis. This is the class \(\mathsf{OO}\); it is alternating.

	\(d(L)=2\). Then \(\bar L\) is free of rank two over \(R_2/2R_2\); by
	Nakayama's lemma \(L\) is generated by two elements over \(R_2\), and as
	it has rank four over \(\Z_2\) it is free of rank two over \(R_2\).
	Suppose that \(L\) is not alternating, and let \(x\in L\) with
	\(H(x,x)\in\Z_2^{\times}\). The map
	\(y\mapsto y-H(y,x)H(x,x)^{-1}x\) is a projection of \(L\) onto
	\(M=\{y\in L:H(y,x)=0\}\), so \(L=R_2x\oplus M\) orthogonally. Then
	\(M\) is unimodular, and projective of rank one over the local ring
	\(R_2\), so \(M=R_2y\) with \(H(y,y)\in\Z_2^{\times}\). By
	Lemma~\ref{lem:norms}(a) we can multiply \(x\) and \(y\) by units so that
	\(H=\langle1,1\rangle\): class \(\mathsf{RR}_{\mathrm n}\). If \(L\) is
	alternating, Lemma~\ref{lem:altR} gives the class
	\(\mathsf{RR}_{\mathrm a}\).

	\(d(L)=3\). By Lemma~\ref{lem:cells}, \(L\) is not alternating, and as
	before \(L=R_2x\oplus M\) orthogonally with \(H(x,x)\in\Z_2^{\times}\).
	The kernel of \(t\) on \(R_2x/2R_2x\) is a line; as \(d(L)=3\), \(t=0\)
	on the two-dimensional space \(M/2M\). Hence \(M\) is an \(O_2\)-module
	with \(M\otimes\Q\) free of rank one over \(\mathbb K_2\), so
	\(M=O_2y\). Let \(c=H(y,y)\in\Z_2\). Then
	\(M^{\#}=\{sy:scO_2\subset R_2\}=2c^{-1}O_2y\) by
	Lemma~\ref{lem:norms}(a), and \(M^{\#}=M\) forces
	\(c\in2\Z_2^{\times}\). Multiplying \(x\) by a unit of \(R_2\) and \(y\)
	by a unit of \(O_2\) we get \(H(x,x)=1\), \(H(y,y)=2\): class
	\(\mathsf{RO}\).

	Let \(\eps=-1\) and \(p\equiv1\pmod4\). By
	Lemma~\ref{lem:norms}(b), \(L\) is free of rank two over the discrete
	valuation ring \(R_2\), so \(\bar L\) is free over \(R_2/2R_2\) and
	\(d(L)=2\). If \(L\) is not alternating, the argument above gives
	\(H=\langle u_1,u_2\rangle\) with \(u_i\in\Z_2^{\times}\), and by
	Lemma~\ref{lem:norms}(b) we may take \(u_i\in\{1,-1\}\). In
	\(\langle1,1\rangle\) the vector \(v=e_1+\varpi e_2\) has
	\(H(v,v)=2+p\equiv3\pmod4\); splitting it off we get
	\(\langle1,1\rangle\cong\langle u,u'\rangle\) with \(u\equiv3\pmod4\),
	and \(uu'\equiv1\pmod4\) because the determinant of a Gram matrix is well
	defined modulo \(\operatorname N(R_2^{\times})\). Hence
	\(\langle1,1\rangle\cong\langle-1,-1\rangle\), whereas
	\(\langle1,1\rangle\) and \(\langle1,-1\rangle\) are not isometric. These
	are the classes \(\mathsf N_{+}\) and \(\mathsf N_{-}\). If \(L\) is
	alternating, Lemma~\ref{lem:altR} gives the class
	\(\mathsf N_{\mathrm a}\). The class
	\(\mathfrak d(L)\) of \(\det H\) in
	\(\Q_2^{\times}/\operatorname N(\mathbb K_2^{\times})\) is an invariant
	of the Hermitian space \(L\otimes\Q\), and
	\(\mathfrak d(\mathsf N_{+})=1\),
	\(\mathfrak d(\mathsf N_{-})=\mathfrak d(\mathsf N_{\mathrm a})=-1\ne1\) by
	Lemma~\ref{lem:norms}(b).

	Let \(\eps=+1\). By Lemma~\ref{lem:cells},
	\(d(L)\in\{2,4\}\) and \(L\) is alternating. If \(d(L)=2\), then \(L\) is
	free of rank two over \(R_2\) as above, and Lemma~\ref{lem:sympl} gives
	the class \(\mathsf R^{+}\). If \(d(L)=4\), then \((\pi-1)/2\) acts on
	\(L\); its minimal polynomial \(X^{2}+X+(1-p)/4\) has coefficients in
	\(\Z_2\), so \(p\equiv1\pmod4\), \(L\) is a free module of rank two over
	\(O_2\), and Lemma~\ref{lem:sympl} gives the class \(\mathsf O^{+}\).

	In each case the cells of the normal forms are those of
	Table~\ref{tab:twoadic}: \(d\) is read off from \(R_2/2R_2\) and
	\(O_2/2O_2\), on which \(t\) has a kernel of dimension one and two, and
	for \(\eps=-1\) the form \(\langle u_1,u_2\rangle\) takes the values
	\(H(x,x)=u_1\operatorname N(x_1)+u_2\operatorname N(x_2)\), while in a
	hyperbolic basis over \(R_2\),
	\(H(x,x)=\operatorname{Tr}(x_1\bar x_2)\in2\Z_2\). So the classes listed
	for given \(\eps\) and \(p\bmod4\) are distinct.

	(b) For Jacobians this is Lemma~\ref{lem:tate2}. For products and Weil
	restrictions the cell is given by Lemma~\ref{lem:nonjaccells}, and by (a)
	it determines the class, up to the ambiguity between \(\mathsf N_{+}\) and
	\(\mathsf N_{-}\).
\end{proof}

\begin{lemma}\label{lem:otherprimes}
	\begin{enumerate}[label=\textup{(\alph*)},leftmargin=*]
		\item Let \(\ell\) be a prime not dividing \(2p\). Any two polarised
		\(\ell\)-adic modules of sign \(\eps\) are isomorphic. For \(\eps=-1\)
		they have a hyperbolic basis over \(O_\ell\); for \(\eps=+1\) they are
		free over \(O_\ell\).
		\item Let \(\ell=p\). For \(\eps=+1\) any two polarised \(p\)-adic
		modules are isomorphic, and they are free over \(O_p\). For
		\(\eps=-1\) a polarised \(p\)-adic module \(L\) has an
		\(R_p\)-basis in which \(H=\langle1,u\rangle\) with
		\(u\in\Z_p^{\times}\); two such modules are isomorphic if and only if
		the Hermitian spaces \(L\otimes\Q\) are isometric. If \(-u\) is a
		square, \(L\) has a hyperbolic basis. If \(-u\) is not a square, then
		\(L=\{v\in L\otimes\Q:H(v,v)\in\Z_p\}\).
	\end{enumerate}
\end{lemma}

\begin{proof}
	(a) Here \(R_\ell=O_\ell\) is \'etale over \(\Z_\ell\), and \(L\), which is
	an \(R_\ell\)-module by Definition~\ref{def:ladic}, is free of rank two
	over it. Apply Lemmas~\ref{lem:herm} and~\ref{lem:unram} for
	\(\eps=-1\) and Lemma~\ref{lem:sympl} for \(\eps=+1\).

	(b) \(R_p=O_p\) is a discrete valuation ring with uniformiser \(\pi\) and
	residue field \(\F_p\), and \(L\) is free of rank two. For \(\eps=+1\)
	apply Lemma~\ref{lem:sympl}. Let \(\eps=-1\). If
	\(H(x,x)\in p\Z_p\) for all \(x\in L\), then
	\(\operatorname{Tr}H(x,y)\in p\Z_p\) for all \(x,y\), but \(H(x,y)=1\)
	occurs and \(\operatorname{Tr}(1)=2\). So, as in the proof of
	Proposition~\ref{prop:twoadic}, \(H=\langle u_1,u_2\rangle\) with units
	\(u_i\). The conic \(u_1a^{2}+u_2b^{2}=1\) has a point over \(\F_p\),
	which lifts to \(\Z_p\); the vector \(ae_1+be_2\) has norm \(1\) and is
	part of a basis, and splitting it off gives \(H=\langle1,u\rangle\). The
	norms of the units of \(R_p\) are the squares of \(\Z_p^{\times}\), and a
	unit which is a norm from \(\mathbb K_p^{\times}\) is the norm of a unit;
	so the class of \(u\) modulo squares is the class of \(\det H\) modulo
	norms, an invariant of the Hermitian space. If \(-u=c^{2}\), replace
	\(e_2\) by \(c^{-1}e_2\) to get \(\langle1,-1\rangle\); then
	\(f_1=e_1+e_2\), \(f_2=(e_1-e_2)/2\) is a hyperbolic basis. Suppose that
	\(-u\) is not a square and let \(v\in L\otimes\Q\), \(v\notin L\). Let
	\(k\ge1\) be minimal with \(\pi^{k}v=a'e_1+b'e_2\in L\); one of
	\(a',b'\) is a unit, and
	\(\operatorname N(a')+u\operatorname N(b')\equiv a_0^{2}+ub_0^{2}
	\pmod p\), where \(a_0,b_0\in\Z_p\) are congruent to \(a',b'\) modulo
	\(\pi\). This is a unit, so
	\(H(v,v)=p^{-k}(\operatorname N(a')+u\operatorname N(b'))\notin\Z_p\).
\end{proof}

\begin{table}[ht]
	\centering
	\caption{The normal forms of Table~\ref{tab:twoadic} in \(\Z_2\)-bases, as
		used in the proof of Proposition~\ref{prop:neighbours} (script
		\texttt{anc/proofs/vertexproof\_pm/ctable.py}). Here \(\xi=(1+\pi)/2\);
		\(P\) is the matrix of \(\pi\), its \(j\)-th column giving \(\pi\) applied
		to the \(j\)-th basis vector, and \(G\) is the matrix
		\(\bigl(E(b_i,b_j)\bigr)\) of \(E\). We write
		\(A_c=\bigl(\begin{smallmatrix}0&c\\1&0\end{smallmatrix}\bigr)\),
		\(B_c=\bigl(\begin{smallmatrix}-1&c\\2&1\end{smallmatrix}\bigr)\),
		\(J=\bigl(\begin{smallmatrix}0&-1\\1&0\end{smallmatrix}\bigr)\),
		\(S=\bigl(\begin{smallmatrix}0&1\\1&0\end{smallmatrix}\bigr)\),
		\(Q=\bigl(\begin{smallmatrix}0&1\\1&1\end{smallmatrix}\bigr)\), and
		\(\operatorname{diag}(X,Y)\), resp.\ \([X,Y]\), for the block matrix with
		\(X\), \(Y\) on the diagonal, resp.\ with upper right block \(X\) and
		lower left block \(Y\).}\label{tab:zforms}
	\renewcommand{\arraystretch}{1.2}
	\begin{tabular}{@{}lllll@{}}
		\toprule
		class & \(\eps\) & basis & \(P\) & \(G\)\\
		\midrule
		\(\mathsf{OO}\) & \(-1\) & \(e_1,\xi e_1,e_2,\xi e_2\) &
		\(\operatorname{diag}(B_{-(p+1)/2},B_{-(p+1)/2})\) & \(\operatorname{diag}(J,J)\)\\
		\(\mathsf{RO}\) & \(-1\) & \(e_1,\pi e_1,e_2,\xi e_2\) &
		\(\operatorname{diag}(A_{-p},B_{-(p+1)/2})\) & \(\operatorname{diag}(J,J)\)\\
		\(\mathsf{RR}_{\mathrm n}\), \(\mathsf N_{+}\) & \(-1\) & \(e_1,\pi e_1,e_2,\pi e_2\) &
		\(\operatorname{diag}(A_{-p},A_{-p})\) & \(\operatorname{diag}(J,J)\)\\
		\(\mathsf N_{-}\) & \(-1\) & \(e_1,\pi e_1,e_2,\pi e_2\) &
		\(\operatorname{diag}(A_{-p},A_{-p})\) & \(\operatorname{diag}(J,-J)\)\\
		\(\mathsf{RR}_{\mathrm a}\), \(\mathsf N_{\mathrm a}\) & \(-1\) & \(e_1,\pi e_1,e_2,\pi e_2\) &
		\(\operatorname{diag}(A_{-p},A_{-p})\) & \([J,J]\)\\
		\(\mathsf O^{+}\) & \(+1\) & \(e_1,\xi e_1,e_2,\xi e_2\) &
		\(\operatorname{diag}(B_{(p-1)/2},B_{(p-1)/2})\) & \([Q,-Q^{t}]\)\\
		\(\mathsf R^{+}\) & \(+1\) & \(e_1,\pi e_1,e_2,\pi e_2\) &
		\(\operatorname{diag}(A_{p},A_{p})\) & \([S,-S]\)\\
		\bottomrule
	\end{tabular}
\end{table}

\begin{proof}[Proof of Proposition~\textup{\ref{prop:neighbours}}]
	Let \((A,\theta)\in\mathcal A_\eps\), \((L,\pi,E)=T_2(A,\theta)\), and let
	\(K\subset A[2]=L/2L\) be maximal isotropic and rational, that is,
	\(\pi K=K\). By Lemma~\ref{lem:tate2} the \(2\)-adic module of the
	quotient is
	\[
	L_K=\{x\in\tfrac12L:2x\bmod2L\in K\},\qquad\text{with the form }2E.
	\]
	So the number of \(K\) for which \(A/K\) is of class \(Y\) depends only
	on the class \(X\) of \((A,\theta)\). This is the first assertion of (a).

	(b) Let \(\mathcal G(X,Y)\) be the groupoid of the triples
	\((A,\theta,K)\) with \((A,\theta)\) of class \(X\), \(K\) rational
	maximal isotropic, and \(A/K\) of class \(Y\). The functor \(D\) of
	Lemma~\ref{lem:dual} is an equivalence with \(D\circ D\cong\mathrm{id}\),
	and it maps \(\mathcal G(X,Y)\) to \(\mathcal G(Y,X)\); so these two
	groupoids have the same mass. The functor
	\(\mathcal G(X,Y)\to\mathcal A_\eps[X]\) forgetting \(K\) is faithful,
	and its fibre over \((A,\theta)\) is the set of the \(c(X\to Y)\)
	subgroups \(K\); by Lemma~\ref{lem:fibre},
	\(\Mass\mathcal G(X,Y)=c(X\to Y)\Mass(X)\).

	Table~\ref{tab:neighbours} is the result of a finite
	exact enumeration on the normal forms of Table~\ref{tab:twoadic} (script
	\texttt{anc/proofs/vertexproof\_pm/ctable.py}, Section~\ref{sec:comp});
	this part of the proof is a computation. The \(\Z_2\)-bases used and the
	matrices \(P\) of \(\pi\) and \(G\) of \(E\) in these bases are listed in
	Table~\ref{tab:zforms}; their entries are \(0\), \(\pm1\), \(2\),
	\(\eps p\) and \(\eps(p-\eps)/2\). For each of the \(15\)
	maximal isotropic subgroups \(K\) of \(L/2L\) the enumeration decides
	whether \(\pi K=K\) and, if so, computes the cell of \(L_K\). Let
	\(\tilde L\) be the preimage of \(K\) in \(L\), so that
	\(L_K=\tfrac12\tilde L\). Then
	\(d(L_K)=\dim\ker(\pi-1\mid\tilde L/2\tilde L)\), and
	\(\tilde L/2\tilde L\) is a quotient of \(\tilde L/4L\); and \(L_K\) is
	alternating if and only if \(E(y,\pi y)\equiv0\pmod4\) for all
	\(y\in\tilde L\). For \(y\in\tilde L\) and \(z\in L\),
	\(E(y+2z,\pi(y+2z))-E(y,\pi y)=2(1-\eps)E(y,\pi z)+4E(z,\pi z)\), because
	\(E(z,\pi y)=-\eps E(y,\pi z)\); this is divisible by \(4\), so the residue
	of \(E(y,\pi y)\) modulo \(4\) depends only on \(y\bmod2L\), which lies in
	the finite set \(K\). So the outcome depends only on \(P\) and \(G\) modulo
	\(4\), hence only on \(p\) modulo \(8\), as \(\eps p\) and
	\(\eps(p-\eps)/2\) modulo \(4\) depend only on \(p\bmod8\). The script runs
	all odd \(p=r+64k\) with \(1\le r\le63\), \(0\le k\le5\), \(p\ge5\), and
	confirms that the table depends only on \(p\bmod8\).

	The enumeration records the cell of each \(L_K\), not its class. By
	Proposition~\ref{prop:twoadic}(a) the cell determines the class except
	for \(\mathsf N_{+}\) and \(\mathsf N_{-}\), and these are separated as
	follows. The Hermitian form of \((L_K,\pi,2E)\) is \(2H\) on the same
	space, and \(\det(2H)=\operatorname N(2)\det H\); so
	\(\mathfrak d(L_K)=\mathfrak d(L)\) for the invariant \(\mathfrak d\) defined in the
	proof of Proposition~\ref{prop:twoadic}. As \(\mathfrak d(\mathsf N_{+})=1\) and
	\(\mathfrak d(\mathsf N_{-})=\mathfrak d(\mathsf N_{\mathrm a})=-1\), the quotients
	of a surface of class \(\mathsf N_{+}\) are of class \(\mathsf N_{+}\), and
	the quotients with cell \(1^{2}2^{2}\) of a surface of class
	\(\mathsf N_{-}\) or \(\mathsf N_{\mathrm a}\) are of class
	\(\mathsf N_{-}\).

	Some rows can be checked by hand. For \(\eps=-1\) and
	\(p\equiv3\pmod8\) the determinant of \(H\) has odd valuation for
	\(\mathsf{RO}\) and even valuation for the other three classes, while the
	norms from the unramified extension \(\mathbb K_2\) have even valuation;
	so all quotients of a surface of class \(\mathsf{RO}\) are of class
	\(\mathsf{RO}\). For \(\eps=+1\) and \(p\equiv1\pmod4\): if
	\(L=O_2^{2}\), with \(\Psi=\pi^{-1}\det\), then \(\pi=1\) on \(\bar L\),
	and \(L_K\) is of class \(\mathsf O^{+}\) if and only if it is an
	\(O_2\)-module \(gO_2^{2}\) with \(2\det g\in O_2^{\times}\)
	(Lemma~\ref{lem:sympl}), that is, if and only if \(K\) is a free module
	of rank one over \(O_2/2O_2\); there are \(3\cdot3\) such \(K\) if
	\(O_2/2O_2=\F_2\times\F_2\) and \(5\) if \(O_2/2O_2=\F_4\). If
	\(L=R_2^{2}\), an \(O_2\)-module \(L_K\) contains \(O_2L\), which has
	index four over \(L\); so \(L_K=O_2L\) is the only quotient of class
	\(\mathsf O^{+}\).
\end{proof}

\subsection*{Genera and masses}

For a group \(\mathcal H\) with subgroups \(\Gamma\) and \(U\) such that
\(\Gamma\backslash\mathcal H/U\) is finite and all the groups
\(\Gamma\cap gUg^{-1}\) are finite, we put
\[
m(\Gamma,\mathcal H,U)=\sum_{\Gamma gU\in\Gamma\backslash\mathcal H/U}
\frac{1}{\#(\Gamma\cap gUg^{-1})} .
\]

\begin{lemma}\label{lem:hasse}
	Let \(s\) and \(s'\) be positive definite Hermitian forms on
	\(\mathbb K^{2}\). If \(\det s'/\det s\) is a norm from
	\(\mathbb K_\ell\) for every prime \(\ell\ne p\), then \(s\) and \(s'\)
	are isometric.
\end{lemma}

\begin{proof}
	Let \(d=\det s'/\det s\), a positive rational number. The Hilbert symbol
	\((d,-p)_v\) is \(1\) at the real place and at all primes \(v\ne p\),
	hence also at \(p\) by the product formula
	\cite[Prop.~14.2.1]{Voight2021}. So the conic
	\(x^{2}+py^{2}=dz^{2}\) has points over \(\mathbb R\) and over all
	\(\Q_\ell\), hence over \(\Q\) by the Hasse--Minkowski theorem
	\cite[Thm.~14.3.3]{Voight2021}; as \(x^{2}+py^{2}\) is definite,
	\(z\ne0\) and \(d\in\operatorname N(\mathbb K^{\times})\). Diagonalise
	\(s=\langle a,b\rangle\) and \(s'=\langle a',b'\rangle\) with positive
	rational numbers. The quadratic form
	\(a(x_1^{2}+px_2^{2})+b(y_1^{2}+py_2^{2})-a'w^{2}\) in five variables is
	indefinite, hence it has a non-trivial rational zero by Meyer's theorem
	\cite[Ch.~IV, \S3.2, Cor.~2]{Serre1973}, and \(w\ne0\) at this zero. So
	\(s\) represents \(a'\), and \(s\cong\langle a',c\rangle\) with \(c>0\).
	Comparing determinants, \(c/b'\in\operatorname N(\mathbb K^{\times})\),
	so \(\langle c\rangle\cong\langle b'\rangle\).
\end{proof}

\begin{lemma}\label{lem:genus}
	Let \(X\) be a \(2\)-adic class of sign \(\eps\) and
	\((A_0,\theta_0)\in\mathcal A_\eps[X]\). Let \(G\) be the algebraic group
	over \(\Q\) of the elements \(g\) of \(\End^{0}(A_0)^{\times}\) with
	\(g^{\dagger}g=1\). For every prime \(\ell\), \(G(\Q_\ell)\) is the group
	of automorphisms of \(T_\ell A_0\otimes\Q\) that commute with \(\pi\) and
	preserve \(E_{\theta_0}\); let \(U_\ell\) be the stabiliser of
	\(T_\ell A_0\) in it, and \(U=\prod_\ell U_\ell\subset G(\mathbb A_f)\).
	Then
	\[
	\Mass(X)=m\bigl(G(\Q),G(\mathbb A_f),U\bigr).
	\]
\end{lemma}

\begin{proof}
	The description of \(G(\Q_\ell)\) is Tate's theorem
	\cite[Main Theorem]{Tate1966}, with
	the Dieudonn\'e module for \(\ell=p\); see
	\cite[Thm.~5.3(a) and~5.7(a)]{XueYu2021}. By
	\cite[Thm.~5.8]{XueYu2021}, which rests on Lemma~5.2 and on these two
	theorems of that paper (the case without polarisation is Proposition~5.4
	there), \(m(G(\Q),G(\mathbb A_f),U)\) is the mass of the \emph{genus} of
	\((A_0,\theta_0)\): the groupoid of the abelian varieties \(A\) over
	\(\Fp\) with a \(\Q\)-polarisation \(\theta\) such that there is a
	quasi-isogeny \(\varphi\colon A\to A_0\) with
	\(\varphi^{*}\theta_0=\theta\) and such that, for every prime \(\ell\),
	the polarised \(\ell\)-divisible groups of \((A,\theta)\) and
	\((A_0,\theta_0)\) are isomorphic. The last condition means that
	\(T_\ell(A,\theta)\cong T_\ell(A_0,\theta_0)\) for all \(\ell\),
	including \(\ell=p\). The members of the genus are principally polarised
	\cite[Lemma~2.5]{XueYu2020}, they satisfy \(\pi^{2}=\eps p\), and they
	are of class \(X\). So we have to show that every
	\((A,\theta)\in\mathcal A_\eps[X]\) is in the genus.

	The modules \(T_\ell(A,\theta)\) and \(T_\ell(A_0,\theta_0)\) are
	isomorphic for \(\ell\nmid2p\) by Lemma~\ref{lem:otherprimes}(a), and for
	\(\ell=2\) by hypothesis. For \(\eps=+1\) there is a quasi-isogeny
	\(\varphi\) with \(\varphi^{*}\theta_0=\theta\) by
	\cite[Lemma~3.1.3]{XueYu2020}, and the modules at \(p\) are isomorphic by
	Lemma~\ref{lem:otherprimes}(b).

	Let \(\eps=-1\). By Lemma~\ref{lem:Aeps}(a) and Tate's theorem
	\cite[Thm.~1(c)]{Tate1966}, \(A\) and
	\(A_0\) are isogenous to \(E^{2}\) for a supersingular elliptic curve
	\(E/\Fp\), and \(\End^{0}_{\Fp}(E)=\Q(\pi_E)\cong\mathbb K\). So
	\(\End^{0}(A_0)\cong M_2(\mathbb K)\), with \(\pi_{A_0}\) corresponding
	to \(\sqrt{-p}\). The Rosati involution is positive and induces the
	conjugation on the centre \(\mathbb K\); such an involution of
	\(M_2(\mathbb K)\) is the adjoint involution of a Hermitian form \(s_0\)
	on \(\mathbb K^{2}\), unique up to a rational factor, and \(s_0\) can be
	taken positive definite. Indeed, with \(b^{*}={}^{t}\bar b\), the map
	\(b\mapsto(b^{\dagger})^{*}\) is a \(\mathbb K\)-linear automorphism of
	\(M_2(\mathbb K)\), hence inner by the Skolem--Noether theorem
	\cite[Cor.~7.7.4]{Voight2021}: \(b^{\dagger}=s_0^{-1}b^{*}s_0\). As
	\(\dagger\) is an involution, \(s_0^{*}=cs_0\) with \(c\in\mathbb K\),
	\(c\bar c=1\), and replacing \(s_0\) by \(as_0\) with \(a/\bar a=c\)
	(Hilbert's Theorem~90) makes \(s_0\) Hermitian; it is unique up to
	\(\Q^{\times}\). In a basis in which \(s_0=\langle\alpha,\beta\rangle\),
	positivity applied to \(b=e_{12}\) gives \(\operatorname{tr}(bb^{\dagger})
	=\alpha/\beta>0\), so \(s_0\) is definite. Choose a quasi-isogeny
	\(\varphi_1\colon A\to A_0\). The \(\Q\)-polarisation
	\((\varphi_1^{-1})^{*}\theta\) of \(A_0\) is \(\theta_0a\) with
	\(a\in\End^{0}(A_0)\) symmetric and totally positive
	\cite[\S21, Application~III, p.~208]{Mumford2008}; in other words \(s_0a\) is a
	positive definite Hermitian form. For \(g\in\End^{0}(A_0)^{\times}\) the
	quasi-isogeny \(\varphi=g\varphi_1\) satisfies
	\(\varphi^{*}\theta_0=\theta\) if and only if \(g^{\dagger}g=a\), that
	is, if and only if \(g\) is an isometry from \(s_0a\) to \(s_0\). Let
	\(\ell\ne p\). The map \(\varphi_1\) is an isometry from
	\(T_\ell A\otimes\Q\) with \(E_\theta\) to \(T_\ell A_0\otimes\Q\) with
	the form \(E_{\theta_0}(x,ay)\); as the lattices \(T_\ell A\) and
	\(T_\ell A_0\) are isometric, there is an automorphism \(g_\ell\) of
	\(T_\ell A_0\otimes\Q\) commuting with \(\pi\) such that
	\(g_\ell^{\dagger}g_\ell=a\), and then
	\(\det a=\operatorname N(\det g_\ell)\) is a norm from \(\mathbb K_\ell\).
	By Lemma~\ref{lem:hasse}, \(s_0a\) and \(s_0\) are isometric, and
	\(\varphi\) exists. Finally \(\varphi\) maps \(T_pA\) to a unimodular
	lattice in the Hermitian space of \(T_pA_0\), so
	\(T_p(A,\theta)\cong T_p(A_0,\theta_0)\) by
	Lemma~\ref{lem:otherprimes}(b).
\end{proof}

\begin{lemma}\label{lem:normal}
	Let \(\mathcal H\) be a group, \(\det\colon\mathcal H\to\mathcal T\) a
	surjective homomorphism onto an abelian group, \(\mathcal S\) its
	kernel, and \(\Gamma\), \(U\) subgroups of \(\mathcal H\) such that
	\(\Gamma\backslash\mathcal H/U\) is finite and the groups
	\(\Gamma\cap gUg^{-1}\) are finite. Put
	\(\Gamma_{\mathcal T}=\det\Gamma\), \(W=\det U\), and assume that
	\(\mathcal T/\Gamma_{\mathcal T}W\) and \(\Gamma_{\mathcal T}\cap W\) are
	finite. For \(\tau\in\mathcal T\) choose \(g_\tau\in\mathcal H\) with
	\(\det g_\tau=\tau\). Then
	\[
	m(\Gamma,\mathcal H,U)=\frac{1}{\#(\Gamma_{\mathcal T}\cap W)}
	\sum_{\tau\in\mathcal T/\Gamma_{\mathcal T}W}
	m\bigl(\Gamma\cap\mathcal S,\ \mathcal S,\
	g_\tau(U\cap\mathcal S)g_\tau^{-1}\bigr).
	\]
\end{lemma}

\begin{proof}
	The class of \(\det g\) in \(\mathcal T/\Gamma_{\mathcal T}W\) depends
	only on \(\Gamma gU\). Fix \(\tau\). If \(\det g\) lies in the class of
	\(\tau\), then \(\Gamma gU=\Gamma sg_\tau U\) with \(s\in\mathcal S\),
	and \(\Gamma xU\mapsto\Gamma xg_\tau^{-1}U^{\tau}\), with
	\(U^{\tau}=g_\tau Ug_\tau^{-1}\), is a bijection from the double cosets
	over \(\tau\) to the double cosets \(\Gamma sU^{\tau}\),
	\(s\in\mathcal S\), which preserves the groups
	\(\Gamma\cap xUx^{-1}\). As \(\det U^{\tau}=W\) and
	\(U^{\tau}\cap\mathcal S=g_\tau(U\cap\mathcal S)g_\tau^{-1}\), it
	suffices to prove that
	\begin{equation}\label{eq:detbal}
		m\bigl(\Gamma\cap\mathcal S,\mathcal S,U\cap\mathcal S\bigr)
		=\#(\Gamma_{\mathcal T}\cap W)\sum_{\Gamma sU,\ s\in\mathcal S}
		\frac{1}{\#(\Gamma\cap sUs^{-1})} .
	\end{equation}
	Consider the map
	\((\Gamma\cap\mathcal S)\,s\,(U\cap\mathcal S)\mapsto\Gamma sU\). Let
	\(s\in\mathcal S\) and \(\Gamma_s=\Gamma\cap sUs^{-1}\). If
	\(s'=\gamma su\in\mathcal S\) with \(\gamma\in\Gamma\), \(u\in U\), then
	\(\det\gamma=(\det u)^{-1}\in\Gamma_{\mathcal T}\cap W\), and \(\gamma\)
	is determined by \(s'\) up to right multiplication by \(\Gamma_s\). This
	defines a map from the fibre over \(\Gamma sU\) to
	\((\Gamma_{\mathcal T}\cap W)/\det\Gamma_s\). It is injective: if
	\(s''=\gamma'su'\) and \(\det\gamma'=\det\gamma\det\eta\) with
	\(\eta\in\Gamma_s\), write \(s''=(\gamma'\eta^{-1})s(s^{-1}\eta su')\);
	so we may assume \(\det\gamma'=\det\gamma\), and then
	\(\gamma'\gamma^{-1}\in\Gamma\cap\mathcal S\),
	\(u^{-1}u'\in U\cap\mathcal S\) and
	\(s''=(\gamma'\gamma^{-1})s'(u^{-1}u')\). It is surjective because
	\(\mathcal T\) is abelian. The group
	\((\Gamma\cap\mathcal S)\cap s'(U\cap\mathcal S)s'^{-1}\) is conjugate
	by \(\gamma\) to \(\Gamma_s\cap\mathcal S\), of order
	\(\#\Gamma_s/\#\det\Gamma_s\). Hence the fibre over \(\Gamma sU\)
	contributes \(\#(\Gamma_{\mathcal T}\cap W)/\#\Gamma_s\) to the left-hand
	side of \eqref{eq:detbal}.
\end{proof}

\begin{lemma}\label{lem:conj}
	Let \(\mathcal S\) be a locally compact, totally disconnected, unimodular
	group and \(\Gamma_{\mathcal S}\) a subgroup such that
	\(m(\Gamma_{\mathcal S},\mathcal S,U')\) is defined for one, hence for
	every, open compact subgroup \(U'\).
	\begin{enumerate}[label=\textup{(\alph*)},leftmargin=*]
		\item For open compact subgroups \(U''\subset U'\),
		\(m(\Gamma_{\mathcal S},\mathcal S,U'')
		=[U':U'']\,m(\Gamma_{\mathcal S},\mathcal S,U')\). Hence
		\(m(\Gamma_{\mathcal S},\mathcal S,U')\operatorname{vol}(U')\) does
		not depend on \(U'\), for a fixed Haar measure.
		\item Suppose that \(\mathcal S\) is a closed normal subgroup of a
		locally compact group \(\mathcal H\) and that the square of every
		element of \(\mathcal H\) lies in the product of \(\mathcal S\) and
		the centre of \(\mathcal H\). Then
		\(\operatorname{vol}(gU'g^{-1})=\operatorname{vol}(U')\) and
		\(m(\Gamma_{\mathcal S},\mathcal S,gU'g^{-1})
		=m(\Gamma_{\mathcal S},\mathcal S,U')\) for all \(g\in\mathcal H\).
	\end{enumerate}
\end{lemma}

\begin{proof}
	(a) The fibre of
	\(\Gamma_{\mathcal S}\backslash\mathcal S/U''\to
	\Gamma_{\mathcal S}\backslash\mathcal S/U'\) over
	\(\Gamma_{\mathcal S}sU'\) is the set of orbits of the finite group
	\(\Gamma_{\mathcal S}\cap sU's^{-1}\) on \(sU'/U''\), the stabiliser of
	\(suU''\) being
	\(\Gamma_{\mathcal S}\cap suU''u^{-1}s^{-1}\). The sum of the inverses of
	the orders of the stabilisers over the orbits is
	\([U':U'']/\#(\Gamma_{\mathcal S}\cap sU's^{-1})\).
	(b) Conjugation by \(g\) multiplies the Haar measure of \(\mathcal S\) by
	a factor \(\Delta(g)\), and \(\Delta\) is a homomorphism from
	\(\mathcal H\) to the positive real numbers, trivial on \(\mathcal S\),
	which is unimodular, and on the centre. By hypothesis
	\(\Delta(g)^{2}=1\), so \(\Delta=1\). The second assertion follows from
	(a), applied to \(U'\cap gU'g^{-1}\) in \(U'\) and in \(gU'g^{-1}\).
\end{proof}

In Lemmas~\ref{lem:quat} and~\ref{lem:hyp}, \(\eps=-1\). We fix \((A_0,\theta_0)\in\mathcal A_{-1}\)
and use the notation of Lemma~\ref{lem:genus} and of its proof:
\(\End^{0}(A_0)=M_2(\mathbb K)\), \(\dagger\) is the adjoint involution of
the positive definite Hermitian form \(s_0\), and \(G\) is the unitary group
of \((\mathbb K^{2},s_0)\). Let \(\mathbb K^{1}\) be the torus of the
elements of norm one of \(\mathbb K\), so that
\(\det\colon G\to\mathbb K^{1}\), and let \(G^{1}\) be the kernel of
\(\det\). For a prime \(\ell\) put \(V_\ell=T_\ell A_0\otimes\Q\) (for
\(\ell=p\) this is the Dieudonn\'e module tensored with \(\Q\)), with the
Hermitian form \(H\) of Lemma~\ref{lem:herm}. By Tate's theorem
\(\End^{0}(A_0)\otimes\Q_\ell\) is the algebra of the
\(\mathbb K_\ell\)-linear endomorphisms of \(V_\ell\), and \(\dagger\) is
the adjoint involution of \(H\); so \(G(\Q_\ell)\) is the unitary group and
\(G^{1}(\Q_\ell)\) the special unitary group of \((V_\ell,H)\). For
\(b\in M_2(\mathbb K)\) let \(\bar b=\operatorname{tr}(b)-b\) be the
adjugate matrix.

\begin{lemma}\label{lem:quat}
	\begin{enumerate}[label=\textup{(\alph*)},leftmargin=*]
		\item \(\mathfrak B=\{b\in M_2(\mathbb K):b^{\dagger}=\bar b\}\) is a
		definite quaternion algebra over \(\Q\) with
		\(\mathfrak B\otimes\mathbb K=M_2(\mathbb K)\); its reduced norm is
		the determinant, and \(G^{1}\) is the group \(\mathfrak B^{1}\) of
		its elements of reduced norm one.
		\item Let \(\ell\) be a prime. If \(V_\ell\) has a hyperbolic basis
		\(f_1,f_2\), then \(\mathfrak B_\ell=\mathfrak B\otimes\Q_\ell\) is
		the algebra of the \(\Q_\ell\)-linear endomorphisms of
		\(V_{0}=\Q_\ell f_1\oplus\Q_\ell\pi f_2\), extended
		\(\mathbb K_\ell\)-linearly to \(V_\ell=V_0\otimes\mathbb K_\ell\).
		\item Let \(\ell\) be a prime. If
		\(M_0=\{v\in V_\ell:H(v,v)\in\Z_\ell\}\) is a lattice, then
		\(\mathfrak B_\ell\) is a division algebra, and
		\(\mathfrak B_\ell^{1}\) is the group of the elements of reduced norm
		one of its maximal order.
		\item The number of primes at which \(\mathfrak B\) is ramified is
		odd.
	\end{enumerate}
\end{lemma}

\begin{proof}
	(a) Both \(b\mapsto b^{\dagger}\) and \(b\mapsto\bar b\) are
	anti-automorphisms of the ring \(M_2(\mathbb K)\); the first induces the
	conjugation on \(\mathbb K\) and satisfies
	\(\operatorname{tr}(b^{\dagger})=\overline{\operatorname{tr}(b)}\), the
	second is \(\mathbb K\)-linear. They commute, so
	\(b\mapsto(\bar b)^{\dagger}\) is a ring automorphism of order two which
	is semilinear for the conjugation, and its fixed ring \(\mathfrak B\)
	satisfies \(\mathfrak B\otimes\mathbb K=M_2(\mathbb K)\) by Galois
	descent. For \(b\in\mathfrak B\),
	\(b^{\dagger}b=\bar bb=\det b\), and \(\det b\) is fixed by the
	conjugation. If \(g^{\dagger}g=1\), then \(\det g=1\) if and only if
	\(\bar g=g^{-1}=g^{\dagger}\), that is, \(g\in\mathfrak B\); this holds
	over every \(\Q\)-algebra. As \(s_0\) is definite, \(G(\mathbb R)\) is
	compact, and so is \(\mathfrak B^{1}(\mathbb R)\).

	(b) For \(v=af_1+c\pi f_2\) and \(w=bf_1+d\pi f_2\) with
	\(a,b,c,d\in\Q_\ell\) we have \(H(v,w)=-\pi(ad-bc)\). Hence, for a
	\(\Q_\ell\)-linear endomorphism \(g\) of \(V_0\),
	\(H(gv,w)=H(v,\bar gw)\) on \(V_0\) and, by sesquilinearity, on
	\(V_\ell\). So \(g^{\dagger}=\bar g\), and the four-dimensional algebra
	of these endomorphisms is contained in \(\mathfrak B_\ell\).

	(c) The unitary group of \(V_\ell\) stabilises \(M_0\), so it is
	compact; hence \(\mathfrak B_\ell^{1}\) is compact and
	\(\mathfrak B_\ell\) is not a matrix algebra. In a division quaternion
	algebra over \(\Q_\ell\) the maximal order consists of the elements of
	integral reduced norm.

	(d) \(\mathfrak B\) is ramified at the real place, and the number of
	places at which a quaternion algebra over \(\Q\) is ramified is even, by
	the product formula for the Hilbert symbol
	\cite[Prop.~14.2.1 and Cor.~14.2.3]{Voight2021}.
\end{proof}

\begin{lemma}\label{lem:hyp}
	Let \(\ell\) be a prime, \(f_1,f_2\) a hyperbolic basis of \(V_\ell\) and
	\(\Lambda=\mathfrak of_1\oplus\mathfrak of_2\), where \(\mathfrak o\) is
	\(R_\ell\) or \(O_\ell\). In the basis \((f_1,\pi f_2)\) of \(V_0\) the
	stabiliser of \(\Lambda\) in
	\(\mathfrak B_\ell^{1}=\operatorname{SL}_2(\Q_\ell)\) is
	\(\operatorname{SL}_2(\Z_\ell)\) if \(\ell\ne p\); if \(\ell=p\) it is the
	group of the matrices
	\(\bigl(\begin{smallmatrix}a&b\\c&d\end{smallmatrix}\bigr)\in
	\operatorname{SL}_2(\Z_p)\) with \(b\in p\Z_p\), of index \(p+1\) in
	\(\operatorname{SL}_2(\Z_p)\).
\end{lemma}

\begin{proof}
	Let \(g=\bigl(\begin{smallmatrix}a&b\\c&d\end{smallmatrix}\bigr)\), so
	that \(gf_1=af_1+c\pi f_2\) and \(gf_2=\pi^{-1}bf_1+df_2\). As \(g\) has
	determinant one, it stabilises \(\Lambda\) if and only if
	\(g\Lambda\subset\Lambda\), that is, \(a,d\in\mathfrak o\cap\Q_\ell\),
	\(c\in\pi^{-1}\mathfrak o\cap\Q_\ell\) and
	\(b\in\pi\mathfrak o\cap\Q_\ell\). If \(\ell\ne p\), \(\pi\) is a unit of
	\(\mathfrak o\). If \(\ell=p\), \(\pi\) is a uniformiser of
	\(\mathfrak o=R_p\), so \(\pi^{-1}R_p\cap\Q_p=\Z_p\) and
	\(\pi R_p\cap\Q_p=p\Z_p\).
\end{proof}

Eichler's mass formula is the following statement, which we take from
\cite[Thm.~26.5.4 and (26.1.8)]{Voight2021}; see also
\cite[(5.1) and (5.5)]{XueYangYu2019} and
\cite[(5.1), (5.2) and Remark~5.2]{KirschmerVoight2010}. The original source
is \cite[Cor.~V.2.3]{Vigneras1980}. Let \(\mathbb M\) be a totally real
number field of degree \(n\), with ring of integers \(\mathcal O\) and class
number \(h_{\mathbb M}\), let \(\mathfrak D\) be a totally definite
quaternion algebra over \(\mathbb M\) and \(\mathfrak O\) a maximal order of
\(\mathfrak D\). Let \(I_1,\dots,I_r\) represent the classes of right ideals
of \(\mathfrak O\), with left orders \(\mathfrak O_i\). Then
\begin{equation}\label{eq:eichler}
	\sum_{i=1}^{r}\frac{1}{[\mathfrak O_i^{\times}:\mathcal O^{\times}]}
	=2^{1-n}\,\lvert\zeta_{\mathbb M}(-1)\rvert\,h_{\mathbb M}
	\prod_{\mathfrak p}(\operatorname N\mathfrak p-1),
\end{equation}
the product being over the primes at which \(\mathfrak D\) is ramified. We
write \(\widehat{\mathfrak D}=\mathfrak D\otimes\mathbb A_f\),
\(\widehat{\mathfrak O}=\mathfrak O\otimes\hat\Z\), and an exponent \(1\)
denotes the elements of reduced norm one.

\begin{lemma}\label{lem:normone}
	With this notation,
	\[
	m\bigl(\mathfrak D^{1},\widehat{\mathfrak D}^{1},
	\widehat{\mathfrak O}^{1}\bigr)
	=2^{-n}\,\lvert\zeta_{\mathbb M}(-1)\rvert
	\prod_{\mathfrak p}(\operatorname N\mathfrak p-1).
	\]
	It is \((q-1)/24\) if \(\mathbb M=\Q\) and \(\mathfrak D\) is ramified
	exactly at a prime \(q\) and at the real place, and it is \(z_p/4\) if
	\(\mathbb M=\mathbb K_{+}\) and \(\mathfrak D\) is ramified exactly at the
	two real places.
\end{lemma}

\begin{proof}
	Let \(Z=\mathcal O^{\times}\), let \(E_{+}\) be the group of totally
	positive units and \(h^{+}\) the narrow class number of \(\mathbb M\).
	The reduced norm maps \(\mathfrak D^{\times}\) onto the group
	\(\mathbb M^{\times}_{+}\) of totally positive elements, by the
	Hasse--Schilling norm theorem \cite[Main Theorem~14.7.4]{Voight2021}, as
	\(\mathfrak D\) is ramified at all real places (for \(\mathbb M=\Q\) this
	also follows from Meyer's theorem, used in the proof of
	Lemma~\ref{lem:hasse}, applied to the quinary form
	\(\operatorname{nrd}(x)-cw^{2}\), \(c>0\)); it maps
	\(\widehat{\mathfrak D}^{\times}\) onto
	\((\mathbb M\otimes\mathbb A_f)^{\times}\) and
	\(\widehat{\mathfrak O}^{\times}\) onto
	\(\widehat{\mathcal O}^{\times}=(\mathcal O\otimes\hat\Z)^{\times}\).
	The classes of right ideals of \(\mathfrak O\) correspond to
	\(\mathfrak D^{\times}\backslash\widehat{\mathfrak D}^{\times}/
	\widehat{\mathfrak O}^{\times}\), the class of \(x\) having left order
	\(\mathfrak O_x=\mathfrak D\cap x\widehat{\mathfrak O}x^{-1}\). Put
	\(\Gamma_x=\mathfrak O_x^{\times}\) and
	\(\Gamma_x^{1}=\Gamma_x\cap\mathfrak D^{1}\); the latter is finite and
	\(Z\) has finite index in \(\Gamma_x\).

	Fix \(x\in\widehat{\mathfrak D}^{\times}\), put
	\(U_x=x\widehat{\mathfrak O}^{\times}x^{-1}\), and let \(C_x\) be the set
	of the double cosets
	\(\mathfrak D^{\times}x'\widehat{\mathfrak O}^{\times}\) such that
	\(\operatorname{nrd}(x')\in\operatorname{nrd}(x)\,\mathbb M^{\times}_{+}
	\widehat{\mathcal O}^{\times}\). The map
	\[
	\beta_x\colon\mathfrak D^{1}\backslash\widehat{\mathfrak D}^{1}/U_x^{1}
	\to C_x,\qquad \mathfrak D^{1}yU_x^{1}\mapsto
	\mathfrak D^{\times}yx\widehat{\mathfrak O}^{\times},
	\]
	is surjective, by the surjectivity of the three reduced norm maps. We
	determine its fibres. An element of \(C_x\) is the class of \(x'=yx\)
	with \(y\in\widehat{\mathfrak D}^{1}\); as \(U_{x'}=yU_xy^{-1}\), right
	multiplication by \(y\) identifies the source of \(\beta_{x'}\) with
	that of \(\beta_x\), compatibly with the two maps and with the
	stabilisers. So it suffices to describe the
	fibre of \(\beta_x\) over the class of \(x\). It consists of the classes of the
	\(y=\gamma w\) with \(\gamma\in\mathfrak D^{\times}\), \(w\in U_x\) and
	\(\operatorname{nrd}(\gamma)\operatorname{nrd}(w)=1\); then
	\(\operatorname{nrd}(\gamma)\in\mathbb M^{\times}_{+}\cap
	\widehat{\mathcal O}^{\times}=E_{+}\), and \(\gamma\) is determined by
	\(y\) up to right multiplication by \(\Gamma_x\). As in the proof of
	Lemma~\ref{lem:normal}, but with the groups \(\Gamma_x\), which are
	infinite when \(\mathbb M\ne\Q\), one checks directly that
	\(y\mapsto\operatorname{nrd}(\gamma)\) induces a bijection from the fibre
	to \(E_{+}/\operatorname{nrd}(\Gamma_x)\): it is well defined on double
	cosets; if \(y'=\gamma'w'\) and
	\(\operatorname{nrd}(\gamma')=\operatorname{nrd}(\gamma\eta)\) with
	\(\eta\in\Gamma_x\), then
	\(y'=(\gamma'\eta^{-1}\gamma^{-1})\,y\,(w^{-1}\eta w')\) with the two
	outer factors of reduced norm one; and every element of \(E_{+}\) is
	the reduced norm of some \(\gamma\) and the inverse of the reduced norm
	of some \(w\). The group
	\(\mathfrak D^{1}\cap yU_x^{1}y^{-1}\) is conjugate by \(\gamma\) to
	\(\Gamma_x^{1}\). Now
	\(Z^{2}\subset\operatorname{nrd}(\Gamma_x)\subset E_{+}\), the kernel of
	\(\Gamma_x\to\operatorname{nrd}(\Gamma_x)/Z^{2}\) is \(\Gamma_x^{1}Z\),
	and \(\Gamma_x^{1}\cap Z=\{\pm1\}\); so
	\[
	\frac{[E_{+}:\operatorname{nrd}\Gamma_x]}{\#\Gamma_x^{1}}
	=\frac{[E_{+}:Z^{2}]}{[\Gamma_x:\Gamma_x^{1}Z]\,\#\Gamma_x^{1}}
	=\frac{[E_{+}:Z^{2}]}{2\,[\Gamma_x:Z]} .
	\]
	Summing over \(C_x\) we obtain
	\[
	m\bigl(\mathfrak D^{1},\widehat{\mathfrak D}^{1},U_x^{1}\bigr)
	=\frac{[E_{+}:Z^{2}]}{2}\sum_{c\in C_x}\frac{1}{[\Gamma_c:Z]},
	\]
	where \(\Gamma_c=\Gamma_{x'}\) for a representative \(x'\) of \(c\). The
	group \(U_x^{1}\) is conjugate to \(\widehat{\mathfrak O}^{1}\) by an
	element of \(\widehat{\mathfrak D}^{\times}\), and the square of every
	element of \(\widehat{\mathfrak D}^{\times}\) is the product of an
	element of reduced norm one and of a central element; so the left-hand
	side is \(m(\mathfrak D^{1},\widehat{\mathfrak D}^{1},
	\widehat{\mathfrak O}^{1})\) by Lemma~\ref{lem:conj}(b). The sets
	\(C_x\) are the fibres of the reduced norm map to the narrow class
	group
	\((\mathbb M\otimes\mathbb A_f)^{\times}/\mathbb M^{\times}_{+}
	\widehat{\mathcal O}^{\times}\), which is surjective. Summing over the
	\(h^{+}\) narrow classes, and using
	\(h^{+}=h_{\mathbb M}[E_{+}:Z^{2}]\), we find that the left-hand side of
	\eqref{eq:eichler} equals
	\(2h_{\mathbb M}\,m(\mathfrak D^{1},\widehat{\mathfrak D}^{1},
	\widehat{\mathfrak O}^{1})\). The two special cases follow from
	\(\zeta_{\Q}(-1)=-1/12\) and from the definition of \(z_p\).
\end{proof}

In Lemma~\ref{lem:conj}(b), applied to
\(\mathcal S=\widehat{\mathfrak D}^{1}\) or to
\(\mathcal S=G^{1}(\mathbb A_f)\), the group \(\mathcal S\) is unimodular
because it is a restricted product of groups each of which is compact or a
product of special linear groups of degree two over local fields.

\begin{lemma}\label{lem:torus}
	Let \(\eps=-1\), let \(\Gamma=\mathbb K^{1}(\Q)\) be the group of the
	elements of norm one of \(\mathbb K\),
	\(W_{\max}=\prod_\ell O_\ell^{1}\subset\mathbb K^{1}(\mathbb A_f)\) the
	group of the units of norm one, and \(W\subset W_{\max}\) a subgroup of
	finite index. Then
	\[
	\frac{\#\bigl(\mathbb K^{1}(\mathbb A_f)/\Gamma W\bigr)}{\#(\Gamma\cap W)}
	=\frac{h_{\mathbb K}}{c_{\mathbb K}}\cdot\frac{[W_{\max}:W]}{2},
	\]
	where \(c_{\mathbb K}=1\) if \(p\equiv3\pmod4\) and \(c_{\mathbb K}=2\)
	if \(p\equiv1\pmod4\).
\end{lemma}

\begin{proof}
	If \(\ell\) is inert or ramified in \(\mathbb K\), every element of norm
	one of \(\mathbb K_\ell\) is a unit. If \(\ell\) splits, the elements of
	norm one are the \((a,a^{-1})\), \(a\in\Q_\ell^{\times}\), and the
	valuation of \(a\) identifies their quotient by the units with \(\Z\).
	Hence \(\mathbb K^{1}(\mathbb A_f)/W_{\max}\) is the group \(I^{-}\) of
	the fractional ideals \(\mathfrak a\) of \(\mathbb K\) with
	\(\mathfrak a\bar{\mathfrak a}=(1)\), and \(\Gamma\) maps onto the group
	\(P^{-}\) of the principal ideals \((c)\) with
	\(\operatorname N(c)=1\). The map
	\(\mathfrak b\mapsto\mathfrak b/\bar{\mathfrak b}\) from the class group
	to \(I^{-}/P^{-}\) is surjective. The class of \(\mathfrak b\) is in its
	kernel if and only if \(\mathfrak b/\bar{\mathfrak b}=(c)\) with
	\(\operatorname N(c)=1\); writing \(c=\bar e/e\) (Hilbert's Theorem~90)
	this means that the ideal \(\mathfrak be\) is equal to its conjugate,
	that is, is the product of a rational ideal and of ramified primes. So
	the kernel is generated by the classes of the ramified primes. If
	\(p\equiv3\pmod4\) the only ramified prime is \((\sqrt{-p})\). If
	\(p\equiv1\pmod4\) there is also the prime above \(2\), which is not
	principal because \(a^{2}+pb^{2}=2\) has no solution, and whose square
	is \((2)\). Hence \(\#(\mathbb K^{1}(\mathbb A_f)/\Gamma W_{\max})
	=h_{\mathbb K}/c_{\mathbb K}\). The group \(\Gamma\cap W_{\max}\)
	consists of the roots of unity of \(\mathbb K\), which are \(\pm1\).
	Finally
	\([\Gamma W_{\max}:\Gamma W]=[W_{\max}:W]/[\{\pm1\}:\{\pm1\}\cap W]\) and
	\(\Gamma\cap W=\{\pm1\}\cap W\).
\end{proof}

\subsection*{Proofs of the theorems}

\begin{lemma}\label{lem:genusmass}
	Let \(\eps=-1\), let \(X\) be a \(2\)-adic class,
	\((A_0,\theta_0)\in\mathcal A_{-1}[X]\), and let \(U=\prod_\ell U_\ell\)
	be as in Lemma~\textup{\ref{lem:genus}}. Put
	\(U_\ell^{1}=U_\ell\cap G^{1}(\Q_\ell)\) and \(W_2=\det U_2\). The
	algebra \(\mathfrak B\) of Lemma~\textup{\ref{lem:quat}} is ramified at
	the real place and at exactly one prime \(q\), and \(q\in\{2,p\}\). Let
	\(i_2\) be the quotient of the volume of the group of the elements of
	reduced norm one of a maximal order of \(\mathfrak B_2\) by the volume
	of \(U_2^{1}\). Then
	\begin{equation}\label{eq:genusmass}
		\Mass(X)=\frac{q-1}{24}\cdot i_2\,i_p\cdot
		\frac{h_{\mathbb K}}{c_{\mathbb K}}\cdot\frac{[O_2^{1}:W_2]}{2},
		\qquad i_p=\begin{cases}1,&q=p,\\p+1,&q=2.\end{cases}
	\end{equation}
\end{lemma}

\begin{proof}
	By Lemma~\ref{lem:genus},
	\(\Mass(X)=m(G(\Q),G(\mathbb A_f),U)\). We apply
	Lemma~\ref{lem:normal} to \(\mathcal H=G(\mathbb A_f)\),
	\(\Gamma=G(\Q)\), \(U\), and the determinant with values in
	\(\mathcal T=\mathbb K^{1}(\mathbb A_f)\). In a basis of
	\(\mathbb K^{2}\) which is orthogonal for \(s_0\) the diagonal matrices
	with entries \(c,1\), where \(c\) has norm one, are unitary; so
	\(\det\) is surjective, \(\det G(\Q)=\mathbb K^{1}(\Q)\), and
	\(\mathcal S=G^{1}(\mathbb A_f)\). The square of an element \(g\) of
	\(\mathcal H\) is the product of the scalar \(\det g\) and of an element
	of \(\mathcal S\). By Lemmas~\ref{lem:normal} and~\ref{lem:conj}(b),
	with \(W=\det U=\prod_\ell\det U_\ell\),
	\[
	\Mass(X)=m\bigl(G^{1}(\Q),G^{1}(\mathbb A_f),U^{1}\bigr)\cdot
	\frac{\#\bigl(\mathbb K^{1}(\mathbb A_f)/\mathbb K^{1}(\Q)W\bigr)}
	{\#(\mathbb K^{1}(\Q)\cap W)},\qquad U^{1}=\prod_\ell U_\ell^{1}.
	\]
	The finiteness conditions of Lemma~\ref{lem:normal} hold because
	\(\mathcal A_{-1}[X]\) is a finite groupoid and by
	Lemma~\ref{lem:torus}, once we know that \(W\) has finite index in
	\(W_{\max}\).

	We determine the local data. Let \(\ell\nmid2p\). By
	Lemma~\ref{lem:otherprimes}(a), \(T_\ell A_0\) has a hyperbolic basis
	over \(O_\ell\), and a basis in which \(H=\langle1,1\rangle\). By
	Lemmas~\ref{lem:quat}(b) and~\ref{lem:hyp}, \(\mathfrak B_\ell\) is
	split and \(U_\ell^{1}=\operatorname{SL}_2(\Z_\ell)\) is the group of
	the elements of reduced norm one of a maximal order; and
	\(\det U_\ell=O_\ell^{1}\), because the diagonal matrices with entries
	\(c,1\), \(c\in O_\ell^{1}\), are in \(U_\ell\). Let \(\ell=p\). By
	Lemma~\ref{lem:otherprimes}(b) the Dieudonn\'e module \(T_pA_0\) has a
	basis in which \(H=\langle1,u\rangle\), so again
	\(\det U_p=O_p^{1}\). If \(-u\) is a square, \(T_pA_0\) has a hyperbolic
	basis, \(\mathfrak B_p\) is split, and \(U_p^{1}\) has index \(p+1\) in
	the group of the elements of reduced norm one of a maximal order
	(Lemma~\ref{lem:hyp}). If \(-u\) is not a square,
	\(T_pA_0=\{v:H(v,v)\in\Z_p\}\), so \(\mathfrak B_p\) is a division
	algebra and \(U_p^{1}=\mathfrak B_p^{1}\) (Lemma~\ref{lem:quat}(c)).
	By Lemma~\ref{lem:quat}(d), \(\mathfrak B\) is ramified at exactly one
	of the primes \(2\) and \(p\).

	For the factor of \(G^{1}\), Lemma~\ref{lem:quat}(a) gives
	\(G^{1}(\Q)=\mathfrak B^{1}\) and
	\(G^{1}(\mathbb A_f)=\widehat{\mathfrak B}^{1}\). The order
	\(\mathfrak O'=\mathfrak B\cap\End(A_0)\) of \(\mathfrak B\) satisfies
	\((\mathfrak O'\otimes\Z_\ell)^{1}=U_\ell^{1}\) for all \(\ell\), by
	Tate's theorem. Let \(\mathfrak O\) be a maximal order containing
	\(\mathfrak O'\). By Lemma~\ref{lem:conj}(a),
	\(m(\mathfrak B^{1},\widehat{\mathfrak B}^{1},U^{1})
	=[\widehat{\mathfrak O}^{1}:U^{1}]\,
	m(\mathfrak B^{1},\widehat{\mathfrak B}^{1},\widehat{\mathfrak O}^{1})\),
	and the index is the product of the local indices, which are quotients
	of volumes. All maximal orders of \(\mathfrak B_\ell\) are conjugate
	under \(\mathfrak B_\ell^{\times}\), and by Lemma~\ref{lem:conj}(b)
	their groups of elements of reduced norm one have the same volume. So
	the local index is \(1\) for \(\ell\nmid2p\), it is \(i_p\) at \(p\) and
	\(i_2\) at \(2\), and
	\(m(\mathfrak B^{1},\widehat{\mathfrak B}^{1},
	\widehat{\mathfrak O}^{1})=(q-1)/24\) by Lemma~\ref{lem:normone}.

	As \(\det U_\ell=O_\ell^{1}\) for
	\(\ell\ne2\), Lemma~\ref{lem:torus} gives the last factor of
	\eqref{eq:genusmass}; the index \([O_2^{1}:W_2]\) is finite because
	\(W_2\) is open in the compact group \(O_2^{1}\).
\end{proof}

\begin{proof}[Proof of Theorem~\textup{\ref{thm:celltotals}}]
	We first compute five masses. The classes concerned are not empty:
	\(\mathsf{OO}\) and \(\mathsf{RO}\) contain products of elliptic curves,
	\(\mathsf{RR}_{\mathrm a}\), \(\mathsf N_{\mathrm a}\) and
	\(\mathsf R^{+}\) contain Weil restrictions (which exist by
	Lemma~\ref{lem:ssp2}), \(\mathsf N_{+}\) contains
	\(E\times E\) for \(E\) of type~II (its form is
	\(\langle u,u\rangle\), of determinant a square), and \(\mathsf O^{+}\)
	contains a quotient of a Weil restriction, since
	\(c(\mathsf R^{+}\to\mathsf O^{+})=1\). In each case we fix
	\((A_0,\theta_0)\) in the class, put \(\Lambda=T_2A_0\),
	\(V_2=\Lambda\otimes\Q\), and use Lemma~\ref{lem:genusmass} and its
	notation when \(\eps=-1\).

	\emph{\(\mathsf{OO}\), \(p\equiv3\pmod4\).} By the proof of
	Proposition~\ref{prop:twoadic}, \(\Lambda\) has an \(O_2\)-basis which is
	hyperbolic for \(H/2\), and one in which \(H/2=\langle1,1\rangle\).
	Replacing \(H\) by \(H/2\) does not change the unitary group, and
	\(\pi\in O_2^{\times}\). By Lemmas~\ref{lem:quat}(b) and~\ref{lem:hyp},
	\(\mathfrak B_2\) is split, so \(q=p\), and \(i_2=1\); also
	\(W_2=O_2^{1}\). As \(c_{\mathbb K}=1\) and \(h_{\mathbb K}=h\),
	\[
	\Mass(\mathsf{OO})=\frac{p-1}{24}\cdot\frac h2=\frac{(p-1)h}{48}.
	\]

	\emph{\(\mathsf{RO}\), \(p\equiv3\pmod8\).} Here
	\(\Lambda=R_2x\oplus O_2y\) with \(H=\langle1,2\rangle\), and
	\(\mathbb K_2\) is the unramified quadratic extension of \(\Q_2\). In
	\(H(ax+by,ax+by)=\operatorname N(a)+2\operatorname N(b)\) the first term
	has even and the second odd valuation, so
	\(M_0=\{v:H(v,v)\in\Z_2\}=O_2x\oplus O_2y\). By
	Lemma~\ref{lem:quat}(c), \(\mathfrak B_2\) is a division algebra, so
	\(q=2\), \(i_p=p+1\), and \(i_2\) is the index of \(U_2^{1}\) in the
	special unitary group, that is, the number of lattices in the orbit of
	\(\Lambda\) under this group. Let \(g\) be unitary. Then \(g\) stabilises
	\(M_0\), \(gx=a'x+b'y\) with \(a'\in O_2^{\times}\) and
	\(gy\in2O_2x+O_2^{\times}y\), because \(H(gx,gx)=1\) and
	\(H(gy,gy)=2\); hence
	\[
	g\Lambda=\Z_2a'x+2O_2x+O_2y,
	\]
	which depends only on the image of \(a'\) in
	\((O_2/2O_2)^{\times}=\F_4^{\times}\). So the orbit of \(\Lambda\) under
	the unitary group has at most three elements. Let \(\zeta\in O_2\) be a
	primitive cube root of unity; \(\operatorname N(\zeta)=1\), and the map
	\(x\mapsto\zeta x\), \(y\mapsto\zeta^{-1}y\) is unitary of determinant
	one and permutes three such lattices cyclically. Hence \(i_2=3\). The
	maps \(x\mapsto x\), \(y\mapsto cy\) with \(c\in O_2^{1}\) stabilise
	\(\Lambda\), so \(W_2=O_2^{1}\). Therefore
	\[
	\Mass(\mathsf{RO})=\frac1{24}\cdot3(p+1)\cdot\frac h2=\frac{(p+1)h}{16}.
	\]

	\emph{\(\mathsf N_{\mathrm a}\), \(p\equiv1\pmod4\).} Here
	\(\Lambda=R_2f_1\oplus R_2f_2\) with a hyperbolic basis, and
	\(\pi\in R_2^{\times}\). By Lemmas~\ref{lem:quat}(b)
	and~\ref{lem:hyp}, \(\mathfrak B_2\) is split, \(q=p\) and \(i_2=1\).
	We claim that \(W_2=\{u/\bar u:u\in R_2^{\times}\}\), which has index
	two in \(O_2^{1}\). For \(a\in\mathbb K_2^{\times}\) let \(D_a\) be the
	map \(f_1\mapsto af_1\), \(f_2\mapsto\bar a^{-1}f_2\); it is unitary of
	determinant \(a/\bar a\), and \(D_u\Lambda=\Lambda\) for
	\(u\in R_2^{\times}\). Conversely let \(g\) be unitary with
	\(g\Lambda=\Lambda\). By Hilbert's Theorem~90, \(\det g=a/\bar a\) with
	\(a\in\mathbb K_2^{\times}\) determined up to a factor in
	\(\Q_2^{\times}\); as \(2\) is \(\varpi^{2}\) times a unit, we may take
	\(a=\varpi^{k}u\) with \(u\in R_2^{\times}\) and \(k\in\{0,1\}\). Then
	\(g_1=gD_a^{-1}\) has determinant one, so it is a
	\(\Q_2\)-linear automorphism of \(V_0=\Q_2f_1\oplus\Q_2\pi f_2\), and
	\(D_a\Lambda=g_1^{-1}\Lambda\). Now
	\(\Lambda=R_2\otimes\Lambda_0\) with
	\(\Lambda_0=\Z_2f_1\oplus\Z_2\pi f_2\), so \(g_1^{-1}\Lambda\) is
	generated over \(R_2\) by its intersection with \(V_0\). If \(k=1\), then
	\(D_a\Lambda=\varpi R_2f_1\oplus\varpi^{-1}R_2f_2\), whose intersection
	with \(V_0\) is \(2\Z_2f_1\oplus\Z_2\pi f_2\), and this does not generate
	it. So \(k=0\) and \(\det g=u/\bar u\). The map \(a\mapsto a/\bar a\)
	identifies \(\mathbb K_2^{\times}/\Q_2^{\times}\) with the group of the
	elements of norm one, which are units, and
	\(\mathbb K_2^{\times}/\Q_2^{\times}R_2^{\times}\) has order two. This
	proves the claim. As \(c_{\mathbb K}=2\) and \(h_{\mathbb K}=h'\),
	\[
	\Mass(\mathsf N_{\mathrm a})=\frac{p-1}{24}\cdot\frac{h'}2\cdot\frac22
	=\frac{(p-1)h'}{48}.
	\]

	\emph{\(\mathsf N_{+}\), \(p\equiv1\pmod4\).} Here
	\(\Lambda=R_2e_1\oplus R_2e_2\) with \(H=\langle1,1\rangle\). By
	Lemma~\ref{lem:norms}(b) the norm of an element of \(\mathbb K_2\) of
	valuation \(0\), \(1\), \(\ge2\) (for the valuation with
	\(v(\varpi)=1\)) is congruent to \(1\) modulo \(4\), to \(2\) modulo
	\(4\), to \(0\) modulo \(4\). Let \(v=ae_1+be_2\), so that
	\(H(v,v)=\operatorname N(a)+\operatorname N(b)\), and let \(-k\) be the
	smaller of the valuations of \(a\) and \(b\). If \(k\le0\) then
	\(v\in\Lambda\). If \(k\ge1\), then
	\(\operatorname N(\varpi)^{k}H(v,v)=\operatorname N(\varpi^{k}a)
	+\operatorname N(\varpi^{k}b)\) is congruent to \(2\) modulo \(4\) if
	both valuations are \(-k\), and is odd otherwise; as
	\(\operatorname N(\varpi)\) is twice an odd number, \(H(v,v)\in\Z_2\)
	if and only if \(k=1\) and both valuations are \(-1\). Since the units
	of \(R_2\) are congruent to \(1\) modulo \(\varpi\), this shows that
	\[
	M_0=\{v:H(v,v)\in\Z_2\}=\Lambda+R_2f,\qquad
	f=\varpi^{-1}(e_1+e_2),
	\]
	a lattice containing \(\Lambda\) with index two. By
	Lemma~\ref{lem:quat}(c), \(\mathfrak B_2\) is a division algebra, so
	\(q=2\), \(i_p=p+1\), and \(i_2\) is the number of lattices in the orbit
	of \(\Lambda\) under the special unitary group. In the basis \(e_1,f\)
	of \(M_0\) the Gram matrix has entries \(1\), \(\bar\varpi^{-1}\),
	\(\varpi^{-1}\), \(2/\operatorname N(\varpi)\); so
	\(H(M_0,M_0)=\varpi^{-1}R_2\), hence
	\(\varpi M_0\subset M_0^{\#}\), and both have index four in \(M_0\)
	because \(M_0^{\#}\subset\Lambda^{\#}=\Lambda\subset M_0\) with indices
	two. So \(M_0^{\#}=\varpi M_0\). A unimodular lattice \(\Lambda'\) of
	\(V_2\) is contained in \(M_0\) and contains \(M_0^{\#}\), with index two
	in each direction. Conversely each of the three lattices \(\Lambda'\)
	with \(\varpi M_0\subset\Lambda'\subset M_0\) and of index two is an
	\(R_2\)-module, because \(R_2/\varpi R_2=\F_2\); it is of the form
	\(R_2v+\varpi M_0\), on which \(H\) has values in \(R_2\); and it has
	the same index in \(M_0\) as \(\Lambda\); so it is unimodular. All
	unimodular lattices of \(V_2\) are of class \(\mathsf N_{+}\), since
	\(\mathfrak d=1\) for this space; so the unitary group permutes the three
	lattices transitively. The diagonal matrices with entries \(c,1\),
	\(c\in O_2^{1}\), stabilise \(\Lambda\); hence \(W_2=O_2^{1}\), the
	unitary group is generated by the special unitary group and the
	stabiliser of \(\Lambda\), and \(i_2=3\). Therefore
	\[
	\Mass(\mathsf N_{+})=\frac1{24}\cdot3(p+1)\cdot\frac{h'}{2}\cdot\frac12
	=\frac{(p+1)h'}{32}.
	\]

	\emph{\(\eps=+1\): \(\mathsf O^{+}\) for \(p\equiv1\pmod4\), and
	\(\mathsf R^{+}\) for \(p\equiv3\pmod4\).} By
	\cite[(1.7) and \S3.1.2]{XueYu2020}, \(\mathfrak D=\End^{0}(A_0)\) is
	the quaternion algebra over \(\mathbb K_{+}\) ramified exactly at the
	two real places, and the Rosati involution is its canonical involution;
	so the group \(G\) of Lemma~\ref{lem:genus} satisfies
	\(G(\Q)=\mathfrak D^{1}\) and
	\(G(\mathbb A_f)=\widehat{\mathfrak D}^{1}\). For every prime \(\ell\),
	including \(\ell=2\) and \(\ell=p\), the module \(T_\ell A_0\) is free
	of rank two over the maximal order \(O_\ell\)
	(Proposition~\ref{prop:twoadic} and Lemma~\ref{lem:otherprimes}), and by
	Lemma~\ref{lem:sympl} its automorphism group \(U_\ell\) is the group of
	the elements of determinant one of
	\(\End_{O_\ell}(T_\ell A_0)\cong M_2(O_\ell)\). By Tate's theorem
	\(\End(A_0)\otimes\Z_\ell=\End_{O_\ell}(T_\ell A_0)\). Hence
	\(\mathfrak O=\End(A_0)\) is a maximal order of \(\mathfrak D\) and
	\(U=\widehat{\mathfrak O}^{1}\). By Lemmas~\ref{lem:genus}
	and~\ref{lem:normone} the mass of the class is \(z_p/4\).

	Finally, by Proposition~\ref{prop:neighbours},
	\(\Mass(Y)=\Mass(X)\,c(X\to Y)/c(Y\to X)\) whenever
	\(c(Y\to X)\ne0\). With Table~\ref{tab:neighbours} and
	\(X=\mathsf{OO}\) this gives, for \(p\equiv7\pmod8\), the factors
	\(6/2\), \(3/1\), \(1/1\) for \(\mathsf{RO}\),
	\(\mathsf{RR}_{\mathrm n}\), \(\mathsf{RR}_{\mathrm a}\), and for
	\(p\equiv3\pmod8\) the factors \(9/1\), \(3/1\) for
	\(\mathsf{RR}_{\mathrm n}\), \(\mathsf{RR}_{\mathrm a}\). With
	\(X=\mathsf N_{\mathrm a}\) it gives the factor \(3/2\) for
	\(\mathsf N_{-}\), and with \(X=\mathsf O^{+}\) the factor \(6\) or
	\(10\) for \(\mathsf R^{+}\), according as \(p\equiv1\) or
	\(5\pmod8\), that is, \(8-2\leg2p\).
\end{proof}

\begin{proof}[Proof of Theorem~\textup{\ref{thm:vertexpm}}]
	By Lemma~\ref{lem:Aeps}(b) and Proposition~\ref{prop:twoadic}(b), the
	mass of the superspecial curves with Weil polynomial
	\((T^{2}-\eps p)^{2}\) and Frobenius type \(\lambda\) is the mass of the
	objects of \(\mathcal A_\eps\) with cell \(\lambda\), minus the mass of
	the products and Weil restrictions with this cell. By
	Proposition~\ref{prop:twoadic} and Lemma~\ref{lem:Aeps}(c) the latter is
	\(N_1^{2}/8\), \(N_1N_2/4\), \(N_2^{2}/8\), \((p-1)/48\) for
	\(\lambda=1^{6},1^{4}2,1^{2}2^{2},2^{3}\) if \(\eps=-1\), and it is
	\((p-1)/48\) for \(\lambda=2^{3}\) and zero otherwise if \(\eps=+1\).
	Here \(N_2=h'\), \(N_1=h\) if \(p\equiv3\pmod4\) and
	\(N_1=0\) if \(p\equiv1\pmod4\) (Lemma~\ref{lem:ssp}), and for
	\(p\equiv3\pmod4\), \(h'=h\) if \(p\equiv7\pmod8\) and \(h'=3h\) if
	\(p\equiv3\pmod8\).

	Let \(p\equiv3\pmod4\) and \(\eps=-1\). By
	Theorem~\ref{thm:celltotals}, the cell \(1^{6}\) has mass
	\((p-1)h/48\); the cell \(1^{4}2\) has mass \(3(p-1)h/48\) or
	\((p+1)h/16\), that is, \(\bigl(p-\leg2p\bigr)h/16\); the cell
	\(1^{2}2^{2}\) has mass \(3(p-1)h/48\) or \(9(p-1)h/48\), that is,
	\((p-1)h'/16\); and the cell \(2^{3}\) has mass \((p-1)h/48\) or
	\(3(p-1)h/48\), that is, \((p-1)h'/48\). Subtracting \(h^{2}/8\),
	\(hh'/4\), \(h'^{2}/8\) and \((p-1)/48\) gives the four formulas of
	Theorem~\ref{thm:vertex}(a) for \(P=(T^{2}+p)^{2}\). For
	\(\eps=+1\) the only class is \(\mathsf R^{+}\), with cell \(2^{3}\) and
	mass \(z_p/4\), and we subtract \((p-1)/48\).

	Let \(p\equiv1\pmod4\). For \(\eps=-1\) the cell \(1^{2}2^{2}\) has mass
	\(\Mass(\mathsf N_{+})+\Mass(\mathsf N_{-})=ph'/16\), from which we subtract
	\(h'^{2}/8\), and the cell \(2^{3}\) has mass \((p-1)h'/48\), from which
	we subtract \((p-1)/48\). For \(\eps=+1\) the cell \(1^{6}\) has mass
	\(z_p/4\) and contains no product or Weil restriction, and the cell
	\(2^{3}\) has mass \(\bigl(8-2\leg2p\bigr)z_p/4\), from which we
	subtract \((p-1)/48\). These are the formulas of
	Theorem~\ref{thm:vertex}(b).

	The cells that do not occur in Table~\ref{tab:twoadic} for a given
	\(\eps\) and \(p\bmod4\) are empty; this is the last assertion.
\end{proof}

\begin{proof}[Proof of Corollary~\textup{\ref{cor:split1222}}]
	By Proposition~\ref{prop:twoadic} the Jacobian of \(C\) is of class
	\(\mathsf N_{+}\) or \(\mathsf N_{-}\), and by Table~\ref{tab:neighbours}
	its three quotients are of class \(\mathsf N_{+}\) in the first case,
	while in the second case one is of class \(\mathsf N_{-}\) and two are of
	class \(\mathsf N_{\mathrm a}\). It remains to find the masses of the
	products in the two classes. Let \(E_1,E_2\) be supersingular elliptic
	curves over \(\Fp\); they are of type~II, \(T_2E_i=R_2e_i\)
	(Lemma~\ref{lem:tate2}), and the form of \(E_1\times E_2\) is
	\(\langle-c_1,-c_2\rangle\) with
	\(c_i=\langle e_i,\pi e_i\rangle_i\) as in Proposition~\ref{prop:cells}(ii). So
	\(E_1\times E_2\) is of class \(\mathsf N_{+}\) if and only if
	\(c_1c_2\equiv1\pmod4\), that is, \(\delta(E_1)=\delta(E_2)\) in the
	notation of that proposition. As \(\delta(E^{\mathrm{tw}})=-\delta(E)\),
	each sign occurs for \(h'/2\) of the \(h'\) classes of curves. As in the
	proof of Lemma~\ref{lem:Aeps}(c), the products with
	\(\delta(E_1)=\delta(E_2)\) have mass
	\(\tfrac12\cdot2(h'/2)^{2}/4=h'^{2}/16\), and so have the others. We
	subtract these from \(\Mass(\mathsf N_{+})\) and \(\Mass(\mathsf N_{-})\).
\end{proof}

\begin{remark}\label{rem:vertexlit}
	The following comparisons with the literature are not used in the
	proofs.

	The sum of the masses of
	Theorem~\ref{thm:celltotals} over the classes with \(\eps=-1\) is
	\(8(p-1)h/48=(p-1)h/6\) for \(p\equiv7\pmod8\),
	\(\bigl(13(p-1)+3(p+1)\bigr)h/48=(8p-5)h/24\) for \(p\equiv3\pmod8\), and
	\(\bigl(2(p-1)+3(p-1)+3(p+1)\bigr)h'/96=(4p-1)h'/48\) for
	\(p\equiv1\pmod4\). This is the value \(M(2)\) of
	\cite[Thm.~2]{Ibukiyama1993} recalled in Section~\ref{sec:mass}. The
	ratios \(1:3:3:1\), \(1:9:3\) and \(2:3\) between the masses of the
	classes linked by Proposition~\ref{prop:neighbours} are the quotients of
	the local data at \(2\) listed in \cite[Prop.~4.4]{Ibukiyama1993}, whose
	proof is omitted there.

	For \(\pi^{2}=p\), the classes \(\mathsf O^{+}\) and
	\(\mathsf R^{+}\) are the genera \(\Lambda_1^{pp}\) and
	\(\Lambda_{16}^{pp}\) of \cite[Prop.~3.2.5]{XueYu2020} when
	\(p\equiv1\pmod4\). The class numbers of these genera are computed in
	\cite[(4.33) and (4.34)]{XueYu2020}, and that of the single genus for
	\(p\equiv3\pmod4\) in \cite[Thm.~1.1]{XueYu2020}; their terms in
	\(\zeta_{\mathbb K_{+}}(-1)\) are \(z_p/2\),
	\(\bigl(4-\leg2p\bigr)z_p\) and \(z_p/2\), twice the masses of
	Theorem~\ref{thm:celltotals}(iv). We did not find the masses themselves
	stated in that paper.

	The classes \(\mathsf{OO}\),
	\(\mathsf N_{+}\), \(\mathsf N_{-}\), \(\mathsf N_{\mathrm a}\) consist of
	lattices over the maximal order of \(\mathbb K\), and for these the
	quaternion algebra of Lemma~\ref{lem:quat}, its ramification and the
	local indices at \(p\) and at \(2\) are those of Hashimoto and Koseki
	\cite[\S2: Lemma~2.1, Prop.~2.3, 2.5, 2.7 and~2.10]{HashimotoKoseki1989}.
	The number \(h_{\mathbb K}/c_{\mathbb K}\) of Lemma~\ref{lem:torus} is
	the class number of the torus \(\mathbb K^{1}\) given in
	\cite[Prop.~5.9(1)]{XueYu2021}. For the classes \(\mathsf{RO}\), \(\mathsf{RR}_{\mathrm n}\),
	\(\mathsf{RR}_{\mathrm a}\), which are not modules over the maximal
	order, we did not find the masses in the literature apart from the
	local data of \cite{Ibukiyama1993} mentioned above.
\end{remark}

\begin{remark}\label{rem:vertexchecks}
	The statements of this appendix were compared with the enumeration of
	Section~\ref{sec:comp} for the \(37\) primes \(7\le p\le173\) (script
	\texttt{verify\_census.py}). For every class, the mass of the curves
	(Table~\ref{tab:vertex}) plus the mass of the products and Weil
	restrictions of Lemma~\ref{lem:Aeps}(c) is the mass of
	Theorem~\ref{thm:celltotals}. For each of the \(3754\) orbits of branch
	divisors with \(\pi^{2}=\pm p\), the numbers of rational maximal
	isotropic subgroups whose quotient has a given cell form a row of
	Table~\ref{tab:neighbours}, a decomposed quotient being counted in the
	cell given by Proposition~\ref{prop:twoadic}(b). The masses of
	Corollary~\ref{cor:split1222} agree for the \(19\) primes
	\(p\equiv1\pmod4\). The proofs do not use these comparisons.
\end{remark}

\section{Masses of curves with \texorpdfstring{\(\pi^{2}\neq\pm p\)}{pi\textasciicircum 2 not +-p}}\label{sec:vertexcm}

Parts (c), (d) and (e) of Theorem~\ref{thm:vertex}, proved here as
Theorem~\ref{thm:vertexcmcells}, concern the superspecial curves with \(\pi^{2}\ne\pm p\), whose Weil
polynomial is \(T^{4}+p^{2}\) or \(T^{4}\pm pT^{2}+p^{2}\)
(Theorem~\ref{thm:class}). For such a polynomial \(P\) the abelian surfaces
over \(\Fp\) with Weil polynomial \(P\) form one isogeny class, which is
simple and has a commutative endomorphism algebra, the biquadratic CM field
\(\Q(\pi)\). By a theorem of Centeleghe and Stix these surfaces correspond
to the ideals of the order \(\Z[\pi,p/\pi]\), and by a theorem of Xue, Yang
and Yu the superspecial ones correspond to the ideals of
\(\Z[\pi,\pi^{2}/p]\), an order which is maximal or of index four in the
maximal order. We add the principal polarisations to this description
(Proposition~\ref{prop:dict}): once one polarised surface is fixed, the
principally polarised ones correspond to the pairs of an ideal \(I\) and a
totally positive element \(z\) of the real quadratic subfield such that
\(zI\bar I\) is a fixed ideal, and the classes of these pairs form a
principal homogeneous space under a finite group. The order of this group
is found by class field theory and Kuroda's class number formula. This
gives the number of classes and the mass of the principally polarised
superspecial surfaces with a given Weil polynomial and a given endomorphism
ring (Theorem~\ref{thm:vertexcm}); all of them have the same automorphism
group. The Weil restrictions, known from Proposition~\ref{prop:nonjac}, are
then subtracted.

From the literature we use: the equivalence of categories of Centeleghe and
Stix \cite[Thm.~25]{CentelegheStix2015}; the description of the superspecial
objects and of the orders by Xue, Yang and Yu \cite{XueYangYu2016}, of which
we recall the proofs; Kuroda's class number formula in the form of
\cite[(1.5), p.~247]{Lemmermeyer1994}; the theorems of class field theory as stated
in \cite[\S8]{Cox1989}, and the class number of a quadratic order
\cite[Thm.~7.24]{Cox1989}; the properties of the Weil pairing and the
Riemann--Roch theorem for abelian varieties
\cite[I, Thm.~11.1 and \S13]{MilneAV}, \cite[\S16]{Mumford2008}; the criterion of
Oort for superspeciality \cite[Thm.~2 and Remark~3]{Oort1975},
\cite[proof of Prop.~5.1]{XueYangYu2016}; and, for the polynomial
\(T^{4}-pT^{2}+p^{2}\) with \(p\equiv1\pmod3\), the theorem of Howe, Maisner,
Nart and Ritzenthaler that the isogeny class contains no principally
polarised surface \cite[Thm.~2]{HoweMaisnerNartRitzenthaler2008}. For
\(p\equiv2\pmod3\) the existence of such a surface is taken from
Proposition~\ref{prop:nonjac}.

\subsection*{Notation}
Throughout this appendix \(p\ge7\) is a prime and \(P\) is one of the
polynomials
\[
P_c=T^{4}+p^{2},\qquad P_d=T^{4}+pT^{2}+p^{2},\qquad
P_e=T^{4}-pT^{2}+p^{2};
\]
we put \(\eps=+1\) for \(P_d\) and \(\eps=-1\) for \(P_e\). Let
\(\mathbb L=\Q[T]/(P)\), a field of degree four by Lemma~\ref{lem:irred},
let \(\pi\) be the class of \(T\), and \(u=\pi^{2}/p\). Since
\(T^{4}P(p/T)=p^{2}P(T)\), there is an automorphism \(x\mapsto\bar x\) of
\(\mathbb L\) with \(\bar\pi=p/\pi\), and it is an involution. A complex root
\(z\) of \(P\) satisfies \(z^{4}=-p^{2}\), or \(z^{2}=\eps p\zeta\) with
\(\zeta^{3}=1\); hence \(\lvert z\rvert^{2}=p\) and \(p/z\) is the complex
conjugate of \(z\). So \(x\mapsto\bar x\) is the complex conjugation in every
embedding of \(\mathbb L\) into the complex numbers, \(\mathbb L\) is a CM
field, and we write \(\mathbb L_{+}\) for the field of the elements fixed by
the conjugation, and \(\operatorname N(x)=x\bar x\). An element of
\(\mathbb L_{+}\) is \emph{totally positive} if it is positive in the two
real embeddings of \(\mathbb L_{+}\). We have
\(\bar u=u^{-1}\) and \(\bar\pi=\pi u^{-1}\).

For \(P_c\) we have \(u^{2}=-1\), \((\pi(1+u))^{2}=2pu^{2}=-2p\) and
\((\pi(1-u))^{2}=-2pu^{2}=2p\), and \(\pi(1-u)\) is fixed by the conjugation.
For \(P_d\) and \(P_e\) the element \(\omega=\eps u\) satisfies
\(\omega^{2}+\omega+1=0\); on a Jacobian it is the \(\omega\) of
Theorem~\ref{thm:class}. We put \(\pi_0=\omega\pi\) and
\(\sqrt{-3}=2\omega+1\); then \(\pi_0^{2}=\omega^{2}\cdot\eps p\,\omega=\eps p\)
and \(\bar\pi_0=\omega^{-1}\pi\,(\eps\omega)^{-1}=\eps\pi_0\). Hence
\(\mathbb L\) is biquadratic. Besides \(\Q(u)\) and \(\mathbb L_{+}\) it
contains one quadratic field \(\mathbb E\), which is imaginary; we write
\(\Delta\) for its discriminant:
\begin{equation}\label{eq:cmfields}
\begin{array}{@{}llll@{}}
P_c\colon & \mathbb L_{+}=\Q(\sqrt{2p}), & \mathbb E=\Q(\sqrt{-2p}), &
\Delta=-8p;\\[2pt]
P_d\colon & \mathbb L_{+}=\Q(\sqrt{p}), & \mathbb E=\Q(\sqrt{-3p}), &
\Delta=-3p\ \ (p\equiv1\bmod4),\quad -12p\ \ (p\equiv3\bmod4);\\[2pt]
P_e\colon & \mathbb L_{+}=\Q(\sqrt{3p}), & \mathbb E=\Q(\sqrt{-p}), &
\Delta=-p\ \ (p\equiv3\bmod4),\quad -4p\ \ (p\equiv1\bmod4).
\end{array}
\end{equation}
The class number of \(\mathbb E\) is \(h(\Delta)\), in the notation of
Section~\ref{sec:mass}; for \(P_e\), \(\mathbb E\) is the field
\(\mathbb K\) of Appendix~\ref{sec:vertexpm}. We write \(h_{\mathbb L}\) and
\(h_{\mathbb L_{+}}\) for the class numbers of \(\mathbb L\) and
\(\mathbb L_{+}\), \(\operatorname{Cl}(\mathbb L)\) for the class group of
\(\mathbb L\), \(\operatorname{Cl}^{+}(\mathbb L_{+})\) for the narrow class
group of \(\mathbb L_{+}\) (fractional ideals modulo principal ideals with a
totally positive generator), \(U_{\mathbb L}\) and \(U_{+}\) for the unit
groups of \(\mathbb L\) and \(\mathbb L_{+}\), and \(U_{++}\) for the group
of totally positive units of \(\mathbb L_{+}\).

A \emph{lattice} in \(\mathbb L\) is a finitely generated subgroup of rank
four. For lattices \(X\), \(Y\) we put \((X:Y)=\{x\in\mathbb L:xY\subset X\}\)
and \(X^{*}=\{x\in\mathbb L:\operatorname{Tr}_{\mathbb L/\Q}(xX)\subset\Z\}\),
the trace dual. We consider the orders
\[
R_{\min}=\Z[\pi,\bar\pi]\ \subset\ R_{\mathrm{sp}}=\Z[\pi,\pi^{2}/p]
=\Z[u]\oplus\Z[u]\pi\ \subset\ \mathcal O_{\mathbb L},
\]
where \(\mathcal O_{\mathbb L}\) is the maximal order of \(\mathbb L\); the
first inclusion holds because \(\bar\pi=\pi u^{-1}\), and all three orders
are stable under the conjugation. A \emph{fractional ideal} of an order
\(\mathcal R\) is a lattice \(X\) with \(\mathcal RX\subset X\); it is
\emph{invertible} if \(XY=\mathcal R\) for some lattice \(Y\), and then
\((X:X)=\mathcal R\).

For an abelian variety \(A\) over \(\Fp\) we write \(\pi_A\) for its
Frobenius endomorphism, \(A^{\vee}\) for the dual abelian variety and
\(\hat f\) for the dual of a homomorphism \(f\); \(\Hom\) and \(\End\) refer
to homomorphisms defined over \(\Fp\), and
\(\Hom^{0}=\Hom\otimes\Q\), \(\End^{0}=\End\otimes\Q\). A \emph{polarisation}
of \(A\) is a homomorphism \(\nu\colon A\to A^{\vee}\) over \(\Fp\) which
over \(\kbar\) is of the form \(\phi_{\mathcal N}\colon
a\mapsto t_a^{*}\mathcal N\otimes\mathcal N^{-1}\) for an ample invertible
sheaf \(\mathcal N\) \cite[I, \S11]{MilneAV}; for a prime \(\ell\ne p\) its
Weil form on \(T_\ell A\) is \(E_\nu(a,a')=e_\ell(a,\nu a')\), with values in
\(\Z_\ell(1)\cong\Z_\ell\), where \(e_\ell\) is the pairing between
\(T_\ell A\) and \(T_\ell A^{\vee}\). As in Appendix~\ref{sec:vertexpm},
\((A,\theta)\) denotes a principally polarised abelian surface over \(\Fp\).

Let \(\mathcal V_P\) be the category of the abelian varieties over \(\Fp\)
with Weil polynomial \(P\), and \(\mathcal A_P\) the groupoid of the
principally polarised abelian surfaces \((A,\theta)\) over \(\Fp\) with Weil
polynomial \(P\), the morphisms being the isomorphisms over \(\Fp\) that
respect the polarisations. For \(B\in\mathcal V_P\) there is a unique
isomorphism \(i_B\colon\mathbb L\to\End^{0}(B)\) with \(i_B(\pi)=\pi_B\)
(see the paragraph before Lemma~\ref{lem:cs}). For an order \(\mathcal R\) of \(\mathbb L\) we write
\(\mathcal V_P[\mathcal R]\) for the full subcategory of the \(B\) with
\(\End(B)=i_B(\mathcal R)\), and \(\mathcal A_P[\mathcal R]\) for the full
subgroupoid of the \((A,\theta)\) with \(A\in\mathcal V_P[\mathcal R]\).

\subsection*{Statements}

\begin{lemma}\label{lem:orders}
	\begin{enumerate}[label=\textup{(\alph*)},leftmargin=*]
		\item \(R_{\mathrm{sp}}=\mathcal O_{\mathbb L}\) for \(P_c\), and for
		\(P_d\), \(P_e\) when \(\eps p\equiv3\pmod4\). If
		\(\eps p\equiv1\pmod4\), that is, for \(P_d\) with \(p\equiv1\pmod4\)
		and for \(P_e\) with \(p\equiv3\pmod4\), put \(\delta=(1+\pi_0)/2\).
		Then \(\mathcal O_{\mathbb L}=\Z[\omega]\oplus\Z[\omega]\delta\),
		\(R_{\mathrm{sp}}=\Z[\omega]+2\mathcal O_{\mathbb L}\),
		\([\mathcal O_{\mathbb L}:R_{\mathrm{sp}}]=4\), and
		\(R_{\mathrm{sp}}\), \(\mathcal O_{\mathbb L}\) are the only orders of
		\(\mathbb L\) containing \(R_{\mathrm{sp}}\).
		\item \(R_{\mathrm{sp}}^{*}=g^{-1}R_{\mathrm{sp}}\), with \(g=4\pi\) for
		\(P_c\) and \(g=2\pi_0\sqrt{-3}\) for \(P_d\), \(P_e\). In particular
		the trace duals of \(R_{\mathrm{sp}}\) and of
		\(\mathcal O_{\mathbb L}\) are invertible ideals.
		\item The roots of unity of \(\mathbb L\) form the group
		\(\langle u,-1\rangle\), of order \(w=4\) for \(P_c\) and \(w=6\) for
		\(P_d\), \(P_e\). It is contained in \(R_{\mathrm{sp}}\).
		\item The finite primes of \(\mathbb L_{+}\) that ramify in
		\(\mathbb L\) are: the prime above \(2\) for \(P_c\); the primes above
		\(3\) for \(P_d\); none for \(P_e\).
	\end{enumerate}
\end{lemma}

\begin{theorem}\label{thm:vertexcm}
	Let \(p\ge7\) and let \(P\) be \(P_c\), \(P_d\) or \(P_e\).
	\begin{enumerate}[label=\textup{(\roman*)},leftmargin=*]
		\item An abelian surface \(A\) over \(\Fp\) with Weil polynomial \(P\) is
		superspecial if and only if \(\End(A)\) contains
		\(i_A(R_{\mathrm{sp}})\). Hence the principally polarised superspecial
		surfaces over \(\Fp\) with Weil polynomial \(P\) are the objects of
		\(\mathcal A_P[R_{\mathrm{sp}}]\) and of
		\(\mathcal A_P[\mathcal O_{\mathbb L}]\).
		\item Let \(\mathcal R\) be \(R_{\mathrm{sp}}\) or
		\(\mathcal O_{\mathbb L}\). Every object of \(\mathcal A_P[\mathcal R]\)
		has automorphism group \(\langle u,-1\rangle\), of order \(w\). The
		number of isomorphism classes of
		\(\mathcal A_P[\mathcal O_{\mathbb L}]\) is
		\[
		\begin{cases}
		h(\Delta)/2, & P=P_c\ \text{or}\ P=P_d,\\
		h(\Delta), & P=P_e\ \text{and}\ p\equiv2\pmod3,\\
		0, & P=P_e\ \text{and}\ p\equiv1\pmod3,
		\end{cases}
		\]
		and, when \(R_{\mathrm{sp}}\ne\mathcal O_{\mathbb L}\), the number of
		isomorphism classes of \(\mathcal A_P[R_{\mathrm{sp}}]\) is given by the
		same expression with \(h(4\Delta)\) in place of \(h(\Delta)\). In each
		case the mass is the number of classes divided by \(w\).
		\item Every object of \(\mathcal A_P[R_{\mathrm{sp}}]\) or of
		\(\mathcal A_P[\mathcal O_{\mathbb L}]\) is either the Jacobian of a
		superspecial curve \(C/\Fp\) of genus two with Weil polynomial \(P\),
		which is unique up to \(\Fp\)-isomorphism and satisfies
		\(\Aut_{\Fp}(C)\cong\Aut_{\Fp}(A,\theta)\), or the Weil restriction of a
		supersingular elliptic curve over \(\F_{p^{2}}\). The Weil restrictions
		have endomorphism ring \(R_{\mathrm{sp}}\). If \(P\) is \(P_d\) or
		\(P_e\), the Frobenius type of \(C\) is \(3^{2}\) when
		\(\End(\Jac C)=\mathcal O_{\mathbb L}\ne R_{\mathrm{sp}}\), and \(6\)
		otherwise.
	\end{enumerate}
\end{theorem}

Assertion (i) and Lemma~\ref{lem:orders}(a) are due to Xue, Yang and Yu:
by \cite[Prop.~5.1]{XueYangYu2016} the superspecial abelian varieties in the
isogeny class of a Weil \(p\)-number without real conjugate are the lattices
over \(\Z[\pi,\pi^{2}/p]\), for which they write \(R_{\mathrm{sp}}\), and the
orders containing it are listed in \cite[(5.8) and the tables before (5.9)
and (5.10)]{XueYangYu2016} for our three Weil numbers, which are
\(\sqrt p\,\zeta_8\), \(\sqrt p\,\zeta_3\) and \(\sqrt p\,\zeta_{12}\) in
their notation. We recall the
proofs below because they are short.

\begin{theorem}\label{thm:vertexcmcells}
	Parts \textup{(c)}, \textup{(d)} and \textup{(e)} of
	Theorem~\textup{\ref{thm:vertex}} hold: for every prime \(p\ge7\) the
	masses of the superspecial curves of genus two over \(\Fp\) with Weil
	polynomial \(T^{4}+p^{2}\), \(T^{4}+pT^{2}+p^{2}\) or
	\(T^{4}-pT^{2}+p^{2}\), by Frobenius type, are those stated there, and the
	pairs of a Frobenius type and a residue class of \(p\) that are not listed
	have mass zero.
\end{theorem}

\subsection*{Orders and ideals}

\begin{proof}[Proof of Lemma~\textup{\ref{lem:orders}}]
	(a) The ring \(\Z[u]\) is \(\Z[\sqrt{-1}]\) for \(P_c\) and \(\Z[\omega]\)
	for \(P_d\), \(P_e\), and \(R_{\mathrm{sp}}=\Z[u][X]/(f)\) with
	\(f=X^{2}-pu\), the class of \(X\) being \(\pi\); for \(P_d\), \(P_e\) we
	may also take \(f=X^{2}-\eps p\), the class of \(X\) being \(\pi_0\),
	because \(\Z[u,\pi]=\Z[\omega,\pi_0]\). Such an order is maximal at a prime
	\(\mathfrak q\) of \(\Z[u]\) if \(f\) or \(f(X+1)\) is an Eisenstein
	polynomial at \(\mathfrak q\), or if the discriminant of \(f\) is a unit at
	\(\mathfrak q\). At the primes above \(p\), which is unramified in
	\(\Z[u]\), \(f\) is an Eisenstein polynomial. At the primes not dividing
	\(2p\) the discriminant \(4pu\), resp.\ \(4\eps p\), is a unit. For \(P_c\)
	the prime above \(2\) is generated by \(1+u\), and
	\(f(X+1)=X^{2}+2X+(1-pu)\) is an Eisenstein polynomial there:
	\(2=-u(1+u)^{2}\), and \(1-pu=-u(1+u)-(p-1)u\) is divisible by \(1+u\)
	exactly once, \(p-1\) being even. For \(P_d\),
	\(P_e\) the prime \(2\) is inert in \(\Z[\omega]\), and
	\(f(X+1)=X^{2}+2X+(1-\eps p)\) is an Eisenstein polynomial at \(2\) if and
	only if \(\eps p\equiv3\pmod4\). This proves the first assertion. Let
	\(\eps p\equiv1\pmod4\). Then \(\delta\) is a root of
	\(X^{2}-X-(\eps p-1)/4\), whose discriminant \(\eps p\) is odd, so
	\(\Z[\omega][\delta]\) is maximal at \(2\); it coincides with
	\(R_{\mathrm{sp}}\) at the other primes, where \(R_{\mathrm{sp}}\) is
	maximal. Hence \(\mathcal O_{\mathbb L}=\Z[\omega]\oplus\Z[\omega]\delta\)
	and, as \(\pi_0=2\delta-1\),
	\(R_{\mathrm{sp}}=\Z[\omega]\oplus2\Z[\omega]\delta
	=\Z[\omega]+2\mathcal O_{\mathbb L}\). The quotient
	\(\mathcal O_{\mathbb L}/R_{\mathrm{sp}}\) is isomorphic to
	\(\Z[\omega]/2\), a simple \(\Z[\omega]\)-module with four elements, and an
	order containing \(R_{\mathrm{sp}}\) is a \(\Z[\omega]\)-submodule of
	\(\mathcal O_{\mathbb L}\).

	(b) By Euler's formula the dual of \(\Z[u][X]/(f)\) for the trace form
	of \(\mathbb L\) over \(\Q(u)\) is generated by the inverse of \(f'\) at the
	class of \(X\), which is \(2\pi\), resp.\ \(2\pi_0\); the trace dual of
	\(\Z[\sqrt{-1}]\) is generated by \(1/2\), that of \(\Z[\omega]\) by
	\(1/\sqrt{-3}\), and \(x\in R_{\mathrm{sp}}^{*}\) if and only if
	\(\operatorname{Tr}_{\mathbb L/\Q(u)}(xR_{\mathrm{sp}})\) lies in the trace
	dual of \(\Z[u]\). The trace dual of the maximal order is the inverse of
	the different.

	(c) The roots of unity of \(\mathbb L\) generate a cyclotomic field
	\(\Q(\zeta_m)\subset\mathbb L\) whose degree divides four. If the degree
	is four, \(\mathbb L=\Q(\zeta_m)\) with \(m\in\{5,8,10,12\}\); but
	\(\Q(\zeta_5)\) is cyclic over \(\Q\), and the real subfields
	\(\Q(\sqrt2)\) and \(\Q(\sqrt3)\) of \(\Q(\zeta_8)\) and \(\Q(\zeta_{12})\)
	are different from \(\mathbb L_{+}\). If the degree is two,
	\(\Q(\zeta_m)\) is \(\Q(\sqrt{-1})\) or \(\Q(\sqrt{-3})\) and is one of
	\(\Q(u)\), \(\mathbb E\); it is \(\Q(u)\) because \(\Delta\ne-3,-4\). So
	the roots of unity of \(\mathbb L\) are those of \(\Q(u)\).

	(d) \(\mathbb L=\mathbb L_{+}(u)\), and \(u\) is a root of \(X^{2}+1\) or
	of \(X^{2}\pm X+1\), of discriminant \(-4\), resp.\ \(-3\); so only primes
	above \(2\), resp.\ \(3\), can ramify. For \(P_c\) the prime \(2\) ramifies
	in the three quadratic subfields of \(\mathbb L\) (their discriminants are
	\(-4\), \(8p\), \(-8p\)); hence its inertia group in
	\(\operatorname{Gal}(\mathbb L/\Q)\) fixes none of them and has order four,
	while \(2\) has ramification index two in \(\mathbb L_{+}\). For \(P_d\)
	the prime \(3\) is unramified in \(\mathbb L_{+}=\Q(\sqrt p)\) and ramified
	in \(\Q(u)\), so the primes of \(\mathbb L_{+}\) above \(3\) ramify in
	\(\mathbb L\). For \(P_e\) the prime \(3\) is unramified in
	\(\mathbb E=\Q(\sqrt{-p})\), so its ramification index in \(\mathbb L\) is
	two, and it is also two in \(\mathbb L_{+}=\Q(\sqrt{3p})\).
\end{proof}

\begin{lemma}\label{lem:classno}
	Let \(q_{\mathbb L}=[U_{\mathbb L}:\langle u,-1\rangle U_{+}]\).
	\begin{enumerate}[label=\textup{(\alph*)},leftmargin=*]
		\item \(\operatorname N(U_{\mathbb L})\) contains \(U_{+}^{2}\) with
		index \(q_{\mathbb L}\), and \(q_{\mathbb L}\in\{1,2\}\).
		\item \(h_{\mathbb L}=\tfrac12\,q_{\mathbb L}\,h_{\mathbb L_{+}}\,h(\Delta)\).
		\item The relative norm of ideals induces a homomorphism
		\(\bar{\operatorname N}\colon\operatorname{Cl}(\mathbb L)\to
		\operatorname{Cl}^{+}(\mathbb L_{+})\), whose image has index \(r\),
		where \(r=2\) for \(P_e\) and \(r=1\) for \(P_c\), \(P_d\).
		\item If \(\Delta\) is odd, then
		\(h(4\Delta)=\bigl(2-\leg{\Delta}{2}\bigr)h(\Delta)\), where
		\(\leg{\Delta}{2}\) is \(1\) if \(\Delta\equiv1\pmod8\) and \(-1\) if
		\(\Delta\equiv5\pmod8\).
	\end{enumerate}
\end{lemma}

\begin{proof}
	(a) \(\operatorname N(\eta)=\eta^{2}\) for \(\eta\in U_{+}\), and
	\(\operatorname N\) is trivial on the roots of unity. If
	\(\operatorname N(a)=\eta^{2}\) with \(a\in U_{\mathbb L}\) and
	\(\eta\in U_{+}\), all the conjugates of the algebraic integer \(a/\eta\)
	have absolute value one, so \(a/\eta\) is a root of unity. Hence
	\(\operatorname N\) induces an isomorphism of
	\(U_{\mathbb L}/\langle u,-1\rangle U_{+}\) with
	\(\operatorname N(U_{\mathbb L})/U_{+}^{2}\), a subgroup of
	\(U_{++}/U_{+}^{2}\). As \(U_{+}\) is the product of \(\{\pm1\}\) and an
	infinite cyclic group and \(-1\notin U_{++}\), the group
	\(U_{++}/U_{+}^{2}\) has at most two elements.

	(b) For an extension of number fields with Galois group
	\((\Z/2)^{2}\), Kuroda's formula in the form of
	\cite[(1.5), p.~247]{Lemmermeyer1994} expresses the class number of the upper
	field through those of the base field \(k\) and of the three intermediate
	fields. For \(k=\Q\) and the upper field \(\mathbb L\), in the notation
	used there, \(d=1\) (the infinite place of \(\Q\) ramifies), \(\kappa=0\)
	and \(v=0\), and the formula reads
	\(h_{\mathbb L}=\tfrac12\,q\,h_{\Q(u)}\,h_{\mathbb E}\,h_{\mathbb L_{+}}\),
	where \(q\) is the index in \(U_{\mathbb L}\) of the group generated by the
	units of the three quadratic subfields. These units are the roots of unity
	of \(\Q(u)\), \(\pm1\) (as \(\Delta<-4\)) and \(U_{+}\), so
	\(q=q_{\mathbb L}\) by Lemma~\ref{lem:orders}(c); and \(h_{\Q(u)}=1\),
	\(h_{\mathbb E}=h(\Delta)\).

	(c) We use the theorems of class field theory as they are stated in
	\cite[\S8]{Cox1989}. Let \(I\) be the group of fractional ideals of
	\(\mathbb L_{+}\), \(I_{+}\) the subgroup of the principal ideals with a
	totally positive generator, and \(Y=\operatorname N(I_{\mathbb L})\,I_{+}\),
	where \(\operatorname N(I_{\mathbb L})\) is the group of the relative norms
	of the fractional ideals of \(\mathbb L\). The norm of a principal ideal
	\((a)\) is generated by \(a\bar a\), which is totally positive; so
	\(\bar{\operatorname N}\) is well defined, and the index of its image is
	\([I:Y]\). Let \(\mathfrak m\) be the modulus of \(\mathbb L_{+}\) formed
	by its two real places, so that \(I_{+}\) is the group denoted by
	\(P_{1}(\mathfrak m)\) in \cite{Cox1989}, and let \(\mathfrak f\) be the
	conductor of \(\mathbb L/\mathbb L_{+}\). It is divisible exactly by the
	primes that ramify \cite[Thm.~8.5]{Cox1989}, among them the two real
	places; so \(\mathfrak m\) divides \(\mathfrak f\). By
	\cite[Thm.~8.5 and~(8.21)]{Cox1989} the kernel of the Artin map of
	\(\mathbb L/\mathbb L_{+}\) on the group \(I(\mathfrak f)\) of the ideals
	prime to \(\mathfrak f\) is
	\(\operatorname N(I_{\mathbb L}(\mathfrak f))\,P_{1}(\mathfrak f)\), where
	\(I_{\mathbb L}(\mathfrak f)\) is the group of the fractional ideals of
	\(\mathbb L\) prime to \(\mathfrak f\) and
	\(P_{1}(\mathfrak f)\subset I_{+}\); this kernel is contained in
	\(Y\cap I(\mathfrak f)\).

	If no finite prime ramifies, then \(\mathfrak f=\mathfrak m\),
	\(I(\mathfrak f)=I\), the kernel of the Artin map is \(Y\), and the Artin
	map is surjective \cite[Thm.~8.2]{Cox1989}; so \([I:Y]=2\). In general, by
	the existence theorem \cite[Thm.~8.6]{Cox1989} there is an abelian
	extension \(\mathbb L'\) of \(\mathbb L_{+}\), unramified at the finite
	primes, such that \(Y\) is the kernel of its Artin map on \(I\); then
	\([\mathbb L':\mathbb L_{+}]=[I:Y]\), this map being surjective, and the
	kernel of its Artin map on
	\(I(\mathfrak f)\) is \(Y\cap I(\mathfrak f)\), which contains the kernel
	of the Artin map of \(\mathbb L\). By \cite[Cor.~8.7]{Cox1989},
	\(\mathbb L'\subset\mathbb L\). If a finite prime ramifies in
	\(\mathbb L\), this forces \(\mathbb L'=\mathbb L_{+}\) and \(Y=I\). By
	Lemma~\ref{lem:orders}(d) a finite prime ramifies for \(P_c\) and \(P_d\),
	and none for \(P_e\).

	(d) This is \cite[Thm.~7.24]{Cox1989} for the order of conductor \(2\) in
	\(\mathbb E\), whose units are \(\pm1\) because \(\Delta<-4\).
\end{proof}

Let \(\mathcal R\) be \(R_{\mathrm{sp}}\) or \(\mathcal O_{\mathbb L}\), and
let \(\mathfrak c\) be an invertible fractional ideal of \(\mathcal R\) with
\(\bar{\mathfrak c}=\mathfrak c\). We write
\(\mathcal I(\mathcal R,\mathfrak c)\) for the groupoid whose objects are the
pairs \((I,z)\), with \(I\) an invertible fractional ideal of \(\mathcal R\)
and \(z\in\mathbb L_{+}\) totally positive such that
\[
z\,I\bar I=\mathfrak c ;
\]
a morphism \((I,z)\to(I',z')\) is an element \(a\in\mathbb L^{\times}\) with
\(aI=I'\) and \(z=a\bar a\,z'\). Let \(\mathcal U(\mathcal R)\) be the set of
the pairs \((J,s)\), with \(J\) an invertible fractional ideal of
\(\mathcal R\) and \(s\in\mathbb L_{+}\) totally positive such that
\(J\bar J=s\mathcal R\), modulo the relation \((J,s)\sim(aJ,a\bar a\,s)\) for
\(a\in\mathbb L^{\times}\). It is an abelian group for the multiplication of
ideals and of elements, and it is finite: there are finitely many classes of
ideals \(J\), and for a given \(J\) the element \(s\) is determined up to a
totally positive unit, while the squares of the units of \(\mathbb L_{+}\)
contained in \(\mathcal R\) are of the form \(a\bar a\).

\begin{lemma}\label{lem:ideals}
	\begin{enumerate}[label=\textup{(\alph*)},leftmargin=*]
		\item Every object of \(\mathcal I(\mathcal R,\mathfrak c)\) has
		automorphism group \(\langle u,-1\rangle\). The set of isomorphism
		classes of \(\mathcal I(\mathcal R,\mathfrak c)\) is empty or a principal
		homogeneous space under \(\mathcal U(\mathcal R)\), for the action
		\((J,s)\cdot(I,z)=(JI,z/s)\). Hence the number of classes is \(0\) or
		\(\#\mathcal U(\mathcal R)\), and the mass is this number divided by
		\(w\).
		\item \(\#\mathcal U(\mathcal O_{\mathbb L})=r\,h(\Delta)/2\), with \(r\)
		as in Lemma~\textup{\ref{lem:classno}(c)}.
		\item Let \(\mathfrak c_{+}\) be a fractional ideal of \(\mathbb L_{+}\).
		Then \(\mathcal I(\mathcal O_{\mathbb L},
		\mathfrak c_{+}\mathcal O_{\mathbb L})\) is non-empty if and only if the
		narrow class of \(\mathfrak c_{+}\) lies in the image of
		\(\bar{\operatorname N}\). This always holds for \(P_c\) and \(P_d\).
		\item Let \(R_{\mathrm{sp}}\ne\mathcal O_{\mathbb L}\), and let
		\(\mathfrak c\) be an invertible fractional ideal of
		\(R_{\mathrm{sp}}\) with \(\bar{\mathfrak c}=\mathfrak c\). The functor
		\[
		\mathcal I(R_{\mathrm{sp}},\mathfrak c)\longrightarrow
		\mathcal I(\mathcal O_{\mathbb L},\mathfrak c\,\mathcal O_{\mathbb L}),
		\qquad (I,z)\longmapsto(I\mathcal O_{\mathbb L},z),
		\]
		is faithful, and each of its fibres, in the sense of
		Lemma~\textup{\ref{lem:fibre}}, has exactly
		\(n=2-\leg{\Delta}{2}\) elements. Hence the number of classes of
		\(\mathcal I(R_{\mathrm{sp}},\mathfrak c)\) is \(n\) times that of
		\(\mathcal I(\mathcal O_{\mathbb L},\mathfrak c\,\mathcal O_{\mathbb L})\).
	\end{enumerate}
\end{lemma}

\begin{proof}
	(a) If \(aI=I\) with \(I\) invertible, then \(a\) is a unit of
	\(\mathcal R\); if moreover \(a\bar a=1\), it is a root of unity, as in the
	proof of Lemma~\ref{lem:classno}(a); and the roots of unity lie in
	\(\mathcal R\) and are automorphisms (Lemma~\ref{lem:orders}(c)). The
	action is well defined on classes. It is transitive: for \((I,z)\) and
	\((I',z')\) put \(J=I'I^{-1}\) and \(s=z/z'\); then
	\(J\bar J=(\mathfrak c/z')(z/\mathfrak c)=s\mathcal R\). It is free: if
	\(a\) is a morphism from \((JI,z/s)\) to \((I,z)\), then
	\(J=a^{-1}\mathcal R\) and \(s=(a\bar a)^{-1}\), so \((J,s)\sim(\mathcal
	R,1)\). The last assertion follows because all automorphism groups have
	order \(w\).

	(b) For a fractional ideal \(J\) of \(\mathcal O_{\mathbb L}\) the ideal
	\(J\bar J\) is generated by the relative norm of \(J\), and the extension
	of ideals from \(\mathbb L_{+}\) to \(\mathbb L\) is injective. So
	\(J\bar J=s\mathcal O_{\mathbb L}\) with \(s\) totally positive means that
	the norm of \(J\) is \((s)\), and the image of the homomorphism
	\(\mathcal U(\mathcal O_{\mathbb L})\to\operatorname{Cl}(\mathbb L)\),
	\((J,s)\mapsto[J]\), is the kernel of \(\bar{\operatorname N}\). The kernel
	of this homomorphism consists of the classes of the pairs
	\((\mathcal O_{\mathbb L},\eta)\) with \(\eta\in U_{++}\), and two of them
	are equivalent if and only if they differ by an element of
	\(\operatorname N(U_{\mathbb L})\). Therefore, with
	Lemma~\ref{lem:classno}(a),(c),
	\[
	\#\mathcal U(\mathcal O_{\mathbb L})
	=[U_{++}:\operatorname N(U_{\mathbb L})]\cdot\#\ker\bar{\operatorname N}
	=\frac{[U_{++}:U_{+}^{2}]}{q_{\mathbb L}}\cdot
	\frac{r\,h_{\mathbb L}}{\#\operatorname{Cl}^{+}(\mathbb L_{+})} .
	\]
	The kernel of
	\(\operatorname{Cl}^{+}(\mathbb L_{+})\to\operatorname{Cl}(\mathbb L_{+})\)
	is the group of principal ideals modulo those with a totally positive
	generator, which is the quotient of the group \(\{\pm1\}^{2}\) of signs by
	the signs of the units; so
	\(\#\operatorname{Cl}^{+}(\mathbb L_{+})
	=4h_{\mathbb L_{+}}/[U_{+}:U_{++}]\). As
	\([U_{+}:U_{++}]\,[U_{++}:U_{+}^{2}]=[U_{+}:U_{+}^{2}]=4\), we get
	\(\#\mathcal U(\mathcal O_{\mathbb L})
	=r\,h_{\mathbb L}/(q_{\mathbb L}h_{\mathbb L_{+}})\), which is
	\(r\,h(\Delta)/2\) by Lemma~\ref{lem:classno}(b).

	(c) \(zI\bar I=\mathfrak c_{+}\mathcal O_{\mathbb L}\) means that \(z\)
	times the norm of \(I\) is \(\mathfrak c_{+}\); such \(I\) and totally
	positive \(z\) exist if and only if the narrow class of \(\mathfrak c_{+}\)
	is that of a norm. The last assertion is Lemma~\ref{lem:classno}(c).

	(d) The functor is well defined, because \(\mathcal O_{\mathbb L}\) is
	stable under the conjugation, and it is the identity on morphisms. Let
	\(x=(J,z)\) be an object of the target. A pair formed by an object
	\((I',z')\) of the source and a morphism \(a\) from
	\((I'\mathcal O_{\mathbb L},z')\) to \(x\) is equivalent to the pair formed
	by \((aI',z)\) and the identity, and two pairs of the latter form are
	equivalent only if they are equal. So \(\Fib(x)\) is the set of the
	invertible ideals \(I\) of \(R_{\mathrm{sp}}\) with
	\(I\mathcal O_{\mathbb L}=J\) and \(zI\bar I=\mathfrak c\). Such an \(I\)
	is determined by its completions. At the primes \(\ell\ne2\) the orders
	\(R_{\mathrm{sp}}\) and \(\mathcal O_{\mathbb L}\) have the same
	completion, the completion of \(I\) is that of \(J\), and the equation
	holds. Write a subscript \(2\) for the tensor product with \(\Z_2\).
	Invertible modules over the semilocal rings \(\mathcal O_{\mathbb L,2}\)
	and \(R_{\mathrm{sp},2}\) are free; let \(J_2=x_0\mathcal O_{\mathbb L,2}\)
	and \(\mathfrak c_2=c_0R_{\mathrm{sp},2}\). Then
	\(\eta=z\,x_0\bar x_0/c_0\) is a unit of \(\mathcal O_{\mathbb L,2}\), and
	\(\eta/\bar\eta=\bar c_0/c_0\) is a unit of \(R_{\mathrm{sp},2}\) because
	\(\bar{\mathfrak c}=\mathfrak c\). The invertible
	\(R_{\mathrm{sp},2}\)-modules that generate \(J_2\) are the
	\(vx_0R_{\mathrm{sp},2}\), with \(v\) a unit of
	\(\mathcal O_{\mathbb L,2}\) modulo the units of \(R_{\mathrm{sp},2}\), and
	the equation reads \(v\bar v\,\eta\in R_{\mathrm{sp},2}\). As
	\(R_{\mathrm{sp}}\) contains \(2\mathcal O_{\mathbb L}\), this depends only
	on the images of \(v\) and \(\eta\) in
	\(\bar{\mathcal O}=\mathcal O_{\mathbb L}/2\mathcal O_{\mathbb L}\), and the
	classes of \(v\) correspond to the units of \(\bar{\mathcal O}\) modulo
	the units of the image of \(R_{\mathrm{sp}}\).

	By Lemma~\ref{lem:orders}(a), \(\bar{\mathcal O}=\F_4[\delta]\) with
	\(\F_4=\Z[\omega]/2\) and \(\delta^{2}+\delta+m=0\), where
	\(m=(\eps p-1)/4\bmod2\), and the image of \(R_{\mathrm{sp}}\) is \(\F_4\).
	The conjugation acts on \(\bar{\mathcal O}\) by a ring automorphism which
	is \(x\mapsto x^{2}\) on \(\F_4\) and sends \(\delta\) to \(\delta\) if
	\(\eps=1\) and to \(1+\delta\) if \(\eps=-1\). The polynomial
	\(X^{2}+X+m\) has two roots \(t_1,t_2\) in \(\F_4\), with \(t_1+t_2=1\):
	they are \(0,1\) if \(m=0\), and they are the elements of
	\(\F_4\smallsetminus\F_2\), exchanged by \(x\mapsto x^{2}\), if \(m=1\).
	The homomorphisms of \(\F_4\)-algebras
	\(\operatorname{pr}_i\colon\bar{\mathcal O}\to\F_4\) with
	\(\operatorname{pr}_i(\delta)=t_i\) identify \(\bar{\mathcal O}\) with
	\(\F_4\times\F_4\) and the image of \(R_{\mathrm{sp}}\) with the diagonal.
	The map \(y\mapsto\operatorname{pr}_i(\bar y)^{2}\) is a homomorphism of
	\(\F_4\)-algebras from \(\bar{\mathcal O}\) to \(\F_4\), hence one of
	\(\operatorname{pr}_1\), \(\operatorname{pr}_2\). So
	\(\operatorname{pr}_i(\bar y)=\operatorname{pr}_j(y)^{2}\) for all \(y\),
	where \(j=i\) if \(\operatorname{pr}_i(\bar\delta)=t_i^{2}\) and \(j=3-i\)
	if \(\operatorname{pr}_i(\bar\delta)=t_{3-i}^{2}\). For \(\eps=1\),
	\(\operatorname{pr}_i(\bar\delta)=t_i\), which is \(t_i^{2}\) if \(m=0\)
	and \(t_{3-i}^{2}\) if \(m=1\); for \(\eps=-1\),
	\(\operatorname{pr}_i(\bar\delta)=1+t_i=t_{3-i}\), which is \(t_{3-i}^{2}\)
	if \(m=0\) and \(t_i^{2}\) if \(m=1\). So \(j=i\) exactly when
	\(\eps=1\) and \(p\equiv1\pmod8\), or \(\eps=-1\) and \(p\equiv3\pmod8\),
	that is, when \(\Delta\equiv5\pmod8\) (recall that \(\Delta=-3p\), resp.\
	\(-p\)); this is the case \(n=3\). Write
	\(v_i=\operatorname{pr}_i(v)\) and \(\eta_i=\operatorname{pr}_i(\eta)\);
	these are elements of \(\F_4^{\times}\), a group of order three.

	If \(j=i\), then \(\operatorname{pr}_i(v\bar v)=v_i^{3}=1\), so
	\(v\bar v=1\) in \(\bar{\mathcal O}\); and \(\eta/\bar\eta\), whose
	components are \(\eta_i^{-1}\), lies on the diagonal, so \(\eta\) does.
	The equation holds for the three classes of \(v\). If \(j=3-i\), then
	\(\operatorname{pr}_1(v\bar v)=v_1v_2^{2}\) and
	\(\operatorname{pr}_2(v\bar v)=v_2v_1^{2}=(v_1v_2^{2})^{2}\). The element
	\(t=v_1v_2^{2}\) depends only on the class of \(v\) modulo the diagonal,
	and every element of \(\F_4^{\times}\) is obtained from exactly one class.
	The equation says that \((t\eta_1,t^{2}\eta_2)\) lies on the diagonal, that
	is, \(t=\eta_1/\eta_2\); it holds for one class. The last assertion
	follows from Lemma~\ref{lem:fibre} and (a), all automorphism groups having
	order \(w\).
\end{proof}

\subsection*{Abelian surfaces and polarisations}

The element \(\pi\) is a Weil \(p\)-number without real conjugate. By the
theorem of Honda and Tate there is a simple abelian variety over \(\Fp\)
whose Frobenius endomorphism is a root of \(P\), and it is unique up to
isogeny; over the prime field its endomorphism algebra is commutative,
because \(\pi\) is not real
(\cite[Ch.~2, Theorem (Honda), p.~527, and Thm.~6.1(1)]{Waterhouse1969},
\cite[Thm.~1(c)]{Tate1966}; see
\cite[\S\S1.2 and~1.3]{CentelegheStix2015}). Hence this algebra is generated
by the Frobenius endomorphism and has degree twice the dimension over \(\Q\):
the variety is a surface with Weil polynomial \(P\), the category
\(\mathcal V_P\) is its isogeny class, and for every \(B\in\mathcal V_P\)
there is a unique isomorphism \(i_B\colon\mathbb L\to\End^{0}(B)\) with
\(i_B(\pi)=\pi_B\). A homomorphism \(f\colon B\to B'\) in \(\mathcal V_P\)
commutes with the Frobenius endomorphisms, so
\(f\,i_B(x)=i_{B'}(x)\,f\) for \(x\in\mathbb L\); we write \(fx=xf\) for this
element. Thus \(\Hom^{0}(B,B')\) is a vector space over \(\mathbb L\), of
dimension one.

\begin{lemma}\label{lem:cs}
	There is \(A_{\min}\in\mathcal V_P\) with
	\(\End(A_{\min})=i_{A_{\min}}(R_{\min})\). For \(B\in\mathcal V_P\) choose a
	non-zero \(g_B\in\Hom^{0}(B,A_{\min})\) and let \(\Lambda_B\) be the set of
	the \(x\in\mathbb L\) with \(xg_B\in\Hom(B,A_{\min})\); for
	\(f\in\Hom^{0}(B,B')\) define \(y_f\in\mathbb L\) by
	\(g_{B'}f=y_fg_B\).
	\begin{enumerate}[label=\textup{(\alph*)},leftmargin=*]
		\item \(\Lambda_B\) is a fractional ideal of \(R_{\min}\). The map
		\(f\mapsto y_f\) is a bijection from \(\Hom(B,B')\) onto
		\((\Lambda_B:\Lambda_{B'})\), with \(y_{f'f}=y_{f'}y_f\) and
		\(y_{i_B(x)}=x\). In particular \(\End(B)=i_B((\Lambda_B:\Lambda_B))\),
		and \(f\in\Hom^{0}(B,B')\) is an isomorphism if and only if
		\(y_f\Lambda_{B'}=\Lambda_B\).
		\item Every fractional ideal of \(R_{\min}\) is equal to \(\Lambda_B\)
		for some \(B\in\mathcal V_P\) and a suitable choice of \(g_B\).
		\item Let \(\mathcal R\) be \(R_{\mathrm{sp}}\) or
		\(\mathcal O_{\mathbb L}\). Then \(B\in\mathcal V_P[\mathcal R]\) if and
		only if \(\Lambda_B\) is an invertible ideal of \(\mathcal R\). For
		\(B,B',B''\in\mathcal V_P[\mathcal R]\), the group \(\Hom(B,B')\) is an
		invertible \(\mathcal R\)-module; every element of \(\Hom(B,B'')\) is a
		finite sum of composites \(f'f\) with \(f\in\Hom(B,B')\) and
		\(f'\in\Hom(B',B'')\); and \(f\in\Hom(B,B')\) is an isomorphism if and
		only if \(\Hom(B,B')=\{fx:x\in\mathcal R\}\).
		\item Let \(A_0\in\mathcal V_P[\mathcal R]\) and let \(I\) be an
		invertible fractional ideal of \(\mathcal R\). There are
		\(B\in\mathcal V_P[\mathcal R]\) and an invertible element \(f\) of
		\(\Hom^{0}(A_0,B)\) such that \(\Hom(A_0,B)=\{fx:x\in I\}\).
	\end{enumerate}
\end{lemma}

\begin{proof}
	By \cite[Prop.~21 and Thm.~25]{CentelegheStix2015}, applied to the set
	consisting of the conjugacy class of \(\pi\), there is \(A_{\min}\) as
	stated, and \(B\mapsto\Hom(B,A_{\min})\) is an anti-equivalence from the
	category of the abelian varieties over \(\Fp\) isogenous to a power of
	\(A_{\min}\) to the category of the finitely generated
	\(R_{\min}\)-modules that are free over \(\Z\)
	\cite[Thm.~11 and Lemma~13]{CentelegheStix2015}; here \(R_{\min}\) acts
	through \(\End(A_{\min})\), and the rank of \(\Hom(B,A_{\min})\) over
	\(\Z\) is \(2\dim B\). The objects of \(\mathcal V_P\) are those of
	dimension two.

	(a) For \(B\in\mathcal V_P\), \(\Hom(B,A_{\min})=\Lambda_Bg_B\) is a
	lattice in the \(\mathbb L\)-line \(\Hom^{0}(B,A_{\min})\), stable under
	\(R_{\min}\). The equivalence says that \(f\mapsto(g\mapsto gf)\) is a
	bijection from \(\Hom(B,B')\) onto the group of the \(R_{\min}\)-linear
	maps \(\Hom(B',A_{\min})\to\Hom(B,A_{\min})\). For
	\(f\in\Hom^{0}(B,B')\) and \(x\in\Lambda_{B'}\) we have
	\((xg_{B'})f=xy_fg_B\); so \(g\mapsto gf\) maps \(\Hom(B',A_{\min})\) into
	\(\Hom(B,A_{\min})\) if and only if
	\(y_f\Lambda_{B'}\subset\Lambda_B\), and then it is induced by some
	\(f_1\in\Hom(B,B')\), with \(g(f-f_1)=0\) for all \(g\), hence \(f=f_1\).
	The multiplicativity is clear, \(g_Bi_B(x)=xg_B\), and an equivalence
	reflects isomorphisms.

	(b) This is the essential surjectivity; an abelian variety \(B\) for which
	\(\Hom(B,A_{\min})\) has rank four is a surface, and a change of \(g_B\)
	multiplies \(\Lambda_B\) by an element of \(\mathbb L^{\times}\).

	(c) By (a), \(B\in\mathcal V_P[\mathcal R]\) if and only if
	\((\Lambda_B:\Lambda_B)=\mathcal R\). Let \(X\) be a lattice with
	\((X:X)=\mathcal R\). For lattices one has \((YZ)^{*}=(Y^{*}:Z)\) and
	\(Y^{**}=Y\); hence \((XX^{*})^{*}=(X:X)=\mathcal R\) and
	\(XX^{*}=\mathcal R^{*}\), which is invertible by
	Lemma~\ref{lem:orders}(b). So \(X\) is invertible. If \(\Lambda_B\) and
	\(\Lambda_{B'}\) are invertible ideals of \(\mathcal R\), then
	\((\Lambda_B:\Lambda_{B'})=\Lambda_B\Lambda_{B'}^{-1}\), the composition
	corresponds to the multiplication, and
	\(\Lambda_B\Lambda_{B'}^{-1}\cdot\Lambda_{B'}\Lambda_{B''}^{-1}
	=\Lambda_B\Lambda_{B''}^{-1}\). Finally \(y_f\Lambda_{B'}=\Lambda_B\) if and
	only if \(y_f\) generates \(\Lambda_B\Lambda_{B'}^{-1}\).

	(d) By (b) there is \(B\) with \(\Lambda_B=\Lambda_{A_0}I^{-1}\); then
	\((\Lambda_{A_0}:\Lambda_B)=I\), and we take for \(f\) the element with
	\(y_f=1\).
\end{proof}

\begin{proof}[Proof of Theorem~\textup{\ref{thm:vertexcm}(i)}]
	If \(A\) is superspecial, then \(\pi_A^{2}=u_A[p]\) with
	\(u_A\in\End(A)\) by Lemma~\ref{lem:frob2}, and \(u_A=i_A(u)\). Conversely
	let \(i_A(u)\in\End(A)\); it is an automorphism, because \(u\) is a root
	of unity. On the Dieudonn\'e module \(T_pA\) (see
	Appendix~\ref{sec:vertexpm}) the endomorphism \(\pi_A\) acts as one of the
	operators \(F\), \(V\), and the other one is
	\(p\pi_A^{-1}=u^{-1}\pi_A\). Hence \(F\,T_pA=V\,T_pA\), so \(A_{\kbar}\) has
	\(a\)-number two, and \(A\) is superspecial by the theorem of Oort
	\cite[Thm.~2 and Remark~3, p.~35]{Oort1975}, in the form quoted in
	\cite[\S3.2.1]{XueYu2020}, \cite[proof of Prop.~5.1, p.~1623]{XueYangYu2016}
	and \cite[Prop.~11.1]{HoweNartRitzenthaler2009}. By Lemma~\ref{lem:orders}(a)
	the orders containing \(R_{\mathrm{sp}}\) are \(R_{\mathrm{sp}}\) and
	\(\mathcal O_{\mathbb L}\).
\end{proof}

The category \(\mathcal V_P\) also contains surfaces that are not
superspecial, for instance \(A_{\min}\): the order \(R_{\min}\) has the
\(\Z\)-basis \(1,\pi,\bar\pi,pu\) and does not contain \(u\).

\begin{lemma}\label{lem:conjdual}
	Let \(B\in\mathcal V_P\) and \(x\in\mathbb L\). The dual of \(i_B(x)\) is
	\(i_{B^{\vee}}(\bar x)\). Consequently:
	\begin{enumerate}[label=\textup{(\alph*)},leftmargin=*]
		\item if \(B\in\mathcal V_P[\mathcal R]\) for an order
		\(\mathcal R=\bar{\mathcal R}\), then
		\(B^{\vee}\in\mathcal V_P[\mathcal R]\);
		\item for a polarisation \(\nu\) of \(B\), the dual of \(\nu x\) is
		\(\nu\bar x\), and \(E_\nu(xa,a')=E_\nu(a,\bar xa')\) for
		\(a,a'\in T_\ell B\otimes\Q\), \(\ell\ne p\).
	\end{enumerate}
\end{lemma}

\begin{proof}
	The maps \(x\mapsto\widehat{i_B(x)}\) and
	\(x\mapsto i_{B^{\vee}}(\bar x)\) are homomorphisms of \(\Q\)-algebras from
	\(\mathbb L\) to \(\End^{0}(B^{\vee})\), so it is enough to take \(x=\pi\).
	Every abelian variety over \(\Fp\) is projective
	\cite[I, Thm.~6.4]{MilneAV}, hence has a polarisation defined over
	\(\Fp\); let \(\nu\) be one for \(B\). Then \(\nu\pi_B=\pi_{B^{\vee}}\nu\),
	and \(\nu^{-1}\hat\pi_B\nu\,\pi_B=p\) by
	\cite[II, Lemma~1.2]{MilneAV}. Hence
	\(\hat\pi_B=\nu(p/\pi_B)\nu^{-1}=p/\pi_{B^{\vee}}\). (a) \(B^{\vee}\) is
	isogenous to \(B\), and \(\End(B^{\vee})\) consists of the duals of the
	elements of \(\End(B)\). (b) \(\widehat{\nu x}=\hat x\hat\nu
	=i_{B^{\vee}}(\bar x)\nu=\nu\bar x\), because \(\hat\nu=\nu\); and
	\(E_\nu(xa,a')=e_\ell(xa,\nu a')=e_\ell(a,\hat x\nu a')
	=E_\nu(a,\bar xa')\) by \cite[I, Prop.~13.2]{MilneAV}.
\end{proof}

\begin{lemma}\label{lem:positive}
	Let \(B\in\mathcal V_P\), let \(\nu\) be a polarisation of \(B\), and let
	\(s\in\mathbb L_{+}\) be non-zero with \(\nu s\in\Hom(B,B^{\vee})\). Then
	\(\nu s\) is a polarisation if and only if \(s\) is totally positive.
	Moreover, if \(g\colon B'\to B\) is an isogeny, then \(\hat g\nu g\) is a
	polarisation of \(B'\).
\end{lemma}

\begin{proof}
	We work over \(\kbar\) and write \(\nu=\phi_{\mathcal N}\) with
	\(\mathcal N\) ample. For the last assertion,
	\(\hat g\nu g=\phi_{g^{*}\mathcal N}\); if \(k\ge1\) is such that
	\(\mathcal N^{k}\) is very ample, \(g^{*}\mathcal N^{k}\) has a non-zero
	section and \(\phi_{g^{*}\mathcal N^{k}}=k\,\hat g\nu g\) is an isogeny,
	so \(g^{*}\mathcal N^{k}\) is ample by \cite[I, Prop.~8.1]{MilneAV}, and so
	is \(g^{*}\mathcal N\).

	By Lemma~\ref{lem:conjdual}(b) the form
	\(E_{\nu s}(a,a')=E_\nu(a,sa')\) is skew-symmetric, because \(E_\nu\) is
	and \(E_\nu(a',sa)=E_\nu(sa',a)\). Hence
	\(\nu s=\phi_{\mathcal M}\) for an invertible sheaf \(\mathcal M\)
	\cite[I, Prop.~13.6]{MilneAV}, and \(\nu s\) is a polarisation if and only
	if \(\mathcal M\) can be chosen ample. We compute intersection numbers.
	For \(x\in\mathbb L_{+}\) with \(\nu x\in\Hom(B,B^{\vee})\) we have, the
	degree of the zero homomorphism being \(0\),
	\[
	\deg(\nu x)=\deg(\nu)\,\operatorname N_{\mathbb L_{+}/\Q}(x)^{2}:
	\]
	if \(k\ge1\) is an integer with \(kx\in\End(B)\), then
	\(k^{4}\deg(\nu x)=\deg(\nu)\deg(kx)\), and \(\deg(kx)\) is the
	determinant of \(kx\) on \(T_\ell B\otimes\Q\)
	\cite[I, Prop.~10.20]{MilneAV}, a free module of rank one over
	\(\mathbb L\otimes\Q_\ell\), so it is
	\(\operatorname N_{\mathbb L/\Q}(kx)=k^{4}\operatorname
	N_{\mathbb L_{+}/\Q}(x)^{2}\). For integers \(a,b\) the sheaf
	\(\mathcal N^{a}\otimes\mathcal M^{b}\) defines the homomorphism
	\(\nu(a+bs)\). By the Riemann--Roch theorem
	\cite[\S16]{Mumford2008}, \cite[I, Thm.~11.1]{MilneAV}, an invertible sheaf
	\(\mathcal K\) on the surface \(B\) has Euler characteristic
	\(\chi(\mathcal K)=(\mathcal K\cdot\mathcal K)/2\), and
	\(\chi(\mathcal K)^{2}=\deg\phi_{\mathcal K}\). Therefore
	\[
	\bigl(a^{2}(\mathcal N\cdot\mathcal N)+2ab(\mathcal N\cdot\mathcal M)
	+b^{2}(\mathcal M\cdot\mathcal M)\bigr)^{2}
	=\bigl(2\chi(\mathcal N)\bigr)^{2}
	\bigl(a^{2}+ab\operatorname{Tr}(s)+b^{2}\operatorname N(s)\bigr)^{2}
	\]
	for all integers \(a,b\), where \(\operatorname{Tr}\) and
	\(\operatorname N\) are the trace and the norm of \(\mathbb L_{+}\) over
	\(\Q\). Two polynomials in \(a,b\) with equal squares differ by a sign,
	and here the sign is \(+\), because
	\((\mathcal N\cdot\mathcal N)=2\chi(\mathcal N)>0\). Hence
	\[
	(\mathcal N\cdot\mathcal M)=\chi(\mathcal N)\operatorname{Tr}(s),\qquad
	(\mathcal M\cdot\mathcal M)=2\chi(\mathcal N)\operatorname N(s).
	\]
	If \(\mathcal M\) is ample, both numbers are positive, so
	\(\operatorname{Tr}(s)>0\) and \(\operatorname N(s)>0\): \(s\) is totally
	positive. Conversely let \(s\) be totally positive. Then
	\(\chi(\mathcal M)=\chi(\mathcal N)\operatorname N(s)>0\), so
	\(H^{0}(\mathcal M)\ne0\) or \(H^{2}(\mathcal M)\ne0\). By Serre duality,
	the canonical sheaf of \(B\) being trivial, \(H^{2}(\mathcal M)\) is dual
	to \(H^{0}(\mathcal M^{-1})\); if this were non-zero,
	\(\mathcal M^{-1}\) would be the sheaf of an effective divisor and
	\((\mathcal N\cdot\mathcal M^{-1})\ge0\), contradicting
	\((\mathcal N\cdot\mathcal M)=\chi(\mathcal N)\operatorname{Tr}(s)>0\).
	So \(\mathcal M\) has a non-zero section; and
	\(\phi_{\mathcal M}\) is an isogeny because
	\(\deg\phi_{\mathcal M}=\chi(\mathcal M)^{2}\ne0\). Hence \(\mathcal M\) is
	ample \cite[I, Prop.~8.1]{MilneAV}.
\end{proof}

In Appendix~\ref{sec:vertexpm} the corresponding statement is quoted from
\cite[\S21, Application~III, p.~208]{Mumford2008}, where a symmetric endomorphism is called totally
positive if the roots of its characteristic polynomial are positive. For
\(s\in\mathbb L_{+}\) these roots are the conjugates of \(s\), so the two
notions agree.

\begin{proposition}\label{prop:dict}
	Let \(\mathcal R\) be \(R_{\mathrm{sp}}\) or \(\mathcal O_{\mathbb L}\), let
	\(A_0\in\mathcal V_P[\mathcal R]\), and let \(\theta_0\) be a polarisation
	of \(A_0\), not necessarily principal. Put
	\(\mathfrak c=\{x\in\mathbb L:\theta_0x\in\Hom(A_0,A_0^{\vee})\}\).
	\begin{enumerate}[label=\textup{(\alph*)},leftmargin=*]
		\item \(\mathfrak c\) is an invertible fractional ideal of
		\(\mathcal R\), and \(\bar{\mathfrak c}=\mathfrak c\).
		\item Let \(B\in\mathcal V_P[\mathcal R]\), and let \(f\) be an
		invertible element of \(\Hom^{0}(A_0,B)\) and \(I\) the ideal with
		\(\Hom(A_0,B)=\{fx:x\in I\}\). Then
		\[
		\Hom(B,B^{\vee})=\bigl\{\theta_z=\hat f^{-1}\theta_0zf^{-1}:
		z\in\mathfrak c\,(I\bar I)^{-1}\bigr\},
		\]
		and \(\theta_z\) is a principal polarisation if and only if
		\(z\in\mathbb L_{+}\), \(z\) is totally positive and
		\(zI\bar I=\mathfrak c\).
		\item The groupoid \(\mathcal A_P[\mathcal R]\) is equivalent to
		\(\mathcal I(\mathcal R,\mathfrak c)\).
		\item \(\mathfrak c\,\mathcal O_{\mathbb L}
		=\mathfrak c_{+}\mathcal O_{\mathbb L}\) for a fractional ideal
		\(\mathfrak c_{+}\) of \(\mathbb L_{+}\).
	\end{enumerate}
\end{proposition}

\begin{proof}
	(a) \(A_0^{\vee}\in\mathcal V_P[\mathcal R]\) by
	Lemma~\ref{lem:conjdual}(a), so \(\Hom(A_0,A_0^{\vee})=\theta_0\mathfrak c\)
	is an invertible \(\mathcal R\)-module (Lemma~\ref{lem:cs}(c)). It is
	stable under duality, and the dual of \(\theta_0x\) is \(\theta_0\bar x\)
	(Lemma~\ref{lem:conjdual}(b)).

	(b) The objects \(A_0\), \(A_0^{\vee}\), \(B\), \(B^{\vee}\) lie in
	\(\mathcal V_P[\mathcal R]\), and we apply Lemma~\ref{lem:cs}(c) to them.
	First, the composites of the elements of \(\Hom(B,A_0)\) with those of
	\(\Hom(A_0,B)\) generate \(\End(A_0)=\mathcal R\); hence
	\(\Hom(B,A_0)=\{xf^{-1}:x\in I^{-1}\}\). Taking duals, with
	Lemma~\ref{lem:conjdual},
	\(\Hom(A_0^{\vee},B^{\vee})=\{\hat f^{-1}x:x\in\bar I^{-1}\}\). Composing
	with \(\Hom(A_0,A_0^{\vee})=\theta_0\mathfrak c\) gives
	\(\Hom(A_0,B^{\vee})=\{\hat f^{-1}\theta_0x:x\in\bar I^{-1}\mathfrak c\}\),
	and composing with \(\Hom(B,A_0)\) gives the formula for
	\(\Hom(B,B^{\vee})\). The dual of \(\theta_z\) is \(\theta_{\bar z}\), so
	\(\theta_z\) is symmetric if and only if \(z\in\mathbb L_{+}\). It is an
	isomorphism if and only if it generates \(\Hom(B,B^{\vee})\), that is,
	\(z\mathcal R=\mathfrak c(I\bar I)^{-1}\).

	Let \(z\in\mathbb L_{+}\) be non-zero with
	\(\theta_z\in\Hom(B,B^{\vee})\); we show that \(\theta_z\) is a
	polarisation if and only if \(z\) is totally positive. Let \(\nu\) be a
	polarisation of \(B\). It is symmetric, so \(\nu=\theta_{z_0}\) with
	\(z_0\in\mathbb L_{+}\). Let \(k\ge1\) be an integer with
	\(kf\in\Hom(A_0,B)\). By Lemma~\ref{lem:positive} the pull-back
	\(\widehat{kf}\,\nu\,kf=\theta_0\,k^{2}z_0\) is a polarisation of \(A_0\),
	and then \(k^{2}z_0\) is totally positive, by the same lemma applied to
	\((A_0,\theta_0)\). Now \(\theta_z=\nu\,(z/z_0)\), and by
	Lemma~\ref{lem:positive} applied to \((B,\nu)\) this is a polarisation if
	and only if \(z/z_0\) is totally positive.

	(c) For each object \((B,\theta)\) of \(\mathcal A_P[\mathcal R]\) choose
	\(f_B\) and \(I_B\) as in (b) (Lemma~\ref{lem:cs}(c) shows that they
	exist), and let \(z\) be the element with \(\theta=\theta_z\); by (b),
	\((I_B,z)\) is an object of \(\mathcal I(\mathcal R,\mathfrak c)\). Let
	\(g\colon(B,\theta_z)\to(B',\theta_{z'})\) be an isomorphism, and define
	\(a\in\mathbb L^{\times}\) by \(gf_B=f_{B'}a\). Then
	\(\Hom(A_0,B')=g\Hom(A_0,B)\) gives \(I_{B'}=aI_B\), and
	\[
	\hat g\,\theta_{z'}\,g=\hat f_B^{-1}\,\hat a\,\theta_0z'\,a\,f_B^{-1}
	=\hat f_B^{-1}\theta_0\,a\bar az'\,f_B^{-1}
	\]
	gives \(z=a\bar az'\); so \(a\) is a morphism from \((I_B,z)\) to
	\((I_{B'},z')\). Conversely such a morphism \(a\) defines
	\(g=f_{B'}af_B^{-1}\), which generates
	\(\Hom(B,B')=\{f_{B'}xf_B^{-1}:x\in I_{B'}I_B^{-1}\}\), hence is an
	isomorphism (Lemma~\ref{lem:cs}(c)), and \(\hat g\theta_{z'}g=\theta_z\).
	So we have a fully faithful functor. It is essentially surjective: for an
	object \((I,z)\) of \(\mathcal I(\mathcal R,\mathfrak c)\),
	Lemma~\ref{lem:cs}(d) gives \(B\) and \(f\) with
	\(\Hom(A_0,B)=\{fx:x\in I\}\), and \(\theta_z\) is a principal polarisation
	of \(B\) by (b).

	(d) The ideal \(\mathfrak c\,\mathcal O_{\mathbb L}\) is stable under the
	conjugation. Such an ideal comes from \(\mathbb L_{+}\) as soon as its
	valuation is even at the primes of \(\mathbb L\) that are ramified over
	\(\mathbb L_{+}\): at a prime which is inert there is nothing to prove, and
	at two conjugate primes the valuations are equal. By
	Lemma~\ref{lem:orders}(d) there is no ramified prime for \(P_e\), and for
	\(P_c\), \(P_d\) the ramified primes lie above \(\ell=2\), resp.\
	\(\ell=3\). In both cases \(\ell\ne p\) and
	\(\mathcal R\otimes\Z_\ell=\mathcal O_{\mathbb L}\otimes\Z_\ell\), a
	product of discrete valuation rings which we denote by \(\mathcal O_\ell\);
	we put \(\mathbb L_\ell=\mathbb L\otimes\Q_\ell\).

	The algebra \(\mathbb L_\ell\) acts faithfully on
	\(T_\ell A_0\otimes\Q\) \cite[I, Thm.~10.15]{MilneAV}; it is a product of
	fields whose degrees over \(\Q_\ell\) have sum four, so
	\(T_\ell A_0\otimes\Q\) is free of rank one over \(\mathbb L_\ell\), and
	the lattice \(T_\ell A_0\), which is stable under \(\mathcal O_\ell\), is
	of the form \(\mathcal O_\ell t_0\). The form
	\((x,y)\mapsto E_{\theta_0}(xt_0,yt_0)\) on \(\mathbb L_\ell\) is
	alternating, non-degenerate, and the adjoint of the multiplication by
	\(a\) is the multiplication by \(\bar a\)
	(Lemma~\ref{lem:conjdual}(b)). Hence there is a unique
	\(\xi\in\mathbb L_\ell^{\times}\) with
	\[
	E_{\theta_0}(xt_0,yt_0)
	=\operatorname{Tr}_{\mathbb L_\ell/\Q_\ell}(\xi\,\bar xy),
	\]
	and \(\bar\xi=-\xi\), because the form is alternating. An element
	\(x\in\mathbb L_\ell\) lies in \(\mathfrak c\otimes\Z_\ell\) if and only
	if \(\theta_0x\) maps \(T_\ell A_0\) into \(T_\ell A_0^{\vee}\) (the map
	\(\Hom(A_0,A_0^{\vee})\otimes\Z_\ell\to\Hom(T_\ell A_0,T_\ell A_0^{\vee})\)
	is injective with torsion-free cokernel
	\cite[I, Thm.~10.15]{MilneAV}). As \(e_\ell\) is a perfect pairing between
	\(T_\ell A_0\) and \(T_\ell A_0^{\vee}\), this means that
	\(\operatorname{Tr}(\xi\,\bar a\,xb)\in\Z_\ell\) for all
	\(a,b\in\mathcal O_\ell\), that is, \(\xi x\in\mathcal O_\ell^{*}\). So
	\(\mathfrak c\otimes\Z_\ell=\xi^{-1}\mathcal O_\ell^{*}\), where
	\(\mathcal O_\ell^{*}\) is the dual of \(\mathcal O_\ell\) for the trace
	form of \(\mathbb L_\ell\) over \(\Q_\ell\); it is
	\(g^{-1}\mathcal O_\ell\) with \(g\) as in Lemma~\ref{lem:orders}(b).

	For \(P_c\) and \(\ell=2\): \(g=4\pi\) and \(\pi\) is a unit at \(2\);
	the unit \(u\) satisfies \(\bar u=-u\), so \(\xi=ut\) with \(t\) an
	invertible element of \(\mathbb L_{+}\otimes\Q_2\), and
	\(\mathfrak c\otimes\Z_2\) is generated by \((4t)^{-1}\). For \(P_d\) and
	\(\ell=3\): \(g=2\pi_0\sqrt{-3}\), where \(2\pi_0\) is a unit at \(3\);
	\(\xi=\sqrt{-3}\,t\) with \(t\) an invertible element of
	\(\mathbb L_{+}\otimes\Q_3\), and \(\mathfrak c\otimes\Z_3\) is generated
	by \((3t)^{-1}\). In both cases \(\mathfrak c\otimes\Z_\ell\) is generated
	by an element of \(\mathbb L_{+}\otimes\Q_\ell\), and its valuations at
	the ramified primes are even.
\end{proof}

\begin{lemma}\label{lem:descend}
	Let \(R_{\mathrm{sp}}\ne\mathcal O_{\mathbb L}\). If
	\(\mathcal A_P[R_{\mathrm{sp}}]\) is non-empty, then
	\(\mathcal A_P[\mathcal O_{\mathbb L}]\) is non-empty.
\end{lemma}

\begin{proof}
	Let \((A_1,\theta_1)\in\mathcal A_P[R_{\mathrm{sp}}]\). By
	Lemma~\ref{lem:cs} there is \(A_0\in\mathcal V_P\) with
	\(\Lambda_{A_0}=\Lambda_{A_1}\mathcal O_{\mathbb L}\); then
	\(\End(A_0)=\mathcal O_{\mathbb L}\),
	\((\Lambda_{A_0}:\Lambda_{A_1})=\mathcal O_{\mathbb L}\) and
	\((\Lambda_{A_1}:\Lambda_{A_0})
	=(R_{\mathrm{sp}}:\mathcal O_{\mathbb L})=2\mathcal O_{\mathbb L}\): in
	\(\mathcal O_{\mathbb L}/2\mathcal O_{\mathbb L}=\F_4\times\F_4\) the image
	of \(R_{\mathrm{sp}}\) is the diagonal (proof of
	Lemma~\ref{lem:ideals}(d)), which contains no non-zero ideal. Let
	\(\rho\in\Hom(A_0,A_1)\) and
	\(\rho'\in\Hom(A_1,A_0)\) be the elements with \(y_\rho=1\) and
	\(y_{\rho'}=2\); then \(\rho\rho'=2\), \(\rho'\rho=2\), and
	\(\Hom(A_1,A_0)=\{x\rho':x\in\mathcal O_{\mathbb L}\}\). Taking duals and
	composing with the isomorphism \(\theta_1^{-1}\), we see that
	\(\Hom(A_0^{\vee},A_1)=\{\gamma x:x\in\mathcal O_{\mathbb L}\}\) with
	\(\gamma=\theta_1^{-1}\hat\rho'\). As
	\(A_0^{\vee}\in\mathcal V_P[\mathcal O_{\mathbb L}]\)
	(Lemma~\ref{lem:conjdual}), this says that
	\((\Lambda_{A_0^{\vee}}:\Lambda_{A_1})
	=\Lambda_{A_0^{\vee}}(\Lambda_{A_1}\mathcal O_{\mathbb L})^{-1}\) is
	generated by \(y_\gamma\), that is,
	\(\Lambda_{A_0^{\vee}}=y_\gamma\Lambda_{A_0}\). By Lemma~\ref{lem:cs}(a) an
	element \(\varphi\) of \(\Hom^{0}(A_0,A_0^{\vee})\) is therefore an
	isomorphism if and only if \(y_\gamma y_{\varphi}\) is a unit of
	\(\mathcal O_{\mathbb L}\). Put
	\(\varphi_0=\tfrac12\hat\rho\,\theta_1\rho\). Then
	\[
	\gamma\varphi_0=\tfrac12\,\theta_1^{-1}\,\hat\rho'\hat\rho\,\theta_1\rho
	=\tfrac12\,\theta_1^{-1}\,\widehat{\rho\rho'}\,\theta_1\rho=\rho ,
	\]
	so \(y_\gamma y_{\varphi_0}=1\) and \(\varphi_0\) is an isomorphism. By
	Lemma~\ref{lem:positive}, \(\hat\rho\,\theta_1\rho\) is a polarisation of
	\(A_0\), and so is its product \(\varphi_0\) with the totally positive
	element \(1/2\). Hence
	\((A_0,\varphi_0)\in\mathcal A_P[\mathcal O_{\mathbb L}]\).
\end{proof}

\begin{lemma}\label{lem:weilres}
	The objects of \(\mathcal A_P\) that are Weil restrictions of elliptic
	curves over \(\F_{p^{2}}\) are superspecial and have endomorphism ring
	\(R_{\mathrm{sp}}\). Their mass is \(1/4\) for \(P_c\) if
	\(p\equiv3\pmod4\), and \(1/6\) for \(P_d\) and for \(P_e\) if
	\(p\equiv2\pmod3\); there are none for \(P_c\) if \(p\equiv1\pmod4\), and
	none for \(P_d\), \(P_e\) if \(p\equiv1\pmod3\).
\end{lemma}

\begin{proof}
	Let \(A\) be the Weil restriction of an elliptic curve \(B/\F_{p^{2}}\),
	with Weil polynomial \(P\). Over \(\kbar\) it is the product of \(B\) and
	its conjugate \(B'\), and \(\pi_A^{2}\) induces on \(B\) its Frobenius
	endomorphism \(\pi_B\), which satisfies \(\pi_B^{w}=p^{w}\); so \(B\) is
	supersingular and \(A\) is superspecial. By the proof of
	Proposition~\ref{prop:nonjac}, \(\pi_B=pu_B\) with \(u_B\) an automorphism
	of \(B\), of order four for \(P_c\) and of order three or six for \(P_d\),
	\(P_e\); such supersingular curves exist only if \(p\equiv3\pmod4\),
	resp.\ \(p\equiv2\pmod3\), and the masses are those of
	Proposition~\ref{prop:nonjac}. By Theorem~\ref{thm:vertexcm}(i),
	\(\End(A)\) contains \(R_{\mathrm{sp}}\). Let
	\(R_{\mathrm{sp}}\ne\mathcal O_{\mathbb L}\). The endomorphism \(\pi_A\)
	maps \(B\) to \(B'\), and \(i_A(u)\) maps \(B\) to itself, because
	\(p\,i_A(u)=\pi_A^{2}\) does. Hence \(\pi_0=\eps u\pi\) maps \(B[2]\) into
	\(B'[2]\), and it is not the identity on \(A[2]=B[2]\times B'[2]\). If
	\(\delta=(1+\pi_0)/2\) were an endomorphism of \(A\), then
	\(\pi_0=2\delta-1\) would be the identity on \(A[2]\). So
	\(\End(A)=R_{\mathrm{sp}}\) by Lemma~\ref{lem:orders}(a).
\end{proof}

\begin{proof}[Proof of Theorem~\textup{\ref{thm:vertexcm}(ii)}]
	Let \(\mathcal R\) be \(R_{\mathrm{sp}}\) or \(\mathcal O_{\mathbb L}\).
	The category \(\mathcal V_P[\mathcal R]\) is not empty
	(Lemma~\ref{lem:cs}(b),(c), applied to the ideal \(\mathcal R\)), and each
	of its objects has a polarisation over \(\Fp\), as in the proof of
	Lemma~\ref{lem:conjdual}. Fix an object \(A_0\) with a polarisation
	\(\theta_0\), and let \(\mathfrak c\) be the ideal of
	Proposition~\ref{prop:dict}. By Proposition~\ref{prop:dict}(c) and
	Lemma~\ref{lem:ideals}(a), every object of \(\mathcal A_P[\mathcal R]\) has
	automorphism group \(\langle u,-1\rangle\), the number of classes of
	\(\mathcal A_P[\mathcal R]\) is \(0\) or \(\#\mathcal U(\mathcal R)\), and
	the mass is the number of classes divided by \(w\).

	Let \(\mathcal R=\mathcal O_{\mathbb L}\). By
	Proposition~\ref{prop:dict}(d),
	\(\mathfrak c=\mathfrak c_{+}\mathcal O_{\mathbb L}\). For \(P_c\) and
	\(P_d\) the groupoid
	\(\mathcal I(\mathcal O_{\mathbb L},\mathfrak c)\) is not empty by
	Lemma~\ref{lem:ideals}(c), so the number of classes is
	\(\#\mathcal U(\mathcal O_{\mathbb L})=h(\Delta)/2\)
	(Lemma~\ref{lem:ideals}(b)). For \(P_e\) and \(p\equiv1\pmod3\) the
	isogeny class contains no principally polarised surface
	\cite[Thm.~2]{HoweMaisnerNartRitzenthaler2008}, so
	\(\mathcal A_P[\mathcal O_{\mathbb L}]\) and
	\(\mathcal A_P[R_{\mathrm{sp}}]\) are empty. For \(P_e\) and
	\(p\equiv2\pmod3\) there is a Weil restriction in
	\(\mathcal A_P[R_{\mathrm{sp}}]\) (Lemma~\ref{lem:weilres}); hence
	\(\mathcal A_P[\mathcal O_{\mathbb L}]\) is not empty, trivially if
	\(R_{\mathrm{sp}}=\mathcal O_{\mathbb L}\) and by
	Lemma~\ref{lem:descend} otherwise, and its number of classes is
	\(\#\mathcal U(\mathcal O_{\mathbb L})=h(\Delta)\).

	Now let \(\mathcal R=R_{\mathrm{sp}}\ne\mathcal O_{\mathbb L}\). By
	Lemma~\ref{lem:ideals}(d) the number of classes of
	\(\mathcal I(R_{\mathrm{sp}},\mathfrak c)\) is \(n\) times that of
	\(\mathcal I(\mathcal O_{\mathbb L},\mathfrak c\,\mathcal O_{\mathbb L})\),
	which is \(0\) or \(\#\mathcal U(\mathcal O_{\mathbb L})\). For \(P_d\) the
	ideal \(\mathfrak c\,\mathcal O_{\mathbb L}\) comes from
	\(\mathbb L_{+}\) by Proposition~\ref{prop:dict}(d), applied to
	\(\mathcal R=R_{\mathrm{sp}}\), so
	\(\mathcal I(\mathcal O_{\mathbb L},\mathfrak c\,\mathcal O_{\mathbb L})\)
	is not empty by Lemma~\ref{lem:ideals}(c), and the number of classes is
	\(n\,h(\Delta)/2\). For \(P_e\) and \(p\equiv2\pmod3\) the groupoid
	\(\mathcal A_P[R_{\mathrm{sp}}]\) contains a Weil restriction, so both
	groupoids are non-empty and the number of classes is \(n\,h(\Delta)\). In
	both cases \(\Delta\) is odd and \(n\,h(\Delta)=h(4\Delta)\) by
	Lemma~\ref{lem:classno}(d).
\end{proof}

\begin{proof}[Proof of Theorem~\textup{\ref{thm:vertexcm}(iii)}]
	Let \((A,\theta)\) be an object of \(\mathcal A_P[R_{\mathrm{sp}}]\) or of
	\(\mathcal A_P[\mathcal O_{\mathbb L}]\); it is superspecial by (i). A
	principally polarised abelian surface over \(\Fp\) is the Jacobian of a
	curve of genus two over \(\Fp\), or a product of two elliptic curves over
	\(\Fp\) with the product polarisation, or the Weil restriction of an
	elliptic curve over \(\F_{p^{2}}\), and these cases exclude each other
	\cite[Thm.~1.3]{HoweNartRitzenthaler2009}. The Weil polynomial of a
	product is reducible, so the second case does not occur. In the first
	case the curve \(C\) is unique up to \(\Fp\)-isomorphism and
	\(\Aut_{\Fp}(C)\cong\Aut_{\Fp}(A,\theta)\) (Lemmas~\ref{lem:dual}
	and~\ref{lem:torelli}), and \(C\) is superspecial. The assertion on the
	Weil restrictions is Lemma~\ref{lem:weilres}.

	Let \(P\) be \(P_d\) or \(P_e\) and \(J=\Jac(C)\). By
	Theorem~\ref{thm:class} the Frobenius type is \(3^{2}\) or \(6\). The
	endomorphism \(\pi\) acts on \(J[2]\) through the permutation \(\sigma\)
	of the Weierstrass points and \(u=\pi^{2}/p\) through \(\sigma^{2}\), so
	\(\pi_0=\eps u\pi\) acts through \(\sigma^{3}\)
	(Lemma~\ref{lem:J2}(d)). Hence the type is \(3^{2}\) if and only if
	\(\pi_0\) is the identity on \(J[2]\), that is, if and only if
	\(\pi_0-1\) is twice an endomorphism \cite[I, Lemma~10.16]{MilneAV}. This
	means that \((\pi_0-1)/2\) lies in \(\End(J)\); by
	Lemma~\ref{lem:orders}(a) it does if and only if
	\(\End(J)=\mathcal O_{\mathbb L}\ne R_{\mathrm{sp}}\).
\end{proof}

\begin{proof}[Proof of Theorem~\textup{\ref{thm:vertexcmcells}}]
	By Theorem~\ref{thm:vertexcm}(i),(iii), the mass of the superspecial
	curves with Weil polynomial \(P\) and \(\End(\Jac C)=\mathcal R\) is the
	mass of \(\mathcal A_P[\mathcal R]\) minus the mass of the Weil
	restrictions in it, which are all in \(\mathcal A_P[R_{\mathrm{sp}}]\) and
	whose mass is given by Lemma~\ref{lem:weilres}. The masses of the
	groupoids \(\mathcal A_P[\mathcal R]\) are given by
	Theorem~\ref{thm:vertexcm}(ii), with \(\Delta\) as
	in~\eqref{eq:cmfields}.

	For part (c), \(R_{\mathrm{sp}}=\mathcal O_{\mathbb L}\),
	\(w=4\), and the mass of \(\mathcal A_P[\mathcal O_{\mathbb L}]\) is
	\(h(-8p)/8\). By Theorem~\ref{thm:class} the Frobenius type is \(2\,4\) if
	\(p\equiv3\pmod4\) and \(1^{2}4\) if \(p\equiv1\pmod4\). In the first case
	we subtract \(1/4\) and find \((h(-8p)-2)/8\); in the second there is no
	Weil restriction and the mass is \(h(-8p)/8\).

	For part (d), \(w=6\). If \(p\equiv3\pmod4\), then
	\(R_{\mathrm{sp}}=\mathcal O_{\mathbb L}\), all the curves have type
	\(6\), and the mass of \(\mathcal A_P[\mathcal O_{\mathbb L}]\) is
	\(h(-12p)/12\). If \(p\equiv1\pmod4\), the curves of type \(3^{2}\) are
	those with \(\End(\Jac C)=\mathcal O_{\mathbb L}\); there is no Weil
	restriction in \(\mathcal A_P[\mathcal O_{\mathbb L}]\), and their mass is
	\(h(-3p)/12\). The curves of type \(6\) are those with
	\(\End(\Jac C)=R_{\mathrm{sp}}\), and the mass of
	\(\mathcal A_P[R_{\mathrm{sp}}]\) is \(h(-12p)/12\). In both cases the
	Weil restrictions have mass \(1/6\) if \(p\equiv2\pmod3\) and do not
	exist otherwise, so the curves of type \(6\) have mass
	\(\bigl(h(-12p)-2\,[p\equiv2\bmod3]\bigr)/12\).

	For part (e), \(w=6\), and for \(p\equiv1\pmod3\) there
	is nothing. Let \(p\equiv2\pmod3\). If \(p\equiv1\pmod4\), then
	\(R_{\mathrm{sp}}=\mathcal O_{\mathbb L}\), all the curves have type
	\(6\), the mass of \(\mathcal A_P[\mathcal O_{\mathbb L}]\) is
	\(h(-4p)/6\), and we subtract \(1/6\). If \(p\equiv3\pmod4\), that is,
	\(p\equiv11\pmod{12}\), the curves of type \(3^{2}\) have mass
	\(h(-p)/6\), the mass of \(\mathcal A_P[\mathcal O_{\mathbb L}]\), and
	those of type \(6\) have mass \(h(-4p)/6-1/6\).

	These are the masses stated in Theorem~\ref{thm:vertex}(c)--(e), and
	we have seen that the other pairs of a Frobenius type and a residue class
	are empty.
\end{proof}

\begin{remark}\label{rem:vertexcmlit}
	Apart from the first paragraph, which concerns results quoted above, the
	following comparisons with the literature are not used in the proofs.

	As said after
	Theorem~\ref{thm:vertexcm}, assertion (i) of that theorem and
	Lemma~\ref{lem:orders}(a) are in \cite[\S5]{XueYangYu2016}, where the
	numbers of isomorphism classes of the superspecial abelian surfaces
	themselves are computed: they are the class numbers of the orders
	\(R_{\mathrm{sp}}\) and \(\mathcal O_{\mathbb L}\). We did not find there,
	or elsewhere, the principally polarised surfaces with these Weil
	polynomials counted with a given endomorphism ring, as in
	Theorem~\ref{thm:vertexcm}(ii).

	For \(P_c\) and \(P_d\) the
	existence of a principally polarised surface in the isogeny class is
	stated in \cite[Thm.~5 and Table~1]{HoweMaisnerNartRitzenthaler2008}, with
	the same reason as here, the ramification of
	\(\mathbb L/\mathbb L_{+}\) above \(2\), resp.\ \(3\);
	Proposition~\ref{prop:dict}(d) and Lemma~\ref{lem:ideals}(c) give it for
	each of the two endomorphism rings. For \(P_e\) the obstruction of
	\cite[Thm.~6]{HoweMaisnerNartRitzenthaler2008} lies in the cokernel of
	\(\bar{\operatorname N}\), the group of Lemma~\ref{lem:ideals}(c). In
	Section~4 of that paper the obstruction is computed for
	\(p\equiv1\pmod3\), from a polarisation of degree nine on a twist of the
	square of an elliptic curve; that section assumes \(p\equiv1\pmod3\), and
	for \(p\equiv2\pmod3\) we use the Weil restriction instead.

	Summing the masses of Theorem~\ref{thm:vertexcm}(ii)
	over the three polynomials and over the orders gives
	\[
	\frac{h(-8p)}8+\frac{h(-12p)}{12}
	+[p\equiv2\bmod3]\,\frac{\bigl(3-\leg2p\bigr)h(-p)}6
	\]
	for \(p\equiv3\pmod4\), because \(h(-4p)=\bigl(2-\leg2p\bigr)h(-p)\), and
	\[
	\frac{h(-8p)}8+\frac{\bigl(3+\leg{-2}{p}\bigr)h(-3p)}{12}
	+[p\equiv2\bmod3]\,\frac{h(-4p)}6
	\]
	for \(p\equiv1\pmod4\). Written with class numbers of fields, these are
	the terms of the formula of Ibukiyama and Katsura
	\cite[Remark~3]{IbukiyamaKatsura1994} for the number of principally
	polarised superspecial surfaces over \(\kbar\) with a model over \(\Fp\)
	that remain when the terms equal to the masses for \(\pi^{2}=-p\) and
	\(\pi^{2}=p\) (Remark~\ref{rem:vertexlit}) are removed.
\end{remark}

\begin{remark}\label{rem:vertexcmchecks}
	The following comparisons were made (Section~\ref{sec:comp}); the proofs
	do not use them. For the \(37\) primes
	\(7\le p\le173\) and each pair of a polynomial \(P\) and a Frobenius type,
	the number of \(\Fp\)-isomorphism classes of curves found by the
	enumeration, their mass, and the order of \(\Aut_{\Fp}(C)\) for each of
	the \(325\) orbits of branch divisors are those given by
	Theorems~\ref{thm:vertexcm} and~\ref{thm:vertexcmcells}: the automorphism
	groups have order \(4\) for \(P_c\) and \(6\) for \(P_d\), \(P_e\), and
	the number of curves is the number of classes of
	Theorem~\ref{thm:vertexcm}(ii) minus one when there is a Weil restriction
	(script \texttt{structure\_vs\_census.py}). The assertions of
	Lemma~\ref{lem:orders}, of Lemma~\ref{lem:classno}(a)--(c) and of
	Lemma~\ref{lem:ideals}(b) were verified with PARI/GP for the \(165\)
	primes \(7\le p<1000\), the class groups being certified and
	\(\#\mathcal U(\mathcal O_{\mathbb L})\) being computed from its
	definition by running over the ideal classes of \(\mathbb L\) (scripts
	\texttt{biquad\_G.gp}, \texttt{duals\_ram.gp}); the count of
	Lemma~\ref{lem:ideals}(d) was verified by listing the units of
	\(\mathcal O_{\mathbb L}/2\mathcal O_{\mathbb L}\). A second set of scripts,
	written independently, repeats these verifications for \(547\) primes up
	to \(3989\), computes \(\#\mathcal U(R_{\mathrm{sp}})\) from its definition
	as well, and extends the comparison with the curves to the six primes
	\(179\le p\le199\) (scripts \texttt{rv\_lattice.gp},
	\texttt{rv\_census\_recount.py}).
\end{remark}

\clearpage
\section*{Statements and Declarations}

\small
\noindent\textbf{Competing interests.} The author has no competing interests
to declare that are relevant to the content of this article.

\smallskip
\noindent\textbf{Funding.} No funding was received for conducting this study.

\smallskip
\noindent\textbf{Data availability.} The scripts, logs and data of the
computations of Section~\ref{sec:comp} are provided as supplementary material.

\smallskip
\noindent\textbf{Use of AI tools.} The research questions and the plan of the paper are the author's. The author used a large
language model (Claude, Anthropic) as an assistant: to write and run exploratory
and verification computations; to propose approaches, candidate proofs and
arguments, among which the author chose and which the author examined, corrected,
and accepted or rejected; to search the literature; and to draft and edit the
text. The author read every proof in the article, checked the statements,
computations and citations, and obtained reviews of the manuscript, some of them
generated with other AI systems, which the author directed and evaluated. The
author takes full responsibility for the content of the article. No AI tool is
an author of this work.
\normalsize

\bibliographystyle{amsplain}
\bibliography{references,vertex_pm_refs,vertex_cm_refs}

\end{document}